\pdfoutput=1
\documentclass[11pt]{article}
\usepackage[left=1in,right=1in,top=1in,bottom=1in]{geometry}
\usepackage{times}
\usepackage{expl3}
\usepackage{cite}
\usepackage[table]{xcolor}
\usepackage{multirow}
\usepackage{stackengine} 
\usepackage{hhline}
\usepackage{lipsum}
\usepackage{titlesec}
\usepackage{wrapfig}
\usepackage{enumerate}
\usepackage{epsfig}
\usepackage{amsmath}
\usepackage{tabularx}
\usepackage{array}
\usepackage{booktabs}
\usepackage{enumitem}
\usepackage{bbm}
\usepackage{calc}
\usepackage{graphicx}
\usepackage{amsmath}
\usepackage[title]{appendix}
\usepackage{amssymb}
\usepackage{epstopdf}
\usepackage{boldline}
\usepackage{arydshln}
\usepackage{calligra}
\usepackage{bm}
\usepackage{url}
\usepackage{blindtext}
\usepackage{accents}

\newcommand{\define}{\stackrel{\mbox{\tiny def}}{=}}

\newtheorem{theorem}{Theorem}
\newtheorem{proposition}{Proposition}

\newtheorem{lemma}{Lemma}

\newtheorem{example}{Example}

\usepackage{mathtools}
\usepackage{epstopdf}
\usepackage{balance}
\usepackage{thmtools}
\usepackage{thm-restate}
\usepackage{hyperref}
\usepackage{cleveref}
\usepackage[mathscr]{euscript}

\usepackage[ruled,vlined]{algorithm2e}
\include{pythonlisting}

\newcommand{\ostar}{\mathbin{\mathpalette\make@circled\star}}

\makeatletter
\newcommand{\removelatexerror}{\let\@latex@error\@gobble}
\makeatother
\makeatletter
\newcommand*{\rom}[1]{\expandafter\@slowromancap\romannumeral #1@}
\makeatother

\ExplSyntaxOn
\newcommand\latinabbrev[1]{
  \peek_meaning:NTF. {
    #1\@}%
  { \peek_catcode:NTF a {
      #1.\@ }%
    {#1.\@}}}
\ExplSyntaxOff

\titleclass{\subsubsubsection}{straight}[\subsubsection]

\begin{document}
\vspace{1cm}
\title{Generalized Multiple Operator Integrals for Operators with Finite Dimensions}
\vspace{1.8cm}
\author{Shih-Yu~Chang
\thanks{Shih-Yu Chang is with the Department of Applied Data Science,
San Jose State University, San Jose, CA, U. S. A. (e-mail: {\tt
shihyu.chang@sjsu.edu}). 
           }}

\maketitle

\begin{abstract}
Multiple Operator Integrals (MOIs) have played a foundational role in operator theory and functional calculus, particularly for analyzing Hermitian matrices via spectral decomposition. Conventional MOIs rely on the assumption of self-adjointness, making them analytically tractable for computing Fréchet derivatives, establishing trace formulas, and deriving commutator estimates. However, many problems in mathematics, science, and engineering involve matrices that are fundamentally non-Hermitian, arising in contexts such as numerical discretization of differential operators, signal processing, control systems, and non-Hermitian physics. These cases necessitate a more general framework for operator integration. In this paper, we propose and develop the theory of \emph{Generalized Multiple Operator Integrals (GMOIs)}, which extends MOI techniques to arbitrary matrices, including non-Hermitian and non-normal cases. We unify conventional MOIs under the perspective of the Spectral Mapping Theorem and present a rigorous construction of Generalized Double Operator Integrals (GDOIs) in finite-dimensional setting via Jordan decomposition of any matrix. We establish key algebraic properties, derive norm estimates, and prove continuity and Lipschitz-type perturbation formulas. Finally, we demonstrate how GMOIs can be used to compute matrix function derivatives in non-Hermitian environments, thus significantly broadening the applicability of MOI-based methods in both theoretical and practical domains.
\end{abstract}


\section{Introduction}\label{sec: Introduction}

﻿In the study of functional calculus and operator theory, Multiple Operator Integrals (MOIs) have emerged as effective tools for analyzing the properties of operator functions, mainly in perturbation theory and non-commutative analysis~\cite{chang2024generalizedCDJ,chang2022randomMOI,chang2022randomDTI,chang2022randomPDT}. The conventional method of MOIs is basically grounded from the spectral concept of Hermitian (self-adjoint) matrices, where the spectral decomposition permits the representation of operator functions through integration over spectral projections. In this framework, functions with several argument Hermitian matrices are expressed as iterated integrals involving spectral measures associated with each operator, and kernel functions are used to model the interaction among eigenvalues from input matrices. This conventional MOI method plays an important position in deriving high-order Fréchet derivatives, establishing trace formulae, and facilitating the analysis of commutator estimates. The assumption of Hermitian structure guarantees properly-described spectral projections and facilitates analytical tractability, making it a cornerstone in the theoretical basis of present day perturbation analysis in matrix and operator theory.

﻿While Hermitian matrices play an important function in lots of areas of mathematics and quantum physics due to their real eigenvalues and orthogonal eigenvectors, a wide variety of programs throughout mathematics, science, and engineering involve matrices which are fundamentally non-Hermitian. These non-Hermitian matrices, which do not satisfy the situation $\bm{A} = \bm{A}^{\mathrm{H}}$, where $\mathrm{H}$ is the Hermitian operation, rise up naturally in contexts in which dissipation, asymmetry, or directionality is intrinsic to the system. In the field of numerical linear algebra and differential equations, the discretization of non-self-adjoint differential operators (which includes convection-diffusion operators) ends in non-Hermitian matrices~\cite{trefethen1997numerical}. These operators lack the symmetry necessary to guarantee real spectra, posing unique challenges in spectral analysis and numerical stability.

In engineering fields, non-Hermitian matrices frequently seem in control principle, signal processing, and communication dynamics. Systems with feedback delay, non-reciprocal transmission, or anisotropic fabric residences regularly yield non-Hermitian gadget matrices. In particular, the modeling of open or lossy physical systems in optics and electromagnetics, inclusive of in emph non-Hermitian photonics, calls for complicated-symmetric or even complicated-nonsymmetric matrices to explain energy alternate and mode coupling~\cite{el2018non}. Furthermore, in carrying out physics and biology, non-Hermitian Hamiltonians are used to explain systems with gain and loss or irreversible procedures, together with population dynamics and neural networks~\cite{ashida2020non}. As such, the study and application of non-Hermitian matrices are both theoretically wealthy and almost quintessential in many contemporary scientific and engineering fields.

Conventional formulations of Multiple Operator Integrals (MOIs) should rely on the spectral decomposition of Hermitian matrices, restricting their applicability to cases with non Hermitian matrices. However, many practical and theoretical problems in modern analysis, physics, and engineering involve operators or matrices that are not Hermitian. To overcome this limitation, we proposed a framework of \emph{Generalized MOIs}, which extend the theory to arbitrary matrices—including those that are non-Hermitian or even non-normal. This extension preserves the integral representation structure while accounting for the broader spectral behavior of such matrices. In the finite-dimensional context, we developed a rigorous construction for Generalized Double Operator Integrals (GDOIs) based on matrix decomposition techniques and Jordan canonical forms~\cite{chang2025GDOIMatrix}. Moreover, in infinite-dimensional settings involving operators with continuous spectra, we further generalized the integral framework to accommodate spectral measures beyond the discrete case, offering new tools for functional calculus in non-self-adjoint operator theory~\cite{chang2025GDOICont}. These developments open the door for applying MOI-based techniques in diverse applications, including non-Hermitian perturbation theory, control theory, and quantum systems with complex spectra.

The primary contributions of this paper are as follows. First, we establish that conventional MOIs for Hermitian matrices can be unified under the perspective of the Spectral Mapping Theorem, offering a conceptual bridge to more general formulations. We then introduce the novel framework of Generalized Multiple Operator Integrals (GMOIs), which extends the applicability of MOIs to arbitrary square matrices, including non-Hermitian and non-normal cases. This generalization is accompanied by a thorough investigation of their algebraic and analytic properties. We derive explicit norm estimates for GMOIs, providing essential tools for their boundedness and stability analysis. Furthermore, we develop a perturbation formula and corresponding Lipschitz estimates, enabling controlled sensitivity analysis of operator functions under matrix perturbations. The continuity properties of GMOIs are rigorously characterized, ensuring their compatibility with continuous spectral data. Lastly, we illustrate the effectiveness of our framework by applying the perturbation formula to compute derivatives of matrix functions in terms of GMOIs, thereby extending the scope of matrix calculus in finite-dimensional and non-Hermitian settings.

This paper is organized as follows. In Section~\ref{sec: Conventional MOIs are Special Cases of Spectral Mapping Theorem}, we revisit classical Multiple Operator Integrals (MOIs) and demonstrate that they can be understood as specific instances of the Spectral Mapping Theorem applied to Hermitian matrices. Section~\ref{sec: Generalized Multiple Operator Integrals and Their Algebraic Properties} introduces the framework of Generalized Multiple Operator Integrals (GMOIs) applicable to general, possibly non-Hermitian matrices, and explores their fundamental algebraic structures. Section~\ref{sec: GMOI Norm Estimations} is devoted to the derivation of norm estimates for GMOIs, establishing rigorous bounds that are essential for perturbative analysis. In Section~\ref{sec:Perturbation Formula and Lipschitz Estimation}, we develop a perturbation formula for GMOIs and derive corresponding Lipschitz-type continuity estimates. Section~\ref{sec: GMOI Continuity} examines the continuity properties of GMOIs with respect to parameter matrix perturbations, ensuring the robustness of the framework under small variations. Finally, in Section~\ref{sec: Applications: Differentiation of Matrix Functions}, we demonstrate the utility of the proposed theory by applying the perturbation formula to derive differentiation results for matrix-valued functions.

\section{Conventional MOIs are Special Cases of Spectral Mapping Theorem}\label{sec: Conventional MOIs are Special Cases of Spectral Mapping Theorem}

Let us review conventional DOI definitions. Given a function $\beta: \mathbb{R}^{\zeta+1} \rightarrow \mathbb{C}$, $\zeta+1$ Hermitian matrices (parameter matrices) $\bm{X}_1,\bm{X}_2,\ldots, \bm{X}_{\zeta+1} \in\mathbb{C}^{n \times n}$, and any $\zeta$ matrices (argument matrices) $\bm{Y}_1,\bm{Y}_2,\ldots, \bm{Y}_{\zeta} \in\mathbb{C}^{n \times n}$. From spectral mapping theorem, we have parameter matrices
\begin{eqnarray}\label{eq0-1:  conv DOI def}
\bm{X}_1&=&\sum\limits_{k_1=1}^{K_1}\sum\limits_{i_1=1}^{\alpha^{\mathrm{G}}_{k_1}}\lambda_{k_1}\bm{P}_{k_1,i_1};\nonumber \\
\bm{X}_2&=&\sum\limits_{k_2=1}^{K_2}\sum\limits_{i_2=1}^{\alpha^{\mathrm{G}}_{k_2}}\lambda_{k_2}\bm{P}_{k_2,i_2};\nonumber \\
\vdots&&\nonumber \\
\bm{X}_{\zeta+1}&=&\sum\limits_{k_{\zeta+1}=1}^{K_{\zeta+1}}\sum\limits_{i_{\zeta+1}=1}^{\alpha^{\mathrm{G}}_{k_{\zeta+1}}}\lambda_{k_{\zeta+1}}\bm{P}_{k_{\zeta+1},i_{\zeta+1}},
\end{eqnarray}
where $K_1, K_2,\ldots, K_{\zeta+1}$ are the numbers of distinct eigenvalues of the matrices $\bm{X}_1, \bm{X}_2, \ldots,\bm{X}_{\zeta+1}$,\\
$\alpha^{\mathrm{G}}_{k_1}, \alpha^{\mathrm{G}}_{k_2},\ldots,\alpha^{\mathrm{G}}_{k_{\zeta+1}}$ are the geometry multiplicities of distinct eigenvalues $\lambda_{k_1}, \lambda_{k_2},\ldots,\lambda_{k_{\zeta+1}}$ of the matrices $\bm{X}_1, \bm{X}_2,\ldots,\bm{X}_{\zeta+1}$, and $\bm{P}_{k_1,i_1}, \bm{P}_{k_2,i_2},\ldots,\bm{P}_{k_{\zeta+1},i_{\zeta+1}}$ are the projector matrices corresponding to the $i_j$-th geometric component of the $k_j$-th eigenvalue of the matrices $\bm{X}_j$ , where $j=1,2,\ldots,\zeta+1$. Let $\alpha^{\mathrm{A}}_{k_1}, \alpha^{\mathrm{A}}_{k_2},\ldots,\alpha^{\mathrm{A}}_{k_{\zeta+1}}$ be the algebraic multiplicities of distinct eigenvalues $\lambda_{k_1}, \lambda_{k_2},\ldots,\lambda_{k_{\zeta+1}}$ of the matrices $\bm{X}_1, \bm{X}_2, \ldots,\bm{X}_{\zeta+1}$, we have $\sum\limits_{k_j=1}^{K_j}\alpha^{\mathrm{A}}_{k_j}=n$ for $j=1,2,\ldots,\zeta+1$. 

Similarly, from spectral mapping theorem, we have argument matrices
\begin{eqnarray}\label{eq0-2:  conv DOI def}
\bm{Y}_1&=&\sum\limits_{k'_1=1}^{K'_1}\sum\limits_{i'_1=1}^{\alpha^{\mathrm{G}}_{k'_1}}\lambda_{k'_1}\bm{P}_{k'_1,i'_1}+\sum\limits_{k'_{1}=1}^{K'_{1}}\sum\limits_{i'_{1}=1}^{\alpha^{\mathrm{G}}_{k'_{1}}}\bm{N}_{k'_{1},i'_{1}};\nonumber \\
\bm{Y}_2&=&\sum\limits_{k'_2=1}^{K'_2}\sum\limits_{i'_2=1}^{\alpha^{\mathrm{G}}_{k'_2}}\lambda_{k'_2}\bm{P}_{k'_2,i'_2}+\sum\limits_{k'_{2}=1}^{K'_{2}}\sum\limits_{i'_{2}=1}^{\alpha^{\mathrm{G}}_{k'_{2}}}\bm{N}_{k'_{2},i'_{2}};\nonumber \\
\vdots&&\nonumber \\
\bm{Y}_{\zeta}&=&\sum\limits_{k'_{\zeta}=1}^{K'_{\zeta}}\sum\limits_{i'_{\zeta}=1}^{\alpha^{\mathrm{G}}_{k'_{\zeta}}}\lambda_{k'_{\zeta}}\bm{P}_{k'_{\zeta},i'_{\zeta}}+\sum\limits_{k'_{\zeta}=1}^{K'_{\zeta}}\sum\limits_{i'_{\zeta}=1}^{\alpha^{\mathrm{G}}_{k'_{\zeta}}}\bm{N}_{k'_{\zeta},i'_{\zeta}},
\end{eqnarray}
where $K'_1, K'_2,\ldots, K'_{\zeta}$ are the numbers of distinct eigenvalues of the matrices $\bm{X}_1, \bm{X}_2, \ldots,\bm{X}_{\zeta}$,\\
$\alpha^{\mathrm{G}}_{k'_1}, \alpha^{\mathrm{G}}_{k'_2},\ldots,\alpha^{\mathrm{G}}_{k'_{\zeta}}$ are the geometry multiplicities of distinct eigenvalues $\lambda_{k'_1}, \lambda_{k'_2},\ldots,\lambda_{k'_{\zeta}}$ of the matrices $\bm{Y}_1, \bm{Y}_2,\ldots,\bm{Y}_{\zeta}$, and $\bm{P}_{k'_1,i'_1}, \bm{P}_{k'_2,i'_2},\ldots,\bm{P}_{k'_{\zeta},i'_{\zeta}}$ are the projector matrices corresponding to the $i_j$-th geometric component of the $k_j$-th eigenvalue of the matrices $\bm{Y}_j$ , where $j=1,2,\ldots,\zeta$. Let $\alpha^{\mathrm{A}}_{k'_1}, \alpha^{\mathrm{A}}_{k'_2},\ldots,\alpha^{\mathrm{A}}_{k'_{\zeta}}$ be the algebraic multiplicities of distinct eigenvalues $\lambda_{k'_1}, \lambda_{k'_2},\ldots,\lambda_{k'_{\zeta}}$ of the matrices $\bm{Y}_1, \bm{Y}_2,\ldots,\bm{Y}_{\zeta}$, we have $\sum\limits_{k'_j=1}^{K'_j}\alpha^{\mathrm{A}}_{k'_j}=n$ for $j=1,2,\ldots,\zeta$. 

Let us present Theorem 3 in~\cite{chang2024operatorChar} first. Before presenting this theorem, we review several special ntations related to this Theorem 3 in~\cite{chang2024operatorChar}. Given $r$ positive integers $q_1,q_2,\ldots,q_r$, we define $\varrho_{\kappa}(q_1,\ldots,q_r)$ to be the selection of these $r$ arguments $q_1,\ldots,q_r$ to $\kappa$ arguments, i.e., we have
\begin{eqnarray}
\varrho_{\kappa}(q_1,\ldots,q_r)&=&\{q_{\iota_1},q_{\iota_2},\ldots,q_{\iota_\kappa}\}.
\end{eqnarray}
We use $\mbox{Ind}(\varrho_{\kappa}(q_1,\ldots,q_r))$ to obtain indices of those $\kappa$ positive integers $\{q_{\iota_1},q_{\iota_2},\ldots,q_{\iota_\kappa}\}$, i.e., we have
\begin{eqnarray}
\mbox{Ind}(\varrho_{\kappa}(q_1,\ldots,q_r))&=&\{\iota_1,\iota_2,\ldots,\iota_\kappa\}.
\end{eqnarray}
We use $\varrho_{\kappa}(q_1,\ldots,q_r)=1$ to represent $q_{\iota_1}=1,q_{\iota_2}=1,\ldots,q_{\iota_\kappa}=1$. We also use \\
$m_{k_{\mbox{Ind}(\varrho_{\kappa}(q_1,\ldots,q_r))},i_{\mbox{Ind}(\varrho_{\kappa}(q_1,\ldots,q_r))}}-1$ to represent $m_{k_{\iota_1},i_{\iota_1}}-1,m_{k_{\iota_2},i_{\iota_2}}-1,\ldots,m_{k_{\iota_\kappa},i_{\iota_\kappa}}-1$, where $m_{k_{\iota_j},i_{\iota_j}}$ is the order for the nilpotent matrix $\bm{N}_{k_{\iota_j},i_{\iota_j}}$, i.e., $\bm{N}^\ell_{k_{\iota_j},i_{\iota_j}} = \bm{0}$, for $\ell \geq m_{k_{\iota_j},i_{\iota_j}}$ and $j=1,2,\ldots,\kappa$.

Then, Theorem 3 in~\cite{chang2024operatorChar} is given below.
\begin{theorem}\label{thm: Spectral Mapping Theorem for r Variables}
Given an analytic function $f(z_1,z_2,\ldots,z_r)$ within the domain for $|z_l| < R_l$, and the matrix $\bm{X}_l$ with the dimension $m$ and $K_l$ distinct eigenvalues $\lambda_{k_l}$ for $k_l=1,2,\ldots,K_l$ such that
\begin{eqnarray}\label{eq1-1: thm: Spectral Mapping Theorem for r Variables}
\bm{X}_l&=&\sum\limits_{k_l=1}^{K_l}\sum\limits_{i_l=1}^{\alpha_{k_l}^{\mathrm{G}}} \lambda_{k_l} \bm{P}_{k_l,i_l}+
\sum\limits_{k_l=1}^{K_l}\sum\limits_{i_l=1}^{\alpha_{k_l}^{\mathrm{G}}} \bm{N}_{k_l,i_l},
\end{eqnarray}
where $\left\vert\lambda_{k_l}\right\vert<R_l$ for $l=1,2,\ldots,r$.

Then, we have
\begin{eqnarray}\label{eq2: thm: Spectral Mapping Theorem for kappa Variables}
\lefteqn{f(\bm{X}_1,\ldots,\bm{X}_r)=}\nonumber \\
&& \sum\limits_{k_1=\ldots=k_r=1}^{K_1,\ldots,K_r} \sum\limits_{i_1=\ldots=i_r=1}^{\alpha_{k_1}^{(\mathrm{G})},\ldots,\alpha_{k_r}^{(\mathrm{G})}}
f(\lambda_{k_1},\ldots,\lambda_{k_r})\bm{P}_{k_1,i_1}\ldots\bm{P}_{k_r,i_r}\nonumber \\
&&+\sum\limits_{k_1=\ldots=k_r=1}^{K_1,\ldots,K_r} \sum\limits_{i_1=\ldots=i_r=1}^{\alpha_{k_1}^{(\mathrm{G})},\ldots,\alpha_{k_r}^{(\mathrm{G})}}\sum\limits_{\kappa=1}^{r-1}\sum\limits_{\varrho_\kappa(q_1,\ldots,q_r)}\Bigg(\sum\limits_{\varrho_{\kappa}(q_1,\ldots,q_r)=1}^{m_{k_{\mbox{Ind}(\varrho_{\kappa}(q_1,\ldots,q_r))},i_{\mbox{Ind}(\varrho_{\kappa}(q_1,\ldots,q_r))}}-1}\nonumber \\
&&~~~~~ \frac{f^{\varrho_{\kappa}(q_1,\ldots,q_r)}(\lambda_{k_1},\ldots,\lambda_{k_r})}{q_{\iota_1}!q_{\iota_2}!\ldots q_{\iota_\kappa}!}\times \prod\limits_{\substack{\varsigma =\mbox{Ind}(\varrho_{\kappa}(q_1,\ldots,q_r)), \bm{Y}=\bm{N}^{q_\varsigma}_{k_\varsigma,i_\varsigma} \\ \varsigma \neq \mbox{Ind}(\varrho_{\kappa}(q_1,\ldots,q_r)), \bm{Y}=\bm{P}_{k_\varsigma,i_\varsigma} }
}^{r} \bm{Y}\Bigg) 
\nonumber \\
&&+\sum\limits_{k_1=\ldots=k_r=1}^{K_1,\ldots,K_r} \sum\limits_{i_1=\ldots=i_r=1}^{\alpha_{k_1}^{(\mathrm{G})},\ldots,\alpha_{k_r}^{(\mathrm{G})}} \sum\limits_{q_1=\ldots=q_r=1}^{m_{k_1,i_1}-1,\ldots,m_{k_r,i_r}-1}
\frac{f^{(q_1,\ldots,q_r)}(\lambda_{k_1},\ldots,\lambda_{k_r})}{q_1!\cdots q_r!}\bm{N}^{q_1}_{k_1,i_1}\ldots\bm{N}^{q_r}_{k_r,i_r}
\end{eqnarray}
where we have
\begin{itemize}
\item $\sum\limits_{\varrho_{\kappa}(q_1,\ldots,q_r)}$ is the summation running over all selection of $\varrho_{\kappa}(q_1,\ldots,q_r)$ given $\kappa$;
\item $f^{\varrho_{\kappa}(q_1,\ldots,q_r)}(\lambda_1,\ldots,\lambda_r)$ represents the partial derivatives with respect to variables with indices $\iota_1,\iota_2,\ldots,\iota_\kappa$ and the orders of derivatives given by $q_{\iota_1},q_{\iota_2},\ldots,q_{\iota_\kappa}$. 
\end{itemize}
\end{theorem}

The conventional MOI is a matrix, denoted by $T_{\beta}^{\bm{X}_1,\ldots,\bm{X}_{\zeta+1}}(\bm{Y}_1,\ldots,\bm{Y}_{\zeta})$, which can be expressed as~\cite{skripka2019multilinear}:
\begin{eqnarray}\label{eq1:  conv MOI def}
T_{\beta}^{\bm{X}_1,\ldots,\bm{X}_{\zeta+1}}(\bm{Y}_1,\ldots,\bm{Y}_{\zeta})&\define&\sum\limits_{k_1=1}^{K_1}\ldots \sum\limits_{k_{\zeta+1}=1}^{K_{\zeta+1}}\sum\limits_{i_1=1}^{\alpha_{k_1}^{(\mathrm{G})}} \ldots \sum\limits_{i_{\zeta+1}=1}^{\alpha_{k_{\zeta+1}}^{(\mathrm{G})}}\beta(\lambda_{k_1},\ldots, \lambda_{k_{\zeta+1}})
\nonumber \nonumber \\
&& \times \bm{P}_{k_1,i_1}\bm{Y}_1 \bm{P}_{k_2,i_2} \ldots \bm{P}_{k_\zeta,i_\zeta}\bm{Y}_\zeta \bm{P}_{k_{\zeta+1},i_{\zeta+1}}
\end{eqnarray}
From the decomposition of the matrix $\bm{Y}_j$ given by Eq.~\eqref{eq0-2:  conv DOI def}, Eq.~\eqref{eq1:  conv MOI def} can further be expressed as
\begin{eqnarray}\label{eq2:  conv MOI def}
\lefteqn{T_{\beta}^{\bm{X}_1,\ldots,\bm{X}_{\zeta+1}}(\bm{Y}_1,\ldots,\bm{Y}_{\zeta})=\sum\limits_{k_1=1}^{K_1}\ldots \sum\limits_{k_{\zeta+1}=1}^{K_{\zeta+1}}\sum\limits_{i_1=1}^{\alpha_{k_1}^{(\mathrm{G})}} \ldots \sum\limits_{i_{\zeta+1}=1}^{\alpha_{k_{\zeta+1}}^{(\mathrm{G})}}\beta(\lambda_{k_1},\ldots, \lambda_{k_{\zeta+1}})}\nonumber \\
&& \times \bm{P}_{k_1,i_1}\left(\sum\limits_{k'_1=1}^{K'_1}\sum\limits_{i'_1=1}^{\alpha^{\mathrm{G}}_{k'_1}}\lambda_{k'_1}\bm{P}_{k'_1,i'_1}+\sum\limits_{k'_1=1}^{K'_1}\sum\limits_{i'_1=1}^{\alpha^{\mathrm{G}}_{k'_1}}\bm{N}_{k'_1,i'_1}\right)\bm{P}_{k_2,i_2} \ldots \nonumber \\
&& \times \bm{P}_{k_\zeta,i_\zeta}\left(\sum\limits_{k'_\zeta=1}^{K'_\zeta}\sum\limits_{i'_\zeta=1}^{\alpha^{\mathrm{G}}_{k'_\zeta}}\lambda_{k'_\zeta}\bm{P}_{k'_\zeta,i'_\zeta}+\sum\limits_{k'_\zeta=1}^{K'_\zeta}\sum\limits_{i'_\zeta=1}^{\alpha^{\mathrm{G}}_{k'_\zeta}}\bm{N}_{k'_\zeta,i'_\zeta}\right)\bm{P}_{k_{\zeta+1},i_{\zeta+1}}\nonumber \\
&=& \sum\limits_{k_1=\ldots=k_{\zeta+1}=1}^{K_1,\ldots,K_{\zeta+1}} \sum\limits_{i_1=\ldots=i_{\zeta+1}=1}^{\alpha_{k_1}^{(\mathrm{G})},\ldots,\alpha_{k_{\zeta+1}}^{(\mathrm{G})}}  \sum\limits_{k'_1=\ldots=k'_{\zeta}=1}^{K'_1,\ldots,K'_{\zeta}} \sum\limits_{i'_1=\ldots=i'_{\zeta}=1}^{\alpha_{k'_1}^{(\mathrm{G})},\ldots,\alpha_{k'_{\zeta}}^{(\mathrm{G})}}
\beta(\lambda_{k_1},\ldots,\lambda_{k_{\zeta+1}})\lambda_{k'_1}\cdots\lambda_{k'_{\zeta}}\nonumber \\
&& \times \bm{P}_{k_1,i_1}\bm{P}_{k'_1,i'_1}\cdots \bm{P}_{k'_\zeta,i'_\zeta}\bm{P}_{k_{\zeta+1},i_{\zeta+1}}\nonumber \\
&&+\sum\limits_{k_1=\ldots=k_{\zeta+1}=1}^{K_1,\ldots,K_{\zeta+1}} \sum\limits_{i_1=\ldots=i_{\zeta+1}=1}^{\alpha_{k_1}^{(\mathrm{G})},\ldots,\alpha_{k_{\zeta+1}}^{(\mathrm{G})}}  \sum\limits_{k'_1=\ldots=k'_{\zeta}=1}^{K'_1,\ldots,K'_{\zeta}} \sum\limits_{i'_1=\ldots=i'_{\zeta}=1}^{\alpha_{k'_1}^{(\mathrm{G})},\ldots,\alpha_{k'_{\zeta}}^{(\mathrm{G})}}\beta(\lambda_{k_1},\ldots,\lambda_{k_{\zeta+1}})\sum\limits_{\kappa=1}^{\zeta-1}\sum\limits_{\varrho_\kappa(1,\ldots,\zeta)}\nonumber \\
&& \Bigg(\prod\limits_{\substack{\varsigma =\mbox{Ind}(\varrho_{\kappa}(1,\ldots,\zeta)), \bm{Y}=\bm{N}_{k'_\varsigma,i'_\varsigma} \\ \varsigma \neq \mbox{Ind}(\varrho_{\kappa}(1,\ldots,\zeta)), \bm{Y}=\lambda_{k'_\varsigma}\bm{P}_{k'_\varsigma,i'_\varsigma}}}^{\zeta} \bm{Y}\Bigg) 
\nonumber \\
&&+ \sum\limits_{k_1=\ldots=k_{\zeta+1}=1}^{K_1,\ldots,K_{\zeta+1}} \sum\limits_{i_1=\ldots=i_{\zeta+1}=1}^{\alpha_{k_1}^{(\mathrm{G})},\ldots,\alpha_{k_{\zeta+1}}^{(\mathrm{G})}}  \sum\limits_{k'_1=\ldots=k'_{\zeta}=1}^{K'_1,\ldots,K'_{\zeta}} \sum\limits_{i'_1=\ldots=i'_{\zeta}=1}^{\alpha_{k'_1}^{(\mathrm{G})},\ldots,\alpha_{k'_{\zeta}}^{(\mathrm{G})}}\beta(\lambda_{k_1},\ldots,\lambda_{k_{\zeta+1}})\nonumber \\
&& \times \bm{P}_{k_1,i_1}\bm{N}_{k'_1,i'_1}\cdots \bm{N}_{k'_\zeta,i'_\zeta}\bm{P}_{k_{\zeta+1},i_{\zeta+1}},
\end{eqnarray}
where $\varrho_{\kappa}(1,\ldots,\zeta)$ is the selection of $\kappa$ indices from the indices set $\{1,2,\ldots,\zeta\}$, and $\mbox{Ind}(\varrho_{\kappa}(1,\ldots,\zeta))$ represents those $\kappa$ indices being selected. 

If we set the function $f$ in Theorem~\ref{thm: Spectral Mapping Theorem for r Variables} by
\begin{eqnarray}
f(z_1,z_2, \ldots, z_{2\zeta+1}) &=& \beta(z_1, z_3, \ldots,z_{2\zeta+1})z_2 z_4 \ldots z_{2\zeta}, 
\end{eqnarray}
and set 
\begin{eqnarray}
z_1&=&\lambda_{k_1},~~z_3~=~\lambda_{k_2},~~z_5~=~\lambda_{k_3}, ~~\cdots,~~z_{2\zeta+1}~=~\lambda_{k_{\zeta+1}}, \nonumber \\  
z_2&=&\lambda_{k'_1},~~z_4~=~\lambda_{k'_2},~~z_6~=~\lambda_{k'_3}, ~~\cdots,~~z_{2\zeta}~=~\lambda_{k'_\zeta},
\end{eqnarray}
then, by applying Theorem~\ref{thm: Spectral Mapping Theorem for r Variables}, we will obtain Eq.~\eqref{eq2:  conv MOI def}.

Following Example~\ref{exp:TOI is a special case} shows that the conventional MOI is a special case of spectrum mapping theorem when $\zeta=2$.
\begin{example}\label{exp:TOI is a special case}
Let us consider a triple operator integral $T_\beta^{\bm{X}_1,\bm{X}_2,\bm{X}_3}(\bm{Y}_1, \bm{Y}_2)$ with Hermitian matrices $\bm{X}_1,  \bm{X}_2$ and $\bm{X}_3$. Then, we can express $T_\beta^{\bm{X}_1,\bm{X}_2,\bm{X}_3}(\bm{Y}_1, \bm{Y}_2)$ as 
\begin{eqnarray}\label{eq1:exp:TOI is a special case}
\lefteqn{T_{\beta}^{\bm{X}_1,\bm{X}_2,\bm{X}_3}(\bm{Y}_1, \bm{Y}_2)}\nonumber\\
&\define&\sum\limits_{k_1=1}^{K_1}\sum\limits_{k_2=1}^{K_2}\sum\limits_{k_3=1}^{K_3}\sum\limits_{i_1=1}^{\alpha_{k_1}^{(\mathrm{G})}}\sum\limits_{i_2=1}^{\alpha_{k_2}^{(\mathrm{G})}}\sum\limits_{i_3=1}^{\alpha_{k_3}^{(\mathrm{G})}}\beta(\lambda_{k_1}, \lambda_{k_2}, \lambda_{k_3})\bm{P}_{k_1,i_1}\bm{Y}_1\bm{P}_{k_2,i_2}\bm{Y}_2\bm{P}_{k_3,i_3} \nonumber \\
&=&\sum\limits_{k_1=1}^{K_1}\sum\limits_{k_2=1}^{K_2}\sum\limits_{k_3=1}^{K_3}\sum\limits_{i_1=1}^{\alpha_{k_1}^{(\mathrm{G})}}\sum\limits_{i_2=1}^{\alpha_{k_2}^{(\mathrm{G})}}\sum\limits_{i_3=1}^{\alpha_{k_3}^{(\mathrm{G})}}\beta(\lambda_{k_1}, \lambda_{k_2}, \lambda_{k_3})\bm{P}_{k_1,i_1}\left(\sum\limits_{k'_1=1}^{K'_1}\sum\limits_{i'_1=1}^{\alpha^{\mathrm{G}}_{k'_1}}\lambda_{k'_1}\bm{P}_{k'_1,i'_1}+\sum\limits_{k'_{1}=1}^{K'_{1}}\sum\limits_{i'_{1}=1}^{\alpha^{\mathrm{G}}_{k'_{1}}}\bm{N}_{k'_{1},i'_{1}}\right)\nonumber \\
&& \times \bm{P}_{k_2,i_2}\left(\sum\limits_{k'_2=1}^{K'_2}\sum\limits_{i'_2=1}^{\alpha^{\mathrm{G}}_{k'_2}}\lambda_{k'_2}\bm{P}_{k'_2,i'_2}+\sum\limits_{k'_{2}=1}^{K'_{2}}\sum\limits_{i'_{2}=1}^{\alpha^{\mathrm{G}}_{k'_{2}}}\bm{N}_{k'_{2},i'_{2}}\right)\bm{P}_{k_3,i_3} \nonumber \\
&=&\sum\limits_{k_1=1}^{K_1}\sum\limits_{k_2=1}^{K_2}\sum\limits_{k_3=1}^{K_3}\sum\limits_{i_1=1}^{\alpha_{k_1}^{(\mathrm{G})}}\sum\limits_{i_2=1}^{\alpha_{k_2}^{(\mathrm{G})}}\sum\limits_{i_3=1}^{\alpha_{k_3}^{(\mathrm{G})}}\beta(\lambda_{k_1}, \lambda_{k_2}, \lambda_{k_3})\bm{P}_{k_1,i_1}\nonumber \\
&& \times \left(\sum\limits_{k'_1=1}^{K'_1}\sum\limits_{i'_1=1}^{\alpha^{\mathrm{G}}_{k'_1}}\lambda_{k'_1}\bm{P}_{k'_1,i'_1}\right)\bm{P}_{k_2,i_2}\left(\sum\limits_{k'_2=1}^{K'_2}\sum\limits_{i'_2=1}^{\alpha^{\mathrm{G}}_{k'_2}}\lambda_{k'_2}\bm{P}_{k'_2,i'_2}\right)\bm{P}_{k_3,i_3}  \nonumber \\
&&+\sum\limits_{k_1=1}^{K_1}\sum\limits_{k_2=1}^{K_2}\sum\limits_{k_3=1}^{K_3}\sum\limits_{i_1=1}^{\alpha_{k_1}^{(\mathrm{G})}}\sum\limits_{i_2=1}^{\alpha_{k_2}^{(\mathrm{G})}}\sum\limits_{i_3=1}^{\alpha_{k_3}^{(\mathrm{G})}}\beta(\lambda_{k_1}, \lambda_{k_2}, \lambda_{k_3})\bm{P}_{k_1,i_1}\nonumber \\
&& \times \left(\sum\limits_{k'_1=1}^{K'_1}\sum\limits_{i'_1=1}^{\alpha^{\mathrm{G}}_{k'_1}}\lambda_{k'_1}\bm{P}_{k'_1,i'_1}\right)\bm{P}_{k_2,i_2}\left(\sum\limits_{k'_2=1}^{K'_2}\sum\limits_{i'_2=1}^{\alpha^{\mathrm{G}}_{k'_2}}\bm{N}_{k'_2,i'_2}\right)\bm{P}_{k_3,i_3}  \nonumber \\
&& +\sum\limits_{k_1=1}^{K_1}\sum\limits_{k_2=1}^{K_2}\sum\limits_{k_3=1}^{K_3}\sum\limits_{i_1=1}^{\alpha_{k_1}^{(\mathrm{G})}}\sum\limits_{i_2=1}^{\alpha_{k_2}^{(\mathrm{G})}}\sum\limits_{i_3=1}^{\alpha_{k_3}^{(\mathrm{G})}}\beta(\lambda_{k_1}, \lambda_{k_2}, \lambda_{k_3})\bm{P}_{k_1,i_1}\nonumber \\
&& \times \left(\sum\limits_{k'_1=1}^{K'_1}\sum\limits_{i'_1=1}^{\alpha^{\mathrm{G}}_{k'_1}}\bm{N}_{k'_1,i'_1}\right)\bm{P}_{k_2,i_2}\left(\sum\limits_{k'_2=1}^{K'_2}\sum\limits_{i'_2=1}^{\alpha^{\mathrm{G}}_{k'_2}}\lambda_{k'_2}\bm{P}_{k'_2,i'_2}\right)\bm{P}_{k_3,i_3}  \nonumber \\
&&+\sum\limits_{k_1=1}^{K_1}\sum\limits_{k_2=1}^{K_2}\sum\limits_{k_3=1}^{K_3}\sum\limits_{i_1=1}^{\alpha_{k_1}^{(\mathrm{G})}}\sum\limits_{i_2=1}^{\alpha_{k_2}^{(\mathrm{G})}}\sum\limits_{i_3=1}^{\alpha_{k_3}^{(\mathrm{G})}}\beta(\lambda_{k_1}, \lambda_{k_2}, \lambda_{k_3})\bm{P}_{k_1,i_1}\nonumber \\
&& \times \left(\sum\limits_{k'_1=1}^{K'_1}\sum\limits_{i'_1=1}^{\alpha^{\mathrm{G}}_{k'_1}}\bm{N}_{k'_1,i'_1}\right)\bm{P}_{k_2,i_2}\left(\sum\limits_{k'_2=1}^{K'_2}\sum\limits_{i'_2=1}^{\alpha^{\mathrm{G}}_{k'_2}}\bm{N}_{k'_2,i'_2}\right)\bm{P}_{k_3,i_3} 
\end{eqnarray}

If we set the function $f$ in Theorem~\ref{thm: Spectral Mapping Theorem for r Variables} by
\begin{eqnarray}
f(z_1,z_2,z_3,z_4, z_5) &=& \beta(z_1, z_3, \,z_5)z_2 z_4, 
\end{eqnarray}
and set 
\begin{eqnarray}
z_1&=&\lambda_{k_1},~~z_3~=~\lambda_{k_2},~~z_5~=~\lambda_{k_3}, \nonumber \\  
z_2&=&\lambda_{k'_1},~~z_4~=~\lambda_{k'_2},
\end{eqnarray}
then, by applying Theorem~\ref{thm: Spectral Mapping Theorem for r Variables}, we have
\begin{eqnarray}
\frac{\partial f(z_1,z_2,z_3,z_4, z_5)}{\partial z_2} &=&  \beta(z_1, z_3, \,z_5)z_4;\nonumber \\
\frac{\partial f(z_1,z_2,z_3,z_4, z_5)}{\partial z_4} &=&  \beta(z_1, z_3, \,z_5)z_2;\nonumber \\
\frac{\partial f(z_1,z_2,z_3,z_4, z_5)}{\partial z_2 \partial z_4} &=&  \beta(z_1, z_3, \,z_5). 
\end{eqnarray}
Therefore, we obtain:
\begin{eqnarray}\label{eq2:exp:TOI is a special case}
f(\bm{X}_1,\bm{Y}_1,\bm{X}_2,\bm{Y}_2,\bm{X}_3)
&=&\sum\limits_{k_1=1}^{K_1}\sum\limits_{k_2=1}^{K_2}\sum\limits_{k_3=1}^{K_3}\sum\limits_{i_1=1}^{\alpha_{k_1}^{(\mathrm{G})}}\sum\limits_{i_2=1}^{\alpha_{k_2}^{(\mathrm{G})}}\sum\limits_{i_3=1}^{\alpha_{k_3}^{(\mathrm{G})}}\beta(\lambda_{k_1}, \lambda_{k_2}, \lambda_{k_3})\lambda_{k'_1}\lambda_{k'_2}\bm{P}_{k_1,i_1}\nonumber \\
&& \times \left(\sum\limits_{k'_1=1}^{K'_1}\sum\limits_{i'_1=1}^{\alpha^{\mathrm{G}}_{k'_1}}\bm{P}_{k'_1,i'_1}\right)\bm{P}_{k_2,i_2}\left(\sum\limits_{k'_2=1}^{K'_2}\sum\limits_{i'_2=1}^{\alpha^{\mathrm{G}}_{k'_2}}\bm{P}_{k'_2,i'_2}\right)\bm{P}_{k_3,i_3}  \nonumber \\
&&+\sum\limits_{k_1=1}^{K_1}\sum\limits_{k_2=1}^{K_2}\sum\limits_{k_3=1}^{K_3}\sum\limits_{i_1=1}^{\alpha_{k_1}^{(\mathrm{G})}}\sum\limits_{i_2=1}^{\alpha_{k_2}^{(\mathrm{G})}}\sum\limits_{i_3=1}^{\alpha_{k_3}^{(\mathrm{G})}}\beta(\lambda_{k_1}, \lambda_{k_2}, \lambda_{k_3})\lambda_{k'_1}\bm{P}_{k_1,i_1}\nonumber \\
&& \times \left(\sum\limits_{k'_1=1}^{K'_1}\sum\limits_{i'_1=1}^{\alpha^{\mathrm{G}}_{k'_1}}\bm{P}_{k'_1,i'_1}\right)\bm{P}_{k_2,i_2}\left(\sum\limits_{k'_2=1}^{K'_2}\sum\limits_{i'_2=1}^{\alpha^{\mathrm{G}}_{k'_2}}\bm{N}_{k'_2,i'_2}\right)\bm{P}_{k_3,i_3}  \nonumber \\
&& +\sum\limits_{k_1=1}^{K_1}\sum\limits_{k_2=1}^{K_2}\sum\limits_{k_3=1}^{K_3}\sum\limits_{i_1=1}^{\alpha_{k_1}^{(\mathrm{G})}}\sum\limits_{i_2=1}^{\alpha_{k_2}^{(\mathrm{G})}}\sum\limits_{i_3=1}^{\alpha_{k_3}^{(\mathrm{G})}}\beta(\lambda_{k_1}, \lambda_{k_2}, \lambda_{k_3})\lambda_{k'_2}\bm{P}_{k_1,i_1}\nonumber \\
&& \times \left(\sum\limits_{k'_1=1}^{K'_1}\sum\limits_{i'_1=1}^{\alpha^{\mathrm{G}}_{k'_1}}\bm{N}_{k'_1,i'_1}\right)\bm{P}_{k_2,i_2}\left(\sum\limits_{k'_2=1}^{K'_2}\sum\limits_{i'_2=1}^{\alpha^{\mathrm{G}}_{k'_2}}\bm{P}_{k'_2,i'_2}\right)\bm{P}_{k_3,i_3}  \nonumber \\
&&+\sum\limits_{k_1=1}^{K_1}\sum\limits_{k_2=1}^{K_2}\sum\limits_{k_3=1}^{K_3}\sum\limits_{i_1=1}^{\alpha_{k_1}^{(\mathrm{G})}}\sum\limits_{i_2=1}^{\alpha_{k_2}^{(\mathrm{G})}}\sum\limits_{i_3=1}^{\alpha_{k_3}^{(\mathrm{G})}}\beta(\lambda_{k_1}, \lambda_{k_2}, \lambda_{k_3})\bm{P}_{k_1,i_1}\nonumber \\
&& \times \left(\sum\limits_{k'_1=1}^{K'_1}\sum\limits_{i'_1=1}^{\alpha^{\mathrm{G}}_{k'_1}}\bm{N}_{k'_1,i'_1}\right)\bm{P}_{k_2,i_2}\left(\sum\limits_{k'_2=1}^{K'_2}\sum\limits_{i'_2=1}^{\alpha^{\mathrm{G}}_{k'_2}}\bm{N}_{k'_2,i'_2}\right)\bm{P}_{k_3,i_3}. 
\end{eqnarray}
By comparing Eq.~\eqref{eq1:exp:TOI is a special case} and Eq.~\eqref{eq2:exp:TOI is a special case}, we have
\begin{eqnarray}
f(\bm{X}_1,\bm{Y}_1,\bm{X}_2,\bm{Y}_2,\bm{X}_3)&=&T_{\beta}^{\bm{X}_1,\bm{X}_2,\bm{X}_3}(\bm{Y}_1, \bm{Y}_2).
\end{eqnarray}
\end{example}

\section{Generalized Multiple Operator Integrals and Their Algebraic Properties}\label{sec: Generalized Multiple Operator Integrals and Their Algebraic Properties}

We begin by defining the GMOIs in Section~\ref{sec: Generalized Multiple Operator Integrals}, followed by an exploration of the algebraic properties of MOIs in Section~\ref{sec: GMOI Algebraic Properties}.

\subsection{Generalized Multiple Operator Integrals}\label{sec: Generalized Multiple Operator Integrals}

According to Theorem~\ref{thm: Spectral Mapping Theorem for r Variables}, GMOIs can be defined as follows.
\begin{eqnarray}\label{eq1:GMOI def}
\lefteqn{T_{\beta}^{\bm{X}_1,\ldots,\bm{X}_{\zeta+1}}(\bm{Y}_1,\ldots,\bm{Y}_{\zeta})}\nonumber \\
&\define& \sum\limits_{k_1=\ldots=k_{\zeta+1}=1}^{K_1,\ldots,K_{\zeta+1}} \sum\limits_{i_1=\ldots=i_{\zeta+1}=1}^{\alpha_{k_1}^{(\mathrm{G})},\ldots,\alpha_{k_{\zeta+1}}^{(\mathrm{G})}}
\beta(\lambda_{k_1},\ldots,\lambda_{k_{\zeta+1}})\bm{P}_{k_1,i_1}\bm{Y}_1\bm{P}_{k_2,i_2}\bm{Y}_2\ldots\bm{Y}_{\zeta}\bm{P}_{k_{\zeta+1},i_{\zeta+1}}\nonumber \\
&&+\sum\limits_{k_1=\ldots=k_{\zeta+1}=1}^{K_1,\ldots,K_{\zeta+1}} \sum\limits_{i_1=\ldots=i_{\zeta+1}=1}^{\alpha_{k_1}^{(\mathrm{G})},\ldots,\alpha_{k_{\zeta+1}}^{(\mathrm{G})}}\sum\limits_{\kappa=1}^{\zeta}\sum\limits_{\varrho_\kappa(q_1,\ldots,q_{\zeta+1})}\Bigg(\sum\limits_{\varrho_{\kappa}(q_1,\ldots,q_{\zeta+1})=1}^{m_{k_{\mbox{Ind}(\varrho_{\kappa}(q_1,\ldots,q_{\zeta+1}))},i_{\mbox{Ind}(\varrho_{\kappa}(q_1,\ldots,q_{\zeta+1}))}}-1}\nonumber \\
&&~~~~~ \frac{\beta^{\varrho_{\kappa}(q_1,\ldots,q_{\zeta+1})}(\lambda_{k_1},\ldots,\lambda_{k_{\zeta+1}})}{q_{\iota_1}!q_{\iota_2}!\ldots q_{\iota_\kappa}!}\times \prod\limits_{\substack{\varsigma =\mbox{Ind}(\varrho_{\kappa}(q_1,\ldots,q_r)), \bm{Z}_\varsigma=\bm{N}^{q_\varsigma}_{k_\varsigma,i_\varsigma}\bm{Y}_{\varsigma} \\ \varsigma \neq \mbox{Ind}(\varrho_{\kappa}(q_1,\ldots,q_r)), \bm{Z}_\varsigma=\bm{P}_{k_\varsigma,i_\varsigma}\bm{Y}_{\varsigma}}
}^{\zeta+1} \bm{Z}_\varsigma\Bigg) 
\nonumber \\
&&+\sum\limits_{k_1=\ldots=k_{\zeta+1}=1}^{K_1,\ldots,K_{\zeta+1}} \sum\limits_{i_1=\ldots=i_{\zeta+1}=1}^{\alpha_{k_1}^{(\mathrm{G})},\ldots,\alpha_{k_{\zeta+1}}^{(\mathrm{G})}} \sum\limits_{q_1=\ldots=q_{\zeta+1}=1}^{m_{k_1,i_1}-1,\ldots,m_{k_{\zeta+1},i_{\zeta+1}}-1}
\frac{\beta^{(q_1,\ldots,q_{\zeta+1})}(\lambda_{k_1},\ldots,\lambda_{k_{\zeta+1}})}{q_1!\cdots q_{\zeta+1}!}\nonumber \\
&&
~~\times \bm{N}^{q_1}_{k_1,i_1}\bm{Y}_1 \bm{N}^{q_2}_{k_2,i_2}\bm{Y}_2\ldots \bm{Y}_{\zeta}\bm{N}^{q_{\zeta+1}}_{k_{\zeta+1},i_{\zeta+1}},
\end{eqnarray}
where $\bm{Y}_{\zeta+1}$ is set as identity matrix.

\subsection{GMOI Algebraic Properties}\label{sec: GMOI Algebraic Properties}

In this section, we will establish the algebraic properties of the operator $T_{\beta}^{\bm{X}_1,\ldots,\bm{X}_{\zeta+1}}(\bm{Y}_1,\ldots,\bm{Y}_{\zeta})$ defined by Eq.~\eqref{eq1:GMOI def}. 

If the matrices $\bm{X}_p$ for $p=1,2,\ldots,\zeta+1$ are decomposed as:
\begin{eqnarray}\label{eq:X1 decomp p and n parts}
\bm{X}_p&=&\sum\limits_{k_p=1}^{K_p}\sum\limits_{i_p=1}^{\alpha_{k_p}^{\mathrm{G}}} \lambda_{k_p} \bm{P}_{k_p,i_p}+
\sum\limits_{k_p=1}^{K_p}\sum\limits_{i_p=1}^{\alpha_{k_p}^{\mathrm{G}}} \bm{N}_{k_p,i_p}\nonumber \\
&\define&\bm{X}_{p,\bm{P}}+\bm{X}_{p,\bm{N}},
\end{eqnarray}
then, from the definition of $T_{\beta}^{\bm{X}_1,\ldots,\bm{X}_{\zeta+1}}(\bm{Y}_1,\ldots,\bm{Y}_{\zeta})$ provided in Eq.~\eqref{eq1:GMOI def}, we have the following decomposition proposition with respect to parameters matrices $\bm{X}_p$ immediately.
\begin{proposition}\label{prop:GMOI decomp by parameters X P X N}
Given matrices $\bm{X}_p$, which are decomposed as Eq.~\eqref{eq:X1 decomp p and n parts}, respectively, then, we have
\begin{eqnarray}
T_{\beta}^{\bm{X}_1,\ldots,\bm{X}_{\zeta+1}}(\bm{Y}_1,\ldots,\bm{Y}_{\zeta})&=&\sum\limits_{i=1}^{2^{\zeta+1}}T_{\beta}^{[\bm{X}]_{\xi(i)}}(\bm{Y}_1,\ldots,\bm{Y}_{\zeta}),
\end{eqnarray}
where $[\bm{X}]_{\xi(i)}$ is an array with $\zeta+1$ entries such that its $j$-th entry is expressed by
\begin{eqnarray}
([\bm{X}]_{\xi(i)})_j&=& \begin{cases}
     \bm{X}_{j,P}, & \text{if} \left\lfloor \frac{i-1}{2^{j-1}} \right\rfloor \equiv 0 \pmod{2}; \\
     \bm{X}_{j,N}, & \text{otherwise}.
   \end{cases}
\end{eqnarray}
\end{proposition}

For example, if we consider $\zeta=2$, we have GTOI as below:
\begin{eqnarray}\label{eq1:  GTOI def}
\lefteqn{T_{\beta}^{\bm{X}_1,\bm{X}_2,\bm{X}_3}(\bm{Y}_1, \bm{Y}_2)\define}\nonumber\\
&&\sum\limits_{k_1=1}^{K_1}\sum\limits_{k_2=1}^{K_2}\sum\limits_{k_3=1}^{K_3}\sum\limits_{i_1=1}^{\alpha_{k_1}^{(\mathrm{G})}}\sum\limits_{i_2=1}^{\alpha_{k_2}^{(\mathrm{G})}}\sum\limits_{i_3=1}^{\alpha_{k_3}^{(\mathrm{G})}}\beta(\lambda_{k_1}, \lambda_{k_2}, \lambda_{k_3})\bm{P}_{k_1,i_1}\bm{Y}_1\bm{P}_{k_2,i_2}\bm{Y}_2\bm{P}_{k_3,i_3} \nonumber \\
&&+\sum\limits_{k_1=1}^{K_1}\sum\limits_{k_2=1}^{K_2}\sum\limits_{k_3=1}^{K_3}\sum\limits_{i_1=1}^{\alpha_{k_1}^{(\mathrm{G})}}\sum\limits_{i_2=1}^{\alpha_{k_2}^{(\mathrm{G})}}\sum\limits_{i_3=1}^{\alpha_{k_3}^{(\mathrm{G})}}\sum_{q_3=1}^{m_{k_3,i_3}-1}\frac{\beta^{(-,-,q_3)}(\lambda_{k_1},\lambda_{k_2},\lambda_{k_3})}{q_3!}\bm{P}_{k_1,i_1}\bm{Y}_1\bm{P}_{k_2,i_2}\bm{Y}_2\bm{N}_{k_3,i_3}^{q_3}\nonumber \\
&&+\sum\limits_{k_1=1}^{K_1}\sum\limits_{k_2=1}^{K_2}\sum\limits_{k_3=1}^{K_3}\sum\limits_{i_1=1}^{\alpha_{k_1}^{(\mathrm{G})}}\sum\limits_{i_2=1}^{\alpha_{k_2}^{(\mathrm{G})}}\sum\limits_{i_3=1}^{\alpha_{k_3}^{(\mathrm{G})}}\sum_{q_2=1}^{m_{k_2,i_2}-1}\frac{\beta^{(-,q_2,-)}(\lambda_{k_1},\lambda_{k_2},\lambda_{k_3})}{q_2!}\bm{P}_{k_1,i_1}\bm{Y}_1\bm{N}_{k_2,i_2}^{q_2}\bm{Y}_2\bm{P}_{k_3,i_3} \nonumber \\
&&+\sum\limits_{k_1=1}^{K_1}\sum\limits_{k_2=1}^{K_2}\sum\limits_{k_3=1}^{K_3}\sum\limits_{i_1=1}^{\alpha_{k_1}^{(\mathrm{G})}}\sum\limits_{i_2=1}^{\alpha_{k_2}^{(\mathrm{G})}}\sum\limits_{i_3=1}^{\alpha_{k_3}^{(\mathrm{G})}}\sum_{q_1=1}^{m_{k_1,i_1}-1}\frac{\beta^{(q_1,-,-)}(\lambda_{k_1},\lambda_{k_2},\lambda_{k_3})}{q_1!}\bm{N}_{k_1,i_1}^{q_1}\bm{Y}_1\bm{P}_{k_2,i_2}\bm{Y}_2\bm{P}_{k_3,i_3} \nonumber \\
&&+\sum\limits_{k_1=1}^{K_1}\sum\limits_{k_2=1}^{K_2}\sum\limits_{k_3=1}^{K_3}\sum\limits_{i_1=1}^{\alpha_{k_1}^{(\mathrm{G})}}\sum\limits_{i_2=1}^{\alpha_{k_2}^{(\mathrm{G})}}\sum\limits_{i_3=1}^{\alpha_{k_3}^{(\mathrm{G})}}\sum_{q_2=1}^{m_{k_2,i_2}-1}\sum_{q_3=1}^{m_{k_3,i_3}-1}\frac{\beta^{(-,q_2,q_3)}(\lambda_{k_1},\lambda_{k_2},\lambda_{k_3})}{q_2! q_3!}\bm{P}_{k_1,i_1}\bm{Y}_1\bm{N}_{k_2,i_2}^{q_2}\bm{Y}_2\bm{N}_{k_3,i_3}^{q_3}\nonumber \\
&&+\sum\limits_{k_1=1}^{K_1}\sum\limits_{k_2=1}^{K_2}\sum\limits_{k_3=1}^{K_3}\sum\limits_{i_1=1}^{\alpha_{k_1}^{(\mathrm{G})}}\sum\limits_{i_2=1}^{\alpha_{k_2}^{(\mathrm{G})}}\sum\limits_{i_3=1}^{\alpha_{k_3}^{(\mathrm{G})}}\sum_{q_1=1}^{m_{k_1,i_1}-1}\sum_{q_3=1}^{m_{k_3,i_3}-1}\frac{\beta^{(q_1,-,q_3)}(\lambda_{k_1},\lambda_{k_2},\lambda_{k_3})}{q_1!q_3!}\bm{N}_{k_1,i_1}^{q_1}\bm{Y}_1\bm{P}_{k_2,i_2}\bm{Y}_2\bm{N}_{k_3,i_3}^{q_3} \nonumber \\
&&+\sum\limits_{k_1=1}^{K_1}\sum\limits_{k_2=1}^{K_2}\sum\limits_{k_3=1}^{K_3}\sum\limits_{i_1=1}^{\alpha_{k_1}^{(\mathrm{G})}}\sum\limits_{i_2=1}^{\alpha_{k_2}^{(\mathrm{G})}}\sum\limits_{i_3=1}^{\alpha_{k_3}^{(\mathrm{G})}}\sum_{q_1=1}^{m_{k_1,i_1}-1}\sum_{q_2=1}^{m_{k_2,i_2}-1}\frac{\beta^{(q_1,q_2,-)}(\lambda_{k_1},\lambda_{k_2},\lambda_{k_3})}{q_1!q_2!}\bm{N}_{k_1,i_1}^{q_1}\bm{Y}_1\bm{N}_{k_2,i_2}^{q_2}\bm{Y}_2\bm{P}_{k_3,i_3} \nonumber \\
&&+\sum\limits_{k_1=1}^{K_1}\sum\limits_{k_2=1}^{K_2}\sum\limits_{k_3=1}^{K_3}\sum\limits_{i_1=1}^{\alpha_{k_1}^{(\mathrm{G})}}\sum\limits_{i_2=1}^{\alpha_{k_2}^{(\mathrm{G})}}\sum\limits_{i_3=1}^{\alpha_{k_3}^{(\mathrm{G})}}\sum_{q_1=1}^{m_{k_1,i_1}-1}\sum_{q_2=1}^{m_{k_2,i_2}-1}\sum_{q_3=1}^{m_{k_3,i_3}-1}\frac{\beta^{(q_1,q_2,q_3)}(\lambda_{k_1},\lambda_{k_2},\lambda_{k_3})}{q_1!q_2!q_3!}\nonumber \\
&&~~\times\bm{N}_{k_1,i_1}^{q_1}\bm{Y}_1\bm{N}_{k_2,i_2}^{q_2}\bm{Y}_2\bm{N}_{k_3,i_3}^{q_3} \nonumber \\
&=& \sum\limits_{i=1}^{2^3}T_{\beta}^{[\bm{X}]_{\xi(i)}}(\bm{Y}_1,\bm{Y}_2),
\end{eqnarray}
where $=_1$ comes from Proposition~\ref{prop:GMOI decomp by parameters X P X N}.


In GDOI, the composition of two GDOIs is still a GDOI. However, for GMOI, the composition of a GMOI with multiple GMOIs will change the number of parameter matrices for the composed GMOI. We have the following Propotion~\ref{prop:comp GMOIs} to characterize such composition behavior. 
\begin{proposition}\label{prop:comp GMOIs}
Given the following GMOIs: 
\begin{eqnarray}\label{eq1:prop:comp GMOIs}
T_{f}^{\bm{X}_1,\ldots,\bm{X}_{\zeta+1}}(\bm{Y}_1,\ldots,\bm{Y}_{\zeta}),
\end{eqnarray}
and
\begin{eqnarray}\label{eq2:prop:comp GMOIs}
T_{\beta_i}^{\bm{X}_1,\ldots,\bm{X}_{\zeta+1}}(\bm{Y}_1,\ldots,\bm{Y}_{\zeta}),
\end{eqnarray}
where $i=1,2,\ldots,\zeta$ and $[\bm{Y}]^{\zeta}_1$ is an array of matrices given by $[\bm{Y}]_1^{\zeta} \define \bm{Y}_1,\bm{Y}_2,\ldots,\bm{Y}_{\zeta}$. If the matrices $\bm{X}_p$ for $p=1,2,\ldots,\zeta+1$ are decomposed as:
\begin{eqnarray}\label{eq2.1:prop:comp GMOIs}
\bm{X}_p&=&\sum\limits_{k_p=1}^{K_p}\sum\limits_{i_p=1}^{\alpha_{k_p}^{\mathrm{G}}} \lambda_{k_p} \bm{P}_{k_p,i_p}+
\sum\limits_{k_p=1}^{K_p}\sum\limits_{i_p=1}^{\alpha_{k_p}^{\mathrm{G}}} \bm{N}_{k_p,i_p}.
\end{eqnarray}

Then, we have
\begin{eqnarray}\label{eq3:prop:comp GMOIs}
\lefteqn{T_{f \prod\limits_{i=1}^{\zeta}\beta_i}^{\bm{X}_1,[\bm{X}]_1^{\zeta+1},\bm{X}_2,[\bm{X}]_1^{\zeta+1},\ldots,[\bm{X}]_1^{\zeta+1},\bm{X}_{\zeta+1}}( \overbrace{\bm{I},[\bm{Y}]^{\zeta}_1, \bm{I},\ldots,\bm{I},[\bm{Y}]^{\zeta}_1,\bm{I}}^{\mbox{there are $\zeta$ terms of $\bm{I},[\bm{Y}]^{\zeta}_1, \bm{I}$}} )}\nonumber \\
&=&T_{f}^{\bm{X}_1,\ldots,\bm{X}_{\zeta+1}}(T_{\beta_1}^{\bm{X}_1,\ldots,\bm{X}_{\zeta+1}}(\bm{Y}_1,\ldots,\bm{Y}_{\zeta}),\ldots,T_{\beta_\zeta}^{\bm{X}_1,\ldots,\bm{X}_{\zeta+1}}(\bm{Y}_1,\ldots,\bm{Y}_{\zeta})),
\end{eqnarray}
where $[\bm{X}]_1^{\zeta+1}$ is an array of matrices given by $[\bm{X}]_1^{\zeta+1} \define \bm{X}_1,\bm{X}_2,\ldots,\bm{X}_{\zeta+1}$. Note that there are $(\zeta+1)^2$ parameter matrices in GMOI given by L.H.S. of Eq.~\eqref{eq3:prop:comp GMOIs}
\end{proposition}
\textbf{Proof:}
From GMOI definition given by Eq.~\eqref{eq1:GMOI def}, we have 
\begin{eqnarray}\label{eq4:prop:comp GMOIs}
\lefteqn{T_{f}^{\bm{X}_1,\ldots,\bm{X}_{\zeta+1}}(T_{\beta_1}^{\bm{X}_1,\ldots,\bm{X}_{\zeta+1}}(\bm{Y}_1,\ldots,\bm{Y}_{\zeta}),\ldots,T_{\beta_\zeta}^{\bm{X}_1,\ldots,\bm{X}_{\zeta+1}}(\bm{Y}_1,\ldots,\bm{Y}_{\zeta}))}\nonumber \\
&=& \sum\limits_{k_1=\ldots=k_{\zeta+1}=1}^{K_1,\ldots,K_{\zeta+1}} \sum\limits_{i_1=\ldots=i_{\zeta+1}=1}^{\alpha_{k_1}^{(\mathrm{G})},\ldots,\alpha_{k_{\zeta+1}}^{(\mathrm{G})}}
f(\lambda_{k_1},\ldots,\lambda_{k_{\zeta+1}})\bm{P}_{k_1,i_1}T_{\beta_1}^{\bm{X}_1,\ldots,\bm{X}_{\zeta+1}}(\bm{Y}_1,\ldots,\bm{Y}_{\zeta})\nonumber \\
&&~~\times \bm{P}_{k_2,i_2}T_{\beta_2}^{\bm{X}_1,\ldots,\bm{X}_{\zeta+1}}(\bm{Y}_1,\ldots,\bm{Y}_{\zeta}) \ldots T_{\beta_\zeta}^{\bm{X}_1,\ldots,\bm{X}_{\zeta+1}}(\bm{Y}_1,\ldots,\bm{Y}_{\zeta})\bm{P}_{k_{\zeta+1},i_{\zeta+1}}\nonumber \\
&&+\sum\limits_{k_1=\ldots=k_{\zeta+1}=1}^{K_1,\ldots,K_{\zeta+1}} \sum\limits_{i_1=\ldots=i_{\zeta+1}=1}^{\alpha_{k_1}^{(\mathrm{G})},\ldots,\alpha_{k_{\zeta+1}}^{(\mathrm{G})}}\sum\limits_{\kappa=1}^{\zeta}\sum\limits_{\varrho_\kappa(q_1,\ldots,q_{\zeta+1})}\Bigg(\sum\limits_{\varrho_{\kappa}(q_1,\ldots,q_{\zeta+1})=1}^{m_{k_{\mbox{Ind}(\varrho_{\kappa}(q_1,\ldots,q_{\zeta+1}))},i_{\mbox{Ind}(\varrho_{\kappa}(q_1,\ldots,q_{\zeta+1}))}}-1}\nonumber \\
&&~~~~~ \frac{f^{\varrho_{\kappa}(q_1,\ldots,q_{\zeta+1})}(\lambda_{k_1},\ldots,\lambda_{k_{\zeta+1}})}{q_{\iota_1}!q_{\iota_2}!\ldots q_{\iota_\kappa}!}\times \prod\limits_{\substack{\varsigma =\mbox{Ind}(\varrho_{\kappa}(q_1,\ldots,q_r)), \bm{Z}_\varsigma=\bm{N}^{q_\varsigma}_{k_\varsigma,i_\varsigma}T_{\beta_\varsigma}^{\bm{X}_1,\ldots,\bm{X}_{\zeta+1}}(\bm{Y}_1,\ldots,\bm{Y}_{\zeta})\\ \varsigma \neq \mbox{Ind}(\varrho_{\kappa}(q_1,\ldots,q_r)), \bm{Z}_\varsigma=\bm{P}_{k_\varsigma,i_\varsigma}T_{\beta_\varsigma}^{\bm{X}_1,\ldots,\bm{X}_{\zeta+1}}(\bm{Y}_1,\ldots,\bm{Y}_{\zeta})}
}^{\zeta+1} \bm{Z}_\varsigma\Bigg) 
\nonumber \\
&&+\sum\limits_{k_1=\ldots=k_{\zeta+1}=1}^{K_1,\ldots,K_{\zeta+1}} \sum\limits_{i_1=\ldots=i_{\zeta+1}=1}^{\alpha_{k_1}^{(\mathrm{G})},\ldots,\alpha_{k_{\zeta+1}}^{(\mathrm{G})}} \sum\limits_{q_1=\ldots=q_{\zeta+1}=1}^{m_{k_1,i_1}-1,\ldots,m_{k_{\zeta+1},i_{\zeta+1}}-1}
\frac{f^{(q_1,\ldots,q_{\zeta+1})}(\lambda_{k_1},\ldots,\lambda_{k_{\zeta+1}})}{q_1!\cdots q_{\zeta+1}!}\nonumber \\
&&
~~\times \bm{N}^{q_1}_{k_1,i_1}T_{\beta_1}^{\bm{X}_1,\ldots,\bm{X}_{\zeta+1}}(\bm{Y}_1,\ldots,\bm{Y}_{\zeta})\bm{N}^{q_2}_{k_2,i_2}T_{\beta_2}^{\bm{X}_1,\ldots,\bm{X}_{\zeta+1}}(\bm{Y}_1,\ldots,\bm{Y}_{\zeta})\ldots\nonumber \\
&& ~~ \times \bm{N}^{q_{\zeta}}_{k_{\zeta},i_{\zeta}} T_{\beta_\zeta}^{\bm{X}_1,\ldots,\bm{X}_{\zeta+1}}(\bm{Y}_1,\ldots,\bm{Y}_{\zeta})\bm{N}^{q_{\zeta+1}}_{k_{\zeta+1},i_{\zeta+1}},
\end{eqnarray}
where we set $T_{\beta_{\zeta+1}}^{\bm{X}_1,\ldots,\bm{X}_{\zeta+1}}(\bm{Y}_1,\ldots,\bm{Y}_{\zeta})$ as the identity matrix.

By observing each term in Eq.~\eqref{eq4:prop:comp GMOIs}, we have the following format regarding to parameter matrices
\begin{eqnarray}\label{eq5:prop:comp GMOIs}
\left\{
\begin{array}{l}
\bm{P}_{k_1,i_1} \\
\bm{N}^{q_1}_{k_1,i_1}
\end{array}
\right\}, [\bm{X}]_1^{\zeta+1}, \left\{
\begin{array}{l}
\bm{P}_{k_2,i_2} \\
\bm{N}^{q_2}_{k_2,i_2}
\end{array}
\right\}, [\bm{X}]_1^{\zeta+1} \ldots  [\bm{X}]_1^{\zeta+1}, 
\left\{
\begin{array}{l}
\bm{P}_{k_{\zeta+1},i_{\zeta+1}} \\
\bm{N}^{q_{\zeta+1}}_{k_{\zeta+1},i_{\zeta+1}}
\end{array}
\right\},
\end{eqnarray}
where each $\left\{
\begin{array}{l}
\bm{P}_{k_p,i_p} \\
\bm{N}^{q_p}_{k_p,i_p}
\end{array}
\right\}$ comes from the projection or nilpotent parts from the matrix $\bm{X}_p$. From Eq.\eqref{eq5:prop:comp GMOIs}, we obtain the argument matrices $[\bm{Y}]^{\zeta}_1$ corresponding to each occurrence of $[\bm{X}]_1^{\zeta+1}$. To align Eq.\eqref{eq5:prop:comp GMOIs} with each summand term in Eq.\eqref{eq4:prop:comp GMOIs}, we insert an identity matrix both before and after each $[\bm{X}]_1^{\zeta+1}$. Accordingly, the full sequence of argument matrices can be written as:
\begin{eqnarray}
\overbrace{\bm{I},[\bm{Y}]^{\zeta}_1, \bm{I},\ldots,\bm{I},[\bm{Y}]^{\zeta}_1,\bm{I}}^{\mbox{there are $\zeta$ terms of $\bm{I},[\bm{Y}]^{\zeta}_1, \bm{I}$}}.
\end{eqnarray}
Finally, since the complex-valued functions $f$ and $\beta_i$ are commutative, Eq.\eqref{eq4:prop:comp GMOIs} becomes:
\begin{eqnarray}
T_{f \prod\limits_{i=1}^{\zeta}\beta_i}^{\bm{X}_1,[\bm{X}]_1^{\zeta+1},\bm{X}_2,[\bm{X}]_1^{\zeta+1},\ldots,[\bm{X}]1^{\zeta+1},\bm{X}{\zeta+1}}( \overbrace{\bm{I},[\bm{Y}]^{\zeta}_1, \bm{I},\ldots,\bm{I},[\bm{Y}]^{\zeta}_1,\bm{I}}^{\mbox{there are $\zeta$ terms of $\bm{I},[\bm{Y}]^{\zeta}_1, \bm{I}$}}).
\end{eqnarray}

$\hfill\Box$

\section{GMOI Norm Estimations}\label{sec: GMOI Norm Estimations}

The purpose of this section is to provide norm estimates for GMOI. We begin by introducing some index mapping notations. Given a non-negative integer $i$, we use $B(i,\zeta+1)$ to represent the binary representation of the integer $i$ by $\zeta+1$ bits. For example, we have $B(0,2)=00$, $B(1,2)=01$,  $B(2,2)=10$, and $B(3,2)=11$. We define the following map for the binary representation $B(i,\zeta+1)$ associated to matrices $[\bm{X}]_1^{\zeta+1}$ as:
\begin{eqnarray}
\Psi_{[\bm{X}]_1^{\zeta+1}}(B(i,\zeta+1)) \rightarrow \sum\limits_{j=1}^\kappa q_{\iota_j} \bm{e}^{(\zeta+1)}_{\iota_j},
\end{eqnarray}
where $\kappa$ is the number of one in the binary representation $B(i,\zeta+1)$, $q_{\iota_j}$ will be the exponent for the nilpotent part of the matrix $\bm{X}_{\iota_j}$, $\iota_j$ is the position (leftmost position is indexed by $1$) of $B(i,\zeta+1)$ with value one,  and $\bm{e}^{(\zeta+1)}_{\iota_j}$ is the standard basis vector with the length $\zeta+1$ and all zero entries except the value one at the position$\iota_j$. For example, we have $\Psi_{\bm{X}_1, \bm{X}_2}(B(0,2)) =(0,0)$, $\Psi_{\bm{X}_1, \bm{X}_2}(B(1,2))=(0,q_2)$,$\Psi_{\bm{X}_1, \bm{X}_2}(B(2,2))=(q_1,0)$, and $\Psi_{\bm{X}_1, \bm{X}_2}(B(3,2))=(q_1,q_2)$.

From GMOI definition given by Eq.~\eqref{eq1:GMOI def}, we can express GMOI by
\begin{eqnarray}\label{eq1:GMOI sum form by Ai prime}
T_{\beta}^{\bm{X}_1,\ldots,\bm{X}_{\zeta+1}}(\bm{Y}_1,\ldots,\bm{Y}_{\zeta})&=&\sum\limits_{i'=0}^{2^{\zeta+1}-1}\bm{A}_{i'}, 
\end{eqnarray}
where $\bm{A}_{i'}$ is a matrix, which can be represented by
\begin{eqnarray}\label{eq2:GMOI sum form by Ai prime}
\bm{A}_{i'}&=&\sum\limits_{k_1=\ldots=k_{\zeta+1}=1}^{K_1,\ldots,K_{\zeta+1}} \sum\limits_{i_1=\ldots=i_{\zeta+1}=1}^{\alpha_{k_1}^{(\mathrm{G})},\ldots,\alpha_{k_{\zeta+1}}^{(\mathrm{G})}}
\sum\limits_{\tilde{\Psi}_{[\bm{X}]_1^{\zeta+1}}(B(i',\zeta+1))=1}^{m_{k_{\mbox{Ind}(\tilde{\Psi}_{[\bm{X}]_1^{\zeta+1}}(B(i',\zeta+1)))},i_{\mbox{Ind}(\tilde{\Psi}_{[\bm{X}]_1^{\zeta+1}}(B(i',\zeta+1)))}}-1}\nonumber \\
&& 
\frac{\beta^{\tilde{\Psi}_{[\bm{X}]_1^{\zeta+1}}(B(i',\zeta+1))}(\lambda_{k_1},\ldots,\lambda_{k_{\zeta+1}}) }{\tilde{\Psi}_{[\bm{X}]_1^{\zeta+1}}(B(i',\zeta+1))!} \times  \prod\limits_{\substack{\varsigma =\mbox{Ind}(\tilde{\Psi}_{[\bm{X}]_1^{\zeta+1}}(B(i',\zeta+1))), \bm{Z}_\varsigma=\bm{N}^{q_\varsigma}_{k_\varsigma,i_\varsigma}\bm{Y}_{\varsigma} \\ \varsigma \neq \mbox{Ind}(\tilde{\Psi}_{[\bm{X}]_1^{\zeta+1}}(B(i',\zeta+1))), \bm{Z}_\varsigma=\bm{P}_{k_\varsigma,i_\varsigma}\bm{Y}_{\varsigma}}
}^{\zeta+1} \bm{Z}_\varsigma
\end{eqnarray}
where $\tilde{\Psi}_{[\bm{X}]_1^{\zeta+1}}(B(i',\zeta+1))$ are those non-zero entries of the vector $\sum\limits_{j=1}^\kappa q_{\iota_j} \bm{e}^{(\zeta+1)}_{\iota_j}$,\\
$\beta^{\Psi_{[\bm{X}]_1^{\zeta+1}}(B(i',\zeta+1))}(\lambda_{k_1},\ldots,\lambda_{k_{\zeta+1}})$ is the partial derivatives with respect to orders $\tilde{\Psi}_{[\bm{X}]_1^{\zeta+1}}(B(i',\zeta+1))$ for the funtion $\beta(\lambda_{k_1},\ldots,\lambda_{k_{\zeta+1}})$, $\tilde{\Psi}_{[\bm{X}]_1^{\zeta+1}}(B(i',\zeta+1))=1$ indicates those $q_{\iota_j}$ (summand variables) set as $1$, $\tilde{\Psi}_{[\bm{X}]_1^{\zeta+1}}(B(i',\zeta+1))!$ is a product of factorial of those nilpotent exponents $q_{\iota_j}$. We still assume that $\bm{Y}_{\zeta+1}=\bm{I}$.

Let us consider the following Example~\ref{exp:GMOI sum form by Ai prime for GTOI} about expressing $T_{\beta}^{\bm{X}_1,\bm{X}_2,\bm{X}_3}(\bm{Y}_1, \bm{Y}_2)$ by Eq.~\eqref{eq1:GMOI sum form by Ai prime}.
\begin{example}\label{exp:GMOI sum form by Ai prime for GTOI}
From Eq.~\eqref{eq1:  GTOI def}, we have
\begin{eqnarray}\label{eq1:exp:GMOI sum form by Ai prime for GTOI}
\lefteqn{T_{\beta}^{\bm{X}_1,\bm{X}_2,\bm{X}_3}(\bm{Y}_1, \bm{Y}_2)\define}\nonumber\\
&&\sum\limits_{k_1=1}^{K_1}\sum\limits_{k_2=1}^{K_2}\sum\limits_{k_3=1}^{K_3}\sum\limits_{i_1=1}^{\alpha_{k_1}^{(\mathrm{G})}}\sum\limits_{i_2=1}^{\alpha_{k_2}^{(\mathrm{G})}}\sum\limits_{i_3=1}^{\alpha_{k_3}^{(\mathrm{G})}}\beta(\lambda_{k_1}, \lambda_{k_2}, \lambda_{k_3})\bm{P}_{k_1,i_1}\bm{Y}_1\bm{P}_{k_2,i_2}\bm{Y}_2\bm{P}_{k_3,i_3} \nonumber \\
&&+\sum\limits_{k_1=1}^{K_1}\sum\limits_{k_2=1}^{K_2}\sum\limits_{k_3=1}^{K_3}\sum\limits_{i_1=1}^{\alpha_{k_1}^{(\mathrm{G})}}\sum\limits_{i_2=1}^{\alpha_{k_2}^{(\mathrm{G})}}\sum\limits_{i_3=1}^{\alpha_{k_3}^{(\mathrm{G})}}\sum_{q_3=1}^{m_{k_3,i_3}-1}\frac{\beta^{(-,-,q_3)}(\lambda_{k_1},\lambda_{k_2},\lambda_{k_3})}{q_3!}\bm{P}_{k_1,i_1}\bm{Y}_1\bm{P}_{k_2,i_2}\bm{Y}_2\bm{N}_{k_3,i_3}^{q_3}\nonumber \\
&&+\sum\limits_{k_1=1}^{K_1}\sum\limits_{k_2=1}^{K_2}\sum\limits_{k_3=1}^{K_3}\sum\limits_{i_1=1}^{\alpha_{k_1}^{(\mathrm{G})}}\sum\limits_{i_2=1}^{\alpha_{k_2}^{(\mathrm{G})}}\sum\limits_{i_3=1}^{\alpha_{k_3}^{(\mathrm{G})}}\sum_{q_2=1}^{m_{k_2,i_2}-1}\frac{\beta^{(-,q_2,-)}(\lambda_{k_1},\lambda_{k_2},\lambda_{k_3})}{q_2!}\bm{P}_{k_1,i_1}\bm{Y}_1\bm{N}_{k_2,i_2}^{q_2}\bm{Y}_2\bm{P}_{k_3,i_3} \nonumber \\
&&+\sum\limits_{k_1=1}^{K_1}\sum\limits_{k_2=1}^{K_2}\sum\limits_{k_3=1}^{K_3}\sum\limits_{i_1=1}^{\alpha_{k_1}^{(\mathrm{G})}}\sum\limits_{i_2=1}^{\alpha_{k_2}^{(\mathrm{G})}}\sum\limits_{i_3=1}^{\alpha_{k_3}^{(\mathrm{G})}}\sum_{q_1=1}^{m_{k_1,i_1}-1}\frac{\beta^{(q_1,-,-)}(\lambda_{k_1},\lambda_{k_2},\lambda_{k_3})}{q_1!}\bm{N}_{k_1,i_1}^{q_1}\bm{Y}_1\bm{P}_{k_2,i_2}\bm{Y}_2\bm{P}_{k_3,i_3} \nonumber \\
&&+\sum\limits_{k_1=1}^{K_1}\sum\limits_{k_2=1}^{K_2}\sum\limits_{k_3=1}^{K_3}\sum\limits_{i_1=1}^{\alpha_{k_1}^{(\mathrm{G})}}\sum\limits_{i_2=1}^{\alpha_{k_2}^{(\mathrm{G})}}\sum\limits_{i_3=1}^{\alpha_{k_3}^{(\mathrm{G})}}\sum_{q_2=1}^{m_{k_2,i_2}-1}\sum_{q_3=1}^{m_{k_3,i_3}-1}\frac{\beta^{(-,q_2,q_3)}(\lambda_{k_1},\lambda_{k_2},\lambda_{k_3})}{q_2! q_3!}\bm{P}_{k_1,i_1}\bm{Y}_1\bm{N}_{k_2,i_2}^{q_2}\bm{Y}_2\bm{N}_{k_3,i_3}^{q_3}\nonumber \\
&&+\sum\limits_{k_1=1}^{K_1}\sum\limits_{k_2=1}^{K_2}\sum\limits_{k_3=1}^{K_3}\sum\limits_{i_1=1}^{\alpha_{k_1}^{(\mathrm{G})}}\sum\limits_{i_2=1}^{\alpha_{k_2}^{(\mathrm{G})}}\sum\limits_{i_3=1}^{\alpha_{k_3}^{(\mathrm{G})}}\sum_{q_1=1}^{m_{k_1,i_1}-1}\sum_{q_3=1}^{m_{k_3,i_3}-1}\frac{\beta^{(q_1,-,q_3)}(\lambda_{k_1},\lambda_{k_2},\lambda_{k_3})}{q_1!q_3!}\bm{N}_{k_1,i_1}^{q_1}\bm{Y}_1\bm{P}_{k_2,i_2}\bm{Y}_2\bm{N}_{k_3,i_3}^{q_3} \nonumber \\
&&+\sum\limits_{k_1=1}^{K_1}\sum\limits_{k_2=1}^{K_2}\sum\limits_{k_3=1}^{K_3}\sum\limits_{i_1=1}^{\alpha_{k_1}^{(\mathrm{G})}}\sum\limits_{i_2=1}^{\alpha_{k_2}^{(\mathrm{G})}}\sum\limits_{i_3=1}^{\alpha_{k_3}^{(\mathrm{G})}}\sum_{q_1=1}^{m_{k_1,i_1}-1}\sum_{q_2=1}^{m_{k_2,i_2}-1}\frac{\beta^{(q_1,q_2,-)}(\lambda_{k_1},\lambda_{k_2},\lambda_{k_3})}{q_1!q_2!}\bm{N}_{k_1,i_1}^{q_1}\bm{Y}_1\bm{N}_{k_2,i_2}^{q_2}\bm{Y}_2\bm{P}_{k_3,i_3} \nonumber \\
&&+\sum\limits_{k_1=1}^{K_1}\sum\limits_{k_2=1}^{K_2}\sum\limits_{k_3=1}^{K_3}\sum\limits_{i_1=1}^{\alpha_{k_1}^{(\mathrm{G})}}\sum\limits_{i_2=1}^{\alpha_{k_2}^{(\mathrm{G})}}\sum\limits_{i_3=1}^{\alpha_{k_3}^{(\mathrm{G})}}\sum_{q_1=1}^{m_{k_1,i_1}-1}\sum_{q_2=1}^{m_{k_2,i_2}-1}\sum_{q_3=1}^{m_{k_3,i_3}-1}\frac{\beta^{(q_1,q_2,q_3)}(\lambda_{k_1},\lambda_{k_2},\lambda_{k_3})}{q_1!q_2!q_3!}\nonumber \\
&&~~\times\bm{N}_{k_1,i_1}^{q_1}\bm{Y}_1\bm{N}_{k_2,i_2}^{q_2}\bm{Y}_2\bm{N}_{k_3,i_3}^{q_3}.
\end{eqnarray}

Because we have 
\begin{eqnarray}\label{eq2:exp:GMOI sum form by Ai prime for GTOI}
B(0,3)&=&(0,0,0) \rightarrow \Psi_{[\bm{X}]_1^{\zeta+1}}(B(0,3))=(0,0,0); \nonumber \\
B(1,3)&=&(0,0,1) \rightarrow \Psi_{[\bm{X}]_1^{\zeta+1}}(B(1,3))=(0,0,q_3); \nonumber \\
B(2,3)&=&(0,1,0) \rightarrow \Psi_{[\bm{X}]_1^{\zeta+1}}(B(2,3))=(0,q_2,0); \nonumber \\
B(3,3)&=&(0,1,1) \rightarrow \Psi_{[\bm{X}]_1^{\zeta+1}}(B(3,3))=(0,q_2,q_3); \nonumber \\
B(4,3)&=&(1,0,0) \rightarrow \Psi_{[\bm{X}]_1^{\zeta+1}}(B4,3))=(q_1,0,0); \nonumber \\
B(5,3)&=&(1,0,1) \rightarrow \Psi_{[\bm{X}]_1^{\zeta+1}}(B(5,3))=(q_1,0,q_3); \nonumber \\
B(6,3)&=&(1,1,0) \rightarrow \Psi_{[\bm{X}]_1^{\zeta+1}}(B(6,3))=(q_1,q_2,0); \nonumber \\
B(7,3)&=&(1,1,1) \rightarrow \Psi_{[\bm{X}]_1^{\zeta+1}}(B(7,3))=(q_1,q_2,q_3),
\end{eqnarray}
then, 
\begin{eqnarray}\label{eq3:exp:GMOI sum form by Ai prime for GTOI}
B(0,3)&=&(0,0,0) \rightarrow \tilde{\Psi}_{[\bm{X}]_1^{\zeta+1}}(B(0,3))=\emptyset; \nonumber \\
B(1,3)&=&(0,0,1) \rightarrow \tilde{\Psi}_{[\bm{X}]_1^{\zeta+1}}(B(1,3))=(q_3) \rightarrow m_{k_3,i_3}; \nonumber \\
B(2,3)&=&(0,1,0) \rightarrow \tilde{\Psi}_{[\bm{X}]_1^{\zeta+1}}(B(2,3))=(q_2) \rightarrow m_{k_2,i_2};\nonumber \\
B(3,3)&=&(0,1,1) \rightarrow \tilde{\Psi}_{[\bm{X}]_1^{\zeta+1}}(B(3,3))=(q_2,q_3) \rightarrow m_{k_2,i_2}, m_{k_3,i_3}; \nonumber \\
B(4,3)&=&(1,0,0) \rightarrow \tilde{\Psi}_{[\bm{X}]_1^{\zeta+1}}(B4,3))=(q_1) \rightarrow m_{k_1,i_1};\nonumber \\
B(5,3)&=&(1,0,1) \rightarrow \tilde{\Psi}_{[\bm{X}]_1^{\zeta+1}}(B(5,3))=(q_1,q_3) \rightarrow m_{k_1,i_1}, m_{k_3,i_3}; \nonumber \\
B(6,3)&=&(1,1,0) \rightarrow \tilde{\Psi}_{[\bm{X}]_1^{\zeta+1}}(B(6,3))=(q_1,q_2) \rightarrow m_{k_1,i_1}, m_{k_2,i_2}; \nonumber \\
B(7,3)&=&(1,1,1) \rightarrow \tilde{\Psi}_{[\bm{X}]_1^{\zeta+1}}(B(7,3))=(q_1,q_2,q_3) \rightarrow m_{k_1,i_1}, m_{k_2,i_2}, m_{k_3,i_3}.
\end{eqnarray}
According to Eq.~\eqref{eq2:GMOI sum form by Ai prime}, we can express $T_{\beta}^{\bm{X}_1,\bm{X}_2,\bm{X}_3}(\bm{Y}_1, \bm{Y}_2)$ as
\begin{eqnarray}
T_{\beta}^{\bm{X}_1,\bm{X}_2,\bm{X}_3}(\bm{Y}_1, \bm{Y}_2)&=&\bm{A}_0+\bm{A}_1+\bm{A}_2+\bm{A}_4+\bm{A}_3+\bm{A}_5+\bm{A}_6+\bm{A}_7.
\end{eqnarray}
\end{example}

We recall Lemma 3 in~\cite{chang2025GDOIMatrix} for the converse triangle inequality for the Frobenius norm, denoted by $\left\Vert \cdot \right\Vert$.
\begin{lemma}\label{lma:conv triangle for Frob norm}
Given $n$ matrices $\bm{A}_1, \bm{A}_2, \ldots, \bm{A}_n$ such that $\left\Vert\bm{A}_1\right\Vert\geq\left\Vert\bm{A}_2\right\Vert\geq\ldots\geq\left\Vert\bm{A}_n\right\Vert$, then, we have
\begin{eqnarray}\label{eq1:lma:conv triangle for Frob norm}
\left\Vert\sum\limits_{i=1}^n \bm{A}_i \right\Vert \geq \max\left[0, \left\Vert\bm{A}_1\right\Vert - \sum\limits_{i=2}^n \left\Vert\bm{A}_i\right\Vert\right]
\end{eqnarray}
\end{lemma}

Theorem~\ref{thm: GMOI norm Est} will be provided to give the lower bound and the upper bound estimations for the Frobenius norm of $T_{\beta}^{\bm{X}_1,\ldots,\bm{X}_{\zeta+1}}(\bm{Y}_1,\ldots,\bm{Y}_{\zeta})$.
\begin{theorem}\label{thm: GMOI norm Est}
We have the upper bound for the Frobenius norm of $T_{\beta}^{\bm{X}_1,\ldots,\bm{X}_{\zeta+1}}(\bm{Y}_1,\ldots,\bm{Y}_{\zeta})=\sum\limits_{i'=0}^{2^{\zeta+1}-1}\bm{A}_{i'}$, which is given by
\begin{eqnarray}\label{eq1:  thm: GMOI Lip Est}
\left\Vert T_{\beta}^{\bm{X}_1,\ldots,\bm{X}_{\zeta+1}}(\bm{Y}_1,\ldots,\bm{Y}_{\zeta}) \right\Vert\leq
\sum\limits_{i'=0}^{2^{\zeta+1}-1} \left\Vert \bm{A}_{i'} \right\Vert_{up},
\end{eqnarray}
where $\left\Vert \bm{A}_{i'} \right\Vert_{up}$ are upper bounds for  the Frobenius norm of the matrix $\bm{A}_{i'}$, which is defined by
\begin{eqnarray}\label{eq1.0:  thm: GMOI norm Est}
\left\Vert \bm{A}_{i'} \right\Vert_{up}&=& \sum\limits_{ k_\varsigma=1,  \mbox{for}~\varsigma=\mbox{Ind}(\tilde{\Psi}_{[\bm{X}]_1^{\zeta+1}}(B(i',\zeta+1)))}^{K_\varsigma, \mbox{for}~\varsigma=\mbox{Ind}(\tilde{\Psi}_{[\bm{X}]_1^{\zeta+1}}(B(i',\zeta+1)))} \sum\limits_{i_\varsigma=1,  \mbox{for}~\varsigma=\mbox{Ind}(\tilde{\Psi}_{[\bm{X}]_1^{\zeta+1}}(B(i',\zeta+1)))}^{\alpha_{k_\varsigma}^{(\mathrm{G})}, \mbox{for}~\varsigma=\mbox{Ind}(\tilde{\Psi}_{[\bm{X}]_1^{\zeta+1}}(B(i',\zeta+1)))}\nonumber \\
&& 
\sum\limits_{\tilde{\Psi}_{[\bm{X}]_1^{\zeta+1}}(B(i',\zeta+1))=1}^{m_{k_{\mbox{Ind}(\tilde{\Psi}_{[\bm{X}]_1^{\zeta+1}}(B(i',\zeta+1)))},i_{\mbox{Ind}(\tilde{\Psi}_{[\bm{X}]_1^{\zeta+1}}(B(i',\zeta+1)))}}-1}\nonumber \\
&& 
\left[ \max\limits_{\lambda_{k_1} \in \Lambda_{\bm{X}_1},\ldots,\lambda_{k_{\zeta+1}} \in \Lambda_{\bm{X}_{\zeta+1}}} \left\vert \frac{\beta^{(\tilde{\Psi}_{[\bm{X}]_1^{\zeta+1}}(B(i',\zeta+1)))}(\lambda_{k_1},\ldots,\lambda_{k_{\zeta+1}}) }{\tilde{\Psi}_{[\bm{X}]_1^{\zeta+1}}(B(i',\zeta+1))!}\right\vert\right] \nonumber \\
&& \times  \prod\limits_{\substack{\varsigma =\mbox{Ind}(\tilde{\Psi}_{[\bm{X}]_1^{\zeta+1}}(B(i',\zeta+1))), \bm{Z}_\varsigma=\left\Vert \bm{N}^{q_\varsigma}_{k_\varsigma,i_\varsigma} \right\Vert \left\Vert \bm{Y}_{\varsigma}\right\Vert\\ \varsigma \neq \mbox{Ind}(\tilde{\Psi}_{[\bm{X}]_1^{\zeta+1}}(B(i',\zeta+1))), \bm{Z}_\varsigma=\left\Vert\bm{Y}_{\varsigma}\right\Vert}
}^{\zeta+1} \bm{Z}_\varsigma
\end{eqnarray}

On the other hand, we have the lower bound for the Frobenius norm of $T_{\beta}^{\bm{X}_1,\ldots,\bm{X}_{\zeta+1}}(\bm{Y}_1,\ldots,\bm{Y}_{\zeta})$, which is given by
\begin{eqnarray}\label{eq2:  thm: GMOI norm Est}
\left\Vert T_{\beta}^{\bm{X}_1,\ldots,\bm{X}_{\zeta+1}}(\bm{Y}_1,\ldots,\bm{Y}_{\zeta}) \right\Vert\geq \max\left[0, \left\Vert\bm{A}_{\sigma(0)}\right\Vert - \sum\limits_{i'=1}^{2^{\zeta+1}-1} \left\Vert\bm{A}_{\sigma(i')}\right\Vert\right],
\end{eqnarray}
where $\sigma$ is the permutation of matrices $\bm{A}_{i'}$ for $i'=0,1,\ldots,2^{\zeta+1}-1$ such that $\left\Vert\bm{A}_{\sigma(0)}\right\Vert \geq \left\Vert\bm{A}_{\sigma(1)}\right\Vert \geq \ldots \geq \left\Vert\bm{A}_{\sigma(2^{\zeta+1}-1)}\right\Vert$. 

Further, if we have $\left[\min\limits_{\lambda_1 \in \Lambda_{\bm{X}_1},\ldots,\lambda_{\zeta+1} \in \Lambda_{\bm{X}_{\zeta+1}}} \left\vert\beta(\lambda_1,\dots, \lambda_{\zeta+1})\right\vert\right]\left\Vert \prod\limits_{i'=1}^{\zeta}\bm{Y}_{i'}\right\Vert\geq \sum\limits_{i'=1}^{2^{\zeta+1}-1}\left\Vert\bm{A}_{i'}\right\Vert$,  the lower bound for the Frobenius norm of $T_{\beta}^{\bm{X}_1,\bm{X}_2}(\bm{Y})$ can be expressed by
\begin{eqnarray}\label{eq3:  thm: GMOI norm Est}
\left\Vert T_{\beta}^{\bm{X}_1,\ldots,\bm{X}_{\zeta+1}}(\bm{Y}_1,\ldots,\bm{Y}_{\zeta}) \right\Vert\geq \left[\min\limits_{\lambda_1 \in \Lambda_{\bm{X}_1},\ldots,\lambda_{\zeta+1} \in \Lambda_{\bm{X}_{\zeta+1}}} \left\vert\beta(\lambda_1,\dots, \lambda_{\zeta+1})\right\vert\right]\left\Vert  \prod\limits_{i'=1}^{\zeta}\bm{Y}_{i'} \right\Vert- \sum\limits_{i'=1}^{2^{\zeta+1}-1}\left\Vert\bm{A}_{i'}\right\Vert.
\end{eqnarray}
\end{theorem}
\textbf{Proof:}
From Eq.~\eqref{eq1:GMOI sum form by Ai prime}, we have
\begin{eqnarray}\label{eq4:thm: GMOI norm Est}
T_{\beta}^{\bm{X}_1,\ldots,\bm{X}_{\zeta+1}}(\bm{Y}_1,\ldots,\bm{Y}_{\zeta})&=&\sum\limits_{i'=0}^{2^{\zeta+1}-1}\bm{A}_{i'}, 
\end{eqnarray}
and we immediate have the following by the triangle inequality of Frobenius norm:
\begin{eqnarray}\label{eq5:thm: GMOI norm Est}
\left\Vert T_{\beta}^{\bm{X}_1,\ldots,\bm{X}_{\zeta+1}}(\bm{Y}_1,\ldots,\bm{Y}_{\zeta}) \right\Vert&\leq&\sum\limits_{i'=0}^{2^{\zeta+1}-1}\left\Vert \bm{A}_{i'} \right\Vert.  
\end{eqnarray}
From the $\bm{A}_{i'}$ expression given by Eq.~\eqref{eq2:GMOI sum form by Ai prime}, we have
{\small
\begin{eqnarray}\label{eq6:thm: GMOI norm Est}
\left\Vert\bm{A}_{i'}\right\Vert&=&\left\Vert\sum\limits_{k_1=\ldots=k_{\zeta+1}=1}^{K_1,\ldots,K_{\zeta+1}} \sum\limits_{i_1=\ldots=i_{\zeta+1}=1}^{\alpha_{k_1}^{(\mathrm{G})},\ldots,\alpha_{k_{\zeta+1}}^{(\mathrm{G})}}
\sum\limits_{\tilde{\Psi}_{[\bm{X}]_1^{\zeta+1}}(B(i',\zeta+1))=1}^{m_{k_{\mbox{Ind}(\tilde{\Psi}_{[\bm{X}]_1^{\zeta+1}}(B(i',\zeta+1)))},i_{\mbox{Ind}(\tilde{\Psi}_{[\bm{X}]_1^{\zeta+1}}(B(i',\zeta+1)))}}-1} \right.\nonumber \\
&& 
\left. \frac{\beta^{(\Psi_{[\bm{X}]_1^{\zeta+1}}(B(i',\zeta+1)))}(\lambda_{k_1},\ldots,\lambda_{k_{\zeta+1}}) }{\tilde{\Psi}_{[\bm{X}]_1^{\zeta+1}}(B(i',\zeta+1))!} \times  \prod\limits_{\substack{\varsigma =\mbox{Ind}(\tilde{\Psi}_{[\bm{X}]_1^{\zeta+1}}(B(i',\zeta+1))), \bm{Z}_\varsigma=\bm{N}^{q_\varsigma}_{k_\varsigma,i_\varsigma}\bm{Y}_{\varsigma} \\ \varsigma \neq \mbox{Ind}(\tilde{\Psi}_{[\bm{X}]_1^{\zeta+1}}(B(i',\zeta+1))), \bm{Z}_\varsigma=\bm{P}_{k_\varsigma,i_\varsigma}\bm{Y}_{\varsigma}}
}^{\zeta+1} \bm{Z}_\varsigma \right\Vert \nonumber \\
&\leq& \sum\limits_{ k_\varsigma=1,  \mbox{for}~\varsigma=\mbox{Ind}(\tilde{\Psi}_{[\bm{X}]_1^{\zeta+1}}(B(i',\zeta+1)))}^{K_\varsigma, \mbox{for}~\varsigma=\mbox{Ind}(\tilde{\Psi}_{[\bm{X}]_1^{\zeta+1}}(B(i',\zeta+1)))} \sum\limits_{i_\varsigma=1,  \mbox{for}~\varsigma=\mbox{Ind}(\tilde{\Psi}_{[\bm{X}]_1^{\zeta+1}}(B(i',\zeta+1)))}^{\alpha_{k_\varsigma}^{(\mathrm{G})}, \mbox{for}~\varsigma=\mbox{Ind}(\tilde{\Psi}_{[\bm{X}]_1^{\zeta+1}}(B(i',\zeta+1)))} \nonumber \\
&& 
\sum\limits_{\tilde{\Psi}_{[\bm{X}]_1^{\zeta+1}}(B(i',\zeta+1))=1}^{m_{k_{\mbox{Ind}(\tilde{\Psi}_{[\bm{X}]_1^{\zeta+1}}(B(i',\zeta+1)))},i_{\mbox{Ind}(\tilde{\Psi}_{[\bm{X}]_1^{\zeta+1}}(B(i',\zeta+1)))}}-1} \nonumber \\
&& 
 \left[ \max\limits_{\lambda_{k_1} \in \Lambda_{\bm{X}_1},\ldots,\lambda_{k_{\zeta+1}} \in \Lambda_{\bm{X}_{\zeta+1}}} \left\vert \frac{\beta^{(\Psi_{[\bm{X}]_1^{\zeta+1}}(B(i',\zeta+1)))}(\lambda_{k_1},\ldots,\lambda_{k_{\zeta+1}}) }{\tilde{\Psi}_{[\bm{X}]_1^{\zeta+1}}(B(i',\zeta+1))!}\right\vert\right]  \nonumber \\
&& \times \left\Vert  \sum\limits_{ k_\varsigma=1,  \mbox{for}~\varsigma \neq \mbox{Ind}(\tilde{\Psi}_{[\bm{X}]_1^{\zeta+1}}(B(i',\zeta+1)))}^{K_\varsigma, \mbox{for}~\varsigma \neq \mbox{Ind}(\tilde{\Psi}_{[\bm{X}]_1^{\zeta+1}}(B(i',\zeta+1)))} \sum\limits_{i_\varsigma=1,  \mbox{for}~\varsigma \neq \mbox{Ind}(\tilde{\Psi}_{[\bm{X}]_1^{\zeta+1}}(B(i',\zeta+1)))}^{\alpha_{k_\varsigma}^{(\mathrm{G})}, \mbox{for}~\varsigma \neq \mbox{Ind}(\tilde{\Psi}_{[\bm{X}]_1^{\zeta+1}}(B(i',\zeta+1)))} \right. \nonumber \\
&&
\left.  \prod\limits_{\substack{\varsigma =\mbox{Ind}(\tilde{\Psi}_{[\bm{X}]_1^{\zeta+1}}(B(i',\zeta+1))), \bm{Z}_\varsigma= \bm{N}^{q_\varsigma}_{k_\varsigma,i_\varsigma}  \bm{Y}_{\varsigma} \\ \varsigma \neq \mbox{Ind}(\tilde{\Psi}_{[\bm{X}]_1^{\zeta+1}}(B(i',\zeta+1))), \bm{Z}_\varsigma=\bm{P}_{k_\varsigma,i_\varsigma}\bm{Y}_{\varsigma}}
}^{\zeta+1} \bm{Z}_\varsigma \right\Vert \nonumber \\
&\leq_1& \sum\limits_{ k_\varsigma=1,  \mbox{for}~\varsigma=\mbox{Ind}(\tilde{\Psi}_{[\bm{X}]_1^{\zeta+1}}(B(i',\zeta+1)))}^{K_\varsigma, \mbox{for}~\varsigma=\mbox{Ind}(\tilde{\Psi}_{[\bm{X}]_1^{\zeta+1}}(B(i',\zeta+1)))} \sum\limits_{i_\varsigma=1,  \mbox{for}~\varsigma=\mbox{Ind}(\tilde{\Psi}_{[\bm{X}]_1^{\zeta+1}}(B(i',\zeta+1)))}^{\alpha_{k_\varsigma}^{(\mathrm{G})}, \mbox{for}~\varsigma=\mbox{Ind}(\tilde{\Psi}_{[\bm{X}]_1^{\zeta+1}}(B(i',\zeta+1)))}\nonumber \\
&& 
\sum\limits_{\tilde{\Psi}_{[\bm{X}]_1^{\zeta+1}}(B(i',\zeta+1))=1}^{m_{k_{\mbox{Ind}(\tilde{\Psi}_{[\bm{X}]_1^{\zeta+1}}(B(i',\zeta+1)))},i_{\mbox{Ind}(\tilde{\Psi}_{[\bm{X}]_1^{\zeta+1}}(B(i',\zeta+1)))}}-1}\nonumber \\
&& 
\left[ \max\limits_{\lambda_{k_1} \in \Lambda_{\bm{X}_1},\ldots,\lambda_{k_{\zeta+1}} \in \Lambda_{\bm{X}_{\zeta+1}}} \left\vert \frac{\beta^{(\Psi_{[\bm{X}]_1^{\zeta+1}}(B(i',\zeta+1)))}(\lambda_{k_1},\ldots,\lambda_{k_{\zeta+1}}) }{\tilde{\Psi}_{[\bm{X}]_1^{\zeta+1}}(B(i',\zeta+1))!}\right\vert\right] \nonumber \\
&& \times  \prod\limits_{\substack{\varsigma =\mbox{Ind}(\tilde{\Psi}_{[\bm{X}]_1^{\zeta+1}}(B(i',\zeta+1))), \bm{Z}_\varsigma=\left\Vert \bm{N}^{q_\varsigma}_{k_\varsigma,i_\varsigma} \right\Vert \left\Vert \bm{Y}_{\varsigma}\right\Vert\\ \varsigma \neq \mbox{Ind}(\tilde{\Psi}_{[\bm{X}]_1^{\zeta+1}}(B(i',\zeta+1))), \bm{Z}_\varsigma=\left\Vert\bm{Y}_{\varsigma}\right\Vert}
}^{\zeta+1} \bm{Z}_\varsigma,
\end{eqnarray}
}
where we apply $\sum\limits_{k_\varsigma=1}^{K_\varsigma}\sum\limits_{i_\varsigma=1}^{\alpha_{k_\varsigma}^{(\mathrm{G})}}\bm{P}_{k_\varsigma,i_\varsigma}= \bm{I}$ for $\varsigma \neq \mbox{Ind}(\tilde{\Psi}_{[\bm{X}]_1^{\zeta+1}}(B(i',\zeta+1)))$, and the triangle inequality of Frobenius norm again in $\leq_1$. Therefore, we have the upper bound for $T_{\beta}^{\bm{X}_1,\ldots,\bm{X}_{\zeta+1}}(\bm{Y}_1,\ldots,\bm{Y}_{\zeta})$ as the upper bound shown by Eq.~\eqref{eq6:thm: GMOI norm Est} is identical to $\left\Vert \bm{A}_{i'} \right\Vert_{up}$. 

For the lower bound of $T_{\beta}^{\bm{X}_1,\ldots,\bm{X}_{\zeta+1}}(\bm{Y}_1,\ldots,\bm{Y}_{\zeta})$, we have Eq.~\eqref{eq2:  thm: GMOI norm Est} immediatedly from Lemma~\ref{lma:conv triangle for Frob norm}.

If we have $\left[\min\limits_{\lambda_1 \in \Lambda_{\bm{X}_1},\ldots,\lambda_{\zeta+1} \in \Lambda_{\bm{X}_{\zeta+1}}} \left\vert\beta(\lambda_1,\dots, \lambda_{\zeta+1})\right\vert\right]\left\Vert \prod\limits_{i'=1}^{\zeta}\bm{Y}_{i'}\right\Vert\geq \sum\limits_{i'=1}^{2^{\zeta+1}-1}\left\Vert\bm{A}_{i'}\right\Vert$ and Lemma~\ref{lma:conv triangle for Frob norm}, we have
\begin{eqnarray}\label{eq5:  thm: GMOI norm Est}
\lefteqn{\left\Vert T_{\beta}^{\bm{X}_1,\ldots,\bm{X}_{\zeta+1}}(\bm{Y}_1,\ldots,\bm{Y}_{\zeta}) \right\Vert} \nonumber \\
&\geq&\left\Vert\bm{A}_0\right\Vert-(\left\Vert\bm{A}_{1}\right\Vert+\left\Vert\bm{A}_{2}\right\Vert + \ldots + \left\Vert\bm{A}_{2^{\zeta+1}-1}\right\Vert) \nonumber \\
&\geq&\left[\min\limits_{\lambda_1 \in \Lambda_{\bm{X}_1},\ldots,\lambda_{\zeta+1} \in \Lambda_{\bm{X}_{\zeta+1}}} \left\vert\beta(\lambda_1,\dots, \lambda_{\zeta+1})\right\vert\right] \nonumber \\
&& \times \left\Vert \sum\limits_{k_1=\ldots=k_{\zeta+1}=1}^{K_1,\ldots,K_{\zeta+1}} \sum\limits_{i_1=\ldots=i_{\zeta+1}=1}^{\alpha_{k_1}^{(\mathrm{G})},\ldots,\alpha_{k_{\zeta+1}}^{(\mathrm{G})}}
\bm{P}_{k_1,i_1}\bm{Y}_1\bm{P}_{k_2,i_2}\bm{Y}_2\ldots\bm{Y}_{\zeta}\bm{P}_{k_{\zeta+1},i_{\zeta+1}} \right\Vert\nonumber \\
& &-  (\left\Vert\bm{A}_{1}\right\Vert+\left\Vert\bm{A}_{2}\right\Vert + \ldots + \left\Vert\bm{A}_{2^{\zeta+1}-1}\right\Vert)\nonumber \\
&=&\left[\min\limits_{\lambda_1 \in \Lambda_{\bm{X}_1},\ldots,\lambda_{\zeta+1} \in \Lambda_{\bm{X}_{\zeta+1}}} \left\vert\beta(\lambda_1,\dots, \lambda_{\zeta+1})\right\vert\right] \left\Vert \prod\limits_{i'=1}^{\zeta}\bm{Y}_{i'}\right\Vert \nonumber \\
&& - (\left\Vert\bm{A}_{1}\right\Vert+\left\Vert\bm{A}_{2}\right\Vert + \ldots + \left\Vert\bm{A}_{2^{\zeta+1}-1}\right\Vert),
\end{eqnarray}
which is the lower bound of $T_{\beta}^{\bm{X}_1,\ldots,\bm{X}_{\zeta+1}}(\bm{Y}_1,\ldots,\bm{Y}_{\zeta})$ given by Eq.~\eqref{eq3:  thm: GMOI norm Est}.
$\hfill\Box$

\section{Perturbation Formula and Lipschitz Estimations}\label{sec:Perturbation Formula and Lipschitz Estimation}

\subsection{Perturbation Formula and Nilpotent Part Difference Characterization}\label{sec: Perturbation Formula and Nilpotant Part Difference Characterization}

In this section, the perturbation formula of GMOI is derived.  We begin by defining several new notations for later perturbation formula simpler presentation. We have 

\begin{eqnarray}\label{eq1: pert nimpler notation}
\bm{S}_{k_c,i_c}&=&\begin{cases}
\bm{P}_{k_c,i_c}\\
\bm{N}_{k_c,i_c}^{q_c}
\end{cases};\nonumber \\
\bm{S}_{k_d,i_d}&=&\begin{cases}
\bm{P}_{k_d,i_d}\\
\bm{N}_{k_d,i_d}^{q_d}
\end{cases}.
\end{eqnarray}
where $c,d$ are indices for the matrices $\bm{C}$ and $\bm{D}$. Also, we set 
\begin{eqnarray}\label{eq2: pert nimpler notation}
\bm{S}_{k_p,i_p}&=&\begin{cases}
\bm{P}_{k_p,i_p}\\
\bm{N}_{k_p,i_p}^{q_p}
\end{cases};\nonumber \\
\bm{S}_{k_{p'},i_{p'}}&=&\begin{cases}
\bm{P}_{k_{p'},i_{p'}}\\
\bm{N}_{k_{p'},i_{p'}}^{q_{p'}}
\end{cases},
\end{eqnarray}
where $p,p' \in \mathbb{N}$ are indices for the matrices $\bm{X}_p$ and $\bm{X}_{p'}$. We use the following summation notation:
\begin{eqnarray}\label{eq3: pert nimpler notation}
\sum\limits_{\bm{S}_{k_p,i_p}} &=&
\begin{cases}
\sum\limits_{k_p=1}^{K_p}\sum\limits_{i_p=1}^{\alpha_{k_p}^{(\mathrm{G})}}\bm{P}_{k_p,i_p}, \mbox{~~if $\bm{S}_{k_p,i_p}=\bm{P}_{k_p,i_p}$}, \\
\sum\limits_{k_p=1}^{K_p}\sum\limits_{i_p=1}^{\alpha_{k_p}^{(\mathrm{G})}}\sum\limits_{q_p=1}^{m_{k_p,i_p}}\bm{N}_{k_p,i_p}^{q_p}, \mbox{~~if $\bm{S}_{k_p,i_p}=\bm{N}_{k_p,i_p}^{q_p}$}, 
\end{cases}
\end{eqnarray}
therefore, if we write $\sum\limits_{\bm{S}_{k_p,i_p}, p \in [p_1, p_2]}$, which can be expressed by
\begin{eqnarray}\label{eq4: pert nimpler notation}
\sum\limits_{\bm{S}_{k_p,i_p}, p \in [p_1, p_2]} &=&
\sum\limits_{\bm{S}_{k_{p_1},i_{p_1}}}\sum\limits_{\bm{S}_{k_{p_2},i_{p_2}}}\ldots\sum\limits_{\bm{S}_{k_{p_2},i_{p_2}}}
\end{eqnarray}
where $[p_1, p_2]$ are positive integers from $p_1, p_1 + 1, \ldots, p_2$. Note that the expression given by Eq.~\eqref{eq4: pert nimpler notation} will include all summation results for all binary combinations of $\bm{S}_{k_p,i_p}$. Hence, there are $2^{p_2 - p_1 +1}$ summation results given by Eq.~\eqref{eq4: pert nimpler notation}.

Given any differentiable function $f: \mathbb{C}^{\zeta} \rightarrow \mathbb{C}$, we use the following notation to represent partial derivatives with respect to different arguments of the function $f$ and its normalization with respect to differentiation orders:
\begin{eqnarray}\label{eq5: pert nimpler notation}
\frac{f^{(q_1,q_2,\ldots,q_\zeta)}(\lambda_1, \lambda_2, \ldots, \lambda_{\zeta})}{q_1 ! q_2 ! \ldots,q_\zeta! }&=& \frac{f^{([q]_1^{q_\zeta})}([\lambda]_1^{\zeta})}{([q!]_1^{q_\zeta})}
\end{eqnarray}

The Lemma~\ref{lma:first-order divided difference identity} will be provided to show that the $k$-order divided difference identity will be valid for the function and its partial derivatives.
\begin{lemma}\label{lma:first-order divided difference identity}
Let \( f : \mathbb{C}^n \to \mathbb{C} \) be a function of class \( C^{k+1+r} \). Let
\[
D^{(r)} f = \frac{\partial^r f}{\partial x_{i_1} \cdots \partial x_{i_r}}
\]
be an \( r \)-th order partial derivative of \( f \). Then the divided difference identity
\[
D^{(r)} f^{[k]}(\dots, \lambda_{\alpha_j}, \dots) - D^{(r)} f^{[k]}(\dots, \mu_{\alpha_j}, \dots) = (\lambda_{\alpha_j} - \mu_{\alpha_j}) D^{(r)} f^{[k+1]}(\dots, \lambda_{\alpha_j}, \mu_{\alpha_j}, \dots)
\]
holds provided that the divided difference is taken with respect to a variable independent of \( x_{i_1}, \dots, x_{i_r} \).
\end{lemma}
\textbf{Proof:}
Assume \( f \in C^{k+1+r}(\mathbb{C}^n) \), so that all relevant derivatives and divided differences are well-defined and smooth.

Let \( g = D^{(r)} f = \frac{\partial^r f}{\partial x_{i_1} \cdots \partial x_{i_r}} \). Because differentiation and divided difference operators commute when they act on independent variables, we have:
\[
D^{(r)} f^{[k]}(\lambda_1, \dots, \lambda_k) = \left( D^{(r)} f \right)^{[k]}(\lambda_1, \dots, \lambda_k)
\]
Now, apply the standard first-order divided difference identity to the function \( g \):
\[
g^{[k]}(\dots, \lambda_{\alpha_j}, \dots) - g^{[k]}(\dots, \mu_{\alpha_j}, \dots)
= (\lambda_{\alpha_j} - \mu_{\alpha_j}) g^{[k+1]}(\dots, \lambda_{\alpha_j}, \mu_{\alpha_j}, \dots)
\]
Substituting back \( g = D^{(r)} f \), we obtain:
\[
D^{(r)} f^{[k]}(\dots, \lambda_{\alpha_j}, \dots) - D^{(r)} f^{[k]}(\dots, \mu_{\alpha_j}, \dots)
= (\lambda_{\alpha_j} - \mu_{\alpha_j}) D^{(r)} f^{[k+1]}(\dots, \lambda_{\alpha_j}, \mu_{\alpha_j}, \dots)
\]
which completes the proof.
$\hfill\Box$

We are ready to present Theorem~\ref{thm:GMOI Perturbation Formula} about GMOI perturbation formula. 
\begin{theorem}\label{thm:GMOI Perturbation Formula}
Given matrices $\bm{X}_p$ decomposed as:
\begin{eqnarray}\label{eq1:thm:GMOI Perturbation Formula}
\bm{X}_p&=&\sum\limits_{k_p=1}^{K_p}\sum\limits_{i_p=1}^{\alpha_{k_p}^{\mathrm{G}}} \lambda_{k_p} \bm{P}_{k_p,i_p}+
\sum\limits_{k_p=1}^{K_p}\sum\limits_{i_p=1}^{\alpha_{k_p}^{\mathrm{G}}} \bm{N}_{k_p,i_p},
\end{eqnarray}
where $p=1,2,\ldots,\zeta$, and matrices $\bm{C}, \bm{D}$ decomposed as:
\begin{eqnarray}\label{eq1:thm:GMOI Perturbation Formula}
\bm{C}&=&\sum\limits_{k_c=1}^{K_c}\sum\limits_{i_c=1}^{\alpha_{k_c}^{\mathrm{G}}} \lambda_{k_c} \bm{P}_{k_c,i_c}+
\sum\limits_{k_c=1}^{K_c}\sum\limits_{i_c=1}^{\alpha_{k_c}^{\mathrm{G}}} \bm{N}_{k_c,i_c}, \nonumber \\ 
\bm{D}&=&\sum\limits_{k_d=1}^{K_d}\sum\limits_{i_d=1}^{\alpha_{k_d}^{\mathrm{G}}} \lambda_{k_d} \bm{P}_{k_d,i_d}+
\sum\limits_{k_d=1}^{K_d}\sum\limits_{i_d=1}^{\alpha_{k_d}^{\mathrm{G}}} \bm{N}_{k_d,i_d}.
\end{eqnarray}

We have the following formula for GMOI:
{\tiny
\begin{eqnarray}\label{eq1:thm:GMOI Perturbation Formula}
\lefteqn{T_{\beta^{[\zeta+1]}}^{[\bm{X}]_1^{j-1},\bm{C},\bm{D},[\bm{X}]_j^{\zeta}}([\bm{Y}]_1^{j-1},\bm{C}-\bm{D},[\bm{Y}]_j^{\zeta})=}\nonumber \\
&& \Bigg(T_{\beta^{[\zeta]}}^{[\bm{X}]_1^{j-2},\bm{X}_{j-1,P},\bm{C}_P,\bm{X}_{j,P},[\bm{X}]_{j+1}^{\zeta}}([\bm{Y}]_1^{j-1},[\bm{Y}]_j^{\zeta}) \nonumber \\
&-& T_{\beta^{[\zeta]}}^{[\bm{X}]_1^{j-2},\bm{X}_{j-1,P},\bm{D}_P,\bm{X}_{j,P},[\bm{X}]_{j+1}^{\zeta}}([\bm{Y}]_1^{j-1},[\bm{Y}]_j^{\zeta})\Bigg)_I + \mathfrak{X}^{\bm{X}_{j-1,P},\bm{C}_P,\bm{D}_P,\bm{X}_{j,P}} \nonumber \\
&+&\Bigg( T_{\beta^{[\zeta]}}^{[\bm{X}]_1^{j-2},\bm{X}_{j-1,P},\bm{C}_N,\bm{X}_{j,P},[\bm{X}]_{j+1}^{\zeta}}([\bm{Y}]_1^{j-1},[\bm{Y}]_j^{\zeta} \nonumber \\
&-& T_{\beta^{[\zeta]}}^{[\bm{X}]_1^{j-2},\bm{X}_{j-1,P},\bm{D}_N,\bm{X}_{j,P},[\bm{X}]_{j+1}^{\zeta}}([\bm{Y}]_1^{j-1},[\bm{Y}]_j^{\zeta}))\Bigg)_{II} + \mathfrak{X}^{\bm{X}_{j-1,P},\bm{C}_N,\bm{D}_N,\bm{X}_{j,P}} \nonumber \\
&+& \Bigg( T_{\beta^{[\zeta]}}^{[\bm{X}]_1^{j-2},\bm{X}_{j-1,P},\bm{C}_P,\bm{X}_{j,N},[\bm{X}]_{j+1}^{\zeta}}([\bm{Y}]_1^{j-1},\bm{C},[\bm{Y}]_j^{\zeta}) \nonumber \\
&-& T_{\beta^{[\zeta]}}^{[\bm{X}]_1^{j-2},\bm{X}_{j-1,P},\bm{D}_P,\bm{X}_{j,N},[\bm{X}]_{j+1}^{\zeta}}([\bm{Y}]_1^{j-1},\bm{D},[\bm{Y}]_j^{\zeta})\Bigg)_{III}+ \mathfrak{X}^{\bm{X}_{j-1,P},\bm{C}_P,\bm{D}_P,\bm{X}_{j,N}} \nonumber \\
&+&\Bigg( T_{\beta^{[\zeta]}}^{[\bm{X}]_1^{j-2},\bm{X}_{j-1,P},\bm{C}_N,\bm{X}_{j,N},[\bm{X}]_{j+1}^{\zeta}}([\bm{Y}]_1^{j-1},[\bm{Y}]_j^{\zeta}) \nonumber \\
&-& T_{\beta^{[\zeta]}}^{[\bm{X}]_1^{j-2},\bm{X}_{j-1,P},\bm{D}_N,\bm{X}_{j,N},[\bm{X}]_{j+1}^{\zeta}}([\bm{Y}]_1^{j-1},[\bm{Y}]_j^{\zeta})\Bigg)_{IV}+ \mathfrak{X}^{\bm{X}_{j-1,P},\bm{C}_N,\bm{D}_N,\bm{X}_{j,N}} \nonumber \\
&+& \Bigg( T_{\beta^{[\zeta]}}^{[\bm{X}]_1^{j-2},\bm{X}_{j-1,N},\bm{C}_P,\bm{X}_{j,P},[\bm{X}]_{j+1}^{\zeta}}([\bm{Y}]_1^{j-1},[\bm{Y}]_j^{\zeta}) \nonumber \\
&-& T_{\beta^{[\zeta]}}^{[\bm{X}]_1^{j-2},\bm{X}_{j-1,N},\bm{D}_P,\bm{X}_{j,P},[\bm{X}]_{j+1}^{\zeta}}([\bm{Y}]_1^{j-1},[\bm{Y}]_j^{\zeta})\Bigg)_{V}+ \mathfrak{X}^{\bm{X}_{j-1,N},\bm{C}_P,\bm{D}_P,\bm{X}_{j,P}} \nonumber \\
&+&\Bigg(T_{\beta^{[\zeta]}}^{[\bm{X}]_1^{j-2},\bm{X}_{j-1,N},\bm{C}_N,\bm{X}_{j,P},[\bm{X}]_{j+1}^{\zeta}}([\bm{Y}]_1^{j-1},[\bm{Y}]_j^{\zeta}) \nonumber \\
&-& T_{\beta^{[\zeta]}}^{[\bm{X}]_1^{j-2},\bm{X}_{j-1,N},\bm{D}_N,\bm{X}_{j,P},[\bm{X}]_{j+1}^{\zeta}}([\bm{Y}]_1^{j-1},[\bm{Y}]_j^{\zeta})\Bigg)_{VI} + \mathfrak{X}^{\bm{X}_{j-1,N},\bm{C}_N,\bm{D}_N,\bm{X}_{j,P}} \nonumber \\
&+& \Bigg(T_{\beta^{[\zeta]}}^{[\bm{X}]_1^{j-2},\bm{X}_{j-1,N},\bm{C}_P,\bm{X}_{j,N},[\bm{X}]_{j+1}^{\zeta}}([\bm{Y}]_1^{j-1},[\bm{Y}]_j^{\zeta}) \nonumber \\
&-& T_{\beta^{[\zeta]}}^{[\bm{X}]_1^{j-2},\bm{X}_{j-1,N},\bm{D}_P,\bm{X}_{j,N},[\bm{X}]_{j+1}^{\zeta}}([\bm{Y}]_1^{j-1},[\bm{Y}]_j^{\zeta})\Bigg)_{VII} + \mathfrak{X}^{\bm{X}_{j-1,N},\bm{C}_P,\bm{D}_P,\bm{X}_{j,N}} \nonumber \\
&+&\Bigg(T_{\beta^{[\zeta]}}^{[\bm{X}]_1^{j-2},\bm{X}_{j-1,N},\bm{C}_N,\bm{X}_{j,N},[\bm{X}]_{j+1}^{\zeta}}([\bm{Y}]_1^{j-1},[\bm{Y}]_j^{\zeta}) \nonumber \\
&-& T_{\beta^{[\zeta]}}^{[\bm{X}]_1^{j-2},\bm{X}_{j-1,N},\bm{D}_N,\bm{X}_{j,N},[\bm{X}]_{j+1}^{\zeta}}([\bm{Y}]_1^{j-1},[\bm{Y}]_j^{\zeta})\Bigg)_{VIII} + \mathfrak{X}^{\bm{X}_{j-1,N},\bm{C}_N,\bm{D}_N,\bm{X}_{j,N}}\nonumber \\
&+& T_{\beta^{[\zeta+1]}}^{[\bm{X}]_1^{j-2},\bm{X}_{j-1,P},\bm{C}_N,\bm{D}_N,\bm{X}_{j,P},[\bm{X}]_{j+1}^{\zeta}}([\bm{Y}]_1^{j-1},\bm{C}-\bm{D},[\bm{Y}]_j^{\zeta})\nonumber \\
&+& T_{\beta^{[\zeta+1]}}^{[\bm{X}]_1^{j-2},\bm{X}_{j-1,P},\bm{C}_N,\bm{D}_N,\bm{X}_{j,N},[\bm{X}]_{j+1}^{\zeta}}([\bm{Y}]_1^{j-1},\bm{C}-\bm{D},[\bm{Y}]_j^{\zeta})\nonumber \\
&+& T_{\beta^{[\zeta+1]}}^{[\bm{X}]_1^{j-2},\bm{X}_{j-1,N},\bm{C}_N,\bm{D}_N,\bm{X}_{j,P},[\bm{X}]_{j+1}^{\zeta}}([\bm{Y}]_1^{j-1},\bm{C}-\bm{D},[\bm{Y}]_j^{\zeta})\nonumber \\
&+& T_{\beta^{[\zeta+1]}}^{[\bm{X}]_1^{j-2},\bm{X}_{j-1,N},\bm{C}_N,\bm{D}_N,\bm{X}_{j,N},[\bm{X}]_{j+1}^{\zeta}}([\bm{Y}]_1^{j-1},\bm{C}-\bm{D},[\bm{Y}]_j^{\zeta})\nonumber \\
&=&T_{\beta^{[\zeta]}}^{[\bm{X}]_1^{j-1},\bm{C},[\bm{X}]_{j}^{\zeta}}([\bm{Y}]_1^{j-1},[\bm{Y}]_j^{\zeta}) - T_{\beta^{[\zeta]}}^{[\bm{X}]_1^{j-1},\bm{D},[\bm{X}]_{j}^{\zeta}}([\bm{Y}]_1^{j-1},[\bm{Y}]_j^{\zeta}) + \bar{\mathfrak{X}}([\bm{X}]_1^{j-1},\bm{C},\bm{D},[\bm{X}]_j^{\zeta}) \nonumber \\
&+& T_{\beta^{[\zeta+1]}}^{[\bm{X}]_1^{j-2},\bm{X}_{j-1,P},\bm{C}_N,\bm{D}_N,\bm{X}_{j,P},[\bm{X}]_{j+1}^{\zeta}}([\bm{Y}]_1^{j-1},\bm{C}-\bm{D},[\bm{Y}]_j^{\zeta})\nonumber \\
&+& T_{\beta^{[\zeta+1]}}^{[\bm{X}]_1^{j-2},\bm{X}_{j-1,P},\bm{C}_N,\bm{D}_N,\bm{X}_{j,N},[\bm{X}]_{j+1}^{\zeta}}([\bm{Y}]_1^{j-1},\bm{C}-\bm{D},[\bm{Y}]_j^{\zeta})\nonumber \\
&+& T_{\beta^{[\zeta+1]}}^{[\bm{X}]_1^{j-2},\bm{X}_{j-1,N},\bm{C}_N,\bm{D}_N,\bm{X}_{j,P},[\bm{X}]_{j+1}^{\zeta}}([\bm{Y}]_1^{j-1},\bm{C}-\bm{D},[\bm{Y}]_j^{\zeta})\nonumber \\
&+& T_{\beta^{[\zeta+1]}}^{[\bm{X}]_1^{j-2},\bm{X}_{j-1,N},\bm{C}_N,\bm{D}_N,\bm{X}_{j,N},[\bm{X}]_{j+1}^{\zeta}}([\bm{Y}]_1^{j-1},\bm{C}-\bm{D},[\bm{Y}]_j^{\zeta})
\end{eqnarray}
}
where $\beta^{[\ell]}$ is the $\ell$-th order divided difference of a differentiable function $\beta$,  $2 \leq \zeta \in \mathbb{N}$, and $\mathfrak{X}^{\bm{Z}_1,\bm{Z}_2,\bm{Z}_3,\bm{Z}_4}$ are perturbation correction terms with respect to 
\begin{eqnarray}\label{eq2:thm:GMOI Perturbation Formula}
\lefteqn{T_{\beta^{[\zeta+1]}}^{[\bm{X}]_1^{j-2},\bm{Z}_1, \bm{Z}_2,\bm{Z}_3,\bm{Z}_4, [\bm{X}]_{j+1}^{\zeta}}([\bm{Y}]_1^{j-1},\bm{C}-\bm{D},[\bm{Y}]_j^{\zeta})=}\nonumber \\
&&T_{\beta^{[\zeta]}}^{[\bm{X}]_1^{j-2},\bm{Z}_1,\bm{Z}_2,\bm{Z}_4,[\bm{X}]_{j+1}^{\zeta}}([\bm{Y}]_1^{j-1},[\bm{Y}]_j^{\zeta}) \nonumber \\
&&-T_{\beta^{[\zeta]}}^{[\bm{X}]_1^{j-2},\bm{Z}_1,\bm{Z}_3,\bm{Z}_4,[\bm{X}]_{j+1}^{\zeta}}([\bm{Y}]_1^{j-1},[\bm{Y}]_j^{\zeta})+ \mathfrak{X}^{\bm{Z}_1,\bm{Z}_2,\bm{Z}_3,\bm{Z}_4},
\end{eqnarray}
and
\begin{eqnarray}\label{eq2:thm:GMOI Perturbation Formula}
\bar{\mathfrak{X}}([\bm{X}]_1^{j-1},\bm{C},\bm{D},[\bm{X}]_j^{\zeta})&=& \mathfrak{X}^{\bm{X}_{j-1,P},\bm{C}_P,\bm{D}_P,\bm{X}_{j,P}} + \mathfrak{X}^{\bm{X}_{j-1,P},\bm{C}_N,\bm{D}_N,\bm{X}_{j,P}} \nonumber \\
&+& \mathfrak{X}^{\bm{X}_{j-1,P},\bm{C}_P,\bm{D}_P,\bm{X}_{j,N}} + \mathfrak{X}^{\bm{X}_{j-1,P},\bm{C}_N,\bm{D}_N,\bm{X}_{j,N}} \nonumber \\
&+&  \mathfrak{X}^{\bm{X}_{j-1,N},\bm{C}_P,\bm{D}_P,\bm{X}_{j,P}} + \mathfrak{X}^{\bm{X}_{j-1,N},\bm{C}_N,\bm{D}_N,\bm{X}_{j,P}} \nonumber \\
&+& \mathfrak{X}^{\bm{X}_{j-1,N},\bm{C}_P,\bm{D}_P,\bm{X}_{j,N}} + \mathfrak{X}^{\bm{X}_{j-1,N},\bm{C}_N,\bm{D}_N,\bm{X}_{j,N}}. 
\end{eqnarray}

Each term related to $\mathfrak{X}^{\bm{Z}_1,\bm{Z}_2,\bm{Z}_3,\bm{Z}_4}$ can be expressed by
\begin{eqnarray}\label{eq1:cross:thm:GMOI Perturbation Formula}
\lefteqn{\mathfrak{X}^{\bm{X}_{j-1,P},\bm{C}_P,\bm{D}_P,\bm{X}_{j,P}}=}\nonumber \\
&&\sum\limits_{\substack{S_{k_p,i_p}, p \in [1,j-2], \\
S_{k_{p'},i_{p'}}, p' \in [j+1,\zeta]}}\sum\limits_{\substack{k_{j-1},i_{j-1},\\k_c, i_c,q_c,k_d,i_d,q_d\\ k_{j},i_{j}}}\frac{(\beta^{[\zeta+1]}([\lambda_k]_1^{j-2},\lambda_{k_{j-1}},\lambda_{k_c},\lambda_{k_d},\lambda_{k_j},[\lambda_k]_{j+1}^{\zeta}))^{([q]_1^{j-2}[0,0,0,0][q]_{j+1}^{\zeta})}}{ ([q!]_1^{j-2}[0!,0!,0!,0!][q!]_{j+1}^{\zeta})  } \nonumber \\
&& \times \left(\prod\limits_{p=1}^{j-2}S_{k_p,i_p}\bm{Y}_p\right)\bm{P}_{k_{j-1},i_{j-1}}\bm{Y}_{j-1}\bm{P}_{k_c,i_c}(\bm{N}_{k_c,i_c}^{q_c}-\bm{N}_{k_d,i_d}^{q_d})\bm{P}_{k_d,i_d}\bm{Y}_j\bm{P}_{k_{j},i_{j}}\left(\prod\limits_{p'=j+1}^{\zeta}\bm{Y}_{p'}S_{k_{p'},i_{p'}}\right),
\end{eqnarray}
\begin{eqnarray}\label{eq2:cross:thm:GMOI Perturbation Formula}
\lefteqn{\mathfrak{X}^{\bm{X}_{j-1,P},\bm{C}_N,\bm{D}_N,\bm{X}_{j,P}}=}\nonumber \\
&&\sum\limits_{\substack{S_{k_p,i_p}, p \in [1,j-2], \\
S_{k_{p'},i_{p'}}, p' \in [j+1,\zeta]}}\sum\limits_{\substack{k_{j-1},i_{j-1},\\k_c, i_c,q_c,k_d,i_d,q_d\\ k_{j},i_{j}}}\frac{(\beta^{[\zeta+1]}([\lambda_k]_1^{j-2},\lambda_{k_{j-1}},\lambda_{k_c},\lambda_{k_d},\lambda_{k_j},[\lambda_k]_{j+1}^{\zeta}))^{([q]_1^{j-2}[0,0,q_d,0][q]_{j+1}^{\zeta})}}{ ([q!]_1^{j-2}[0!,0!,q_d!,0!][q!]_{j+1}^{\zeta})  } \nonumber \\
&& \times \left(\prod\limits_{p=1}^{j-2}S_{k_p,i_p}\bm{Y}_p\right)\bm{P}_{k_{j-1},i_{j-1}}\bm{Y}_{j-1}\bm{P}_{k_c,i_c}(\bm{N}_{k_c,i_c}^{q_c}-\bm{N}_{k_d,i_d}^{q_d})\bm{N}_{k_d,i_d}^{q_d}\bm{Y}_j\bm{P}_{k_{j},i_{j}}\left(\prod\limits_{p'=j+1}^{\zeta}\bm{Y}_{p'}S_{k_{p'},i_{p'}}\right) \nonumber \\
&+& \sum\limits_{\substack{S_{k_p,i_p}, p \in [1,j-2], \\
S_{k_{p'},i_{p'}}, p' \in [j+1,\zeta]}}\sum\limits_{\substack{k_{j-1},i_{j-1},\\k_c, i_c,q_c,k_d,i_d,q_d\\ k_{j},i_{j}}}\frac{(\beta^{[\zeta+1]}([\lambda_k]_1^{j-2},\lambda_{k_{j-1}},\lambda_{k_c},\lambda_{k_d},\lambda_{k_j},[\lambda_k]_{j+1}^{\zeta}))^{([q]_1^{j-2}[0,q_c,0,0][q]_{j+1}^{\zeta})}}{ ([q!]_1^{j-2}[0!,q_c!,0!,0!][q!]_{j+1}^{\zeta})  } \nonumber \\
&& \times \left(\prod\limits_{p=1}^{j-2}S_{k_p,i_p}\bm{Y}_p\right)\bm{P}_{k_{j-1},i_{j-1}}\bm{Y}_{j-1}\bm{N}_{k_c,i_c}^{q_c}(\bm{N}_{k_c,i_c}^{q_c}-\bm{N}_{k_d,i_d}^{q_d})\bm{P}_{k_d,i_d}\bm{Y}_j\bm{P}_{k_{j},i_{j}}\left(\prod\limits_{p'=j+1}^{\zeta}\bm{Y}_{p'}S_{k_{p'},i_{p'}}\right) \nonumber \\
&+& \sum\limits_{\substack{S_{k_p,i_p}, p \in [1,j-2], \\
S_{k_{p'},i_{p'}}, p' \in [j+1,\zeta]}}\sum\limits_{\substack{k_{j-1},i_{j-1},\\k_c,i_c, k_d,i_d,q_d\\ k_{j},i_{j}}} \frac{(\beta^{[\zeta]}([\lambda_k]_1^{j-2},\lambda_{k_{j-1}},\lambda_{k_c},\lambda_{k_j},[\lambda_k]_{j+1}^{\zeta}))^{([q]_1^{j-2}[0,q_d,0][q]_{j+1}^{\zeta})}}{ ([q!]_1^{j-2}[0!,q_d!,0!][q!]_{j+1}^{\zeta})  } \nonumber \\
&& \times \left(\prod\limits_{p=1}^{j-2}S_{k_p,i_p}\bm{Y}_p\right)\bm{P}_{k_{j-1},i_{j-1}}\bm{Y}_{j-1}\bm{P}_{k_c,i_c}\bm{N}_{k_d,i_d}^{q_d}\bm{Y}_j\bm{P}_{k_{j},i_{j}}\left(\prod\limits_{p'=j+1}^{\zeta}\bm{Y}_{p'}S_{k_{p'},i_{p'}}\right) \nonumber \\  
&-& \sum\limits_{\substack{S_{k_p,i_p}, p \in [1,j-2], \\
S_{k_{p'},i_{p'}}, p' \in [j+1,\zeta]}}\sum\limits_{\substack{k_{j-1},i_{j-1},\\k_c, i_c, q_c, k_d,i_d\\ k_{j},i_{j}}} \frac{(\beta^{[\zeta]}([\lambda_k]_1^{j-2},\lambda_{k_{j-1}},\lambda_{k_d},\lambda_{k_j},[\lambda_k]_{j+1}^{\zeta}))^{([q]_1^{j-2}[0,q_c,0][q]_{j+1}^{\zeta})}}{ ([q!]_1^{j-2}[0!,q_c!,0!][q!]_{j+1}^{\zeta})  } \nonumber \\
&& \times \left(\prod\limits_{p=1}^{j-2}S_{k_p,i_p}\bm{Y}_p\right)\bm{P}_{k_{j-1},i_{j-1}}\bm{Y}_{j-1}\bm{N}_{k_c,i_c}^{q_c}\bm{P}_{k_d,i_d}\bm{Y}_j\bm{P}_{k_{j},i_{j}}\left(\prod\limits_{p'=j+1}^{\zeta}\bm{Y}_{p'}S_{k_{p'},i_{p'}}\right),
\end{eqnarray}
\begin{eqnarray}\label{eq3:cross:thm:GMOI Perturbation Formula}
\lefteqn{\mathfrak{X}^{\bm{X}_{j-1,P},\bm{C}_P,\bm{D}_P,\bm{X}_{j,N}}=}\nonumber \\
&& \sum\limits_{\substack{S_{k_p,i_p}, p \in [1,j-2], \\
S_{k_{p'},i_{p'}}, p' \in [j+1,\zeta]}}\sum\limits_{\substack{k_{j-1},i_{j-1}\\k_c, i_c,q_c,k_d,i_d,q_d\\ k_{j},i_{j},q_j}}\frac{(\beta^{[\zeta+1]}([\lambda_k]_1^{j-2},\lambda_{k_{j-1}},\lambda_{k_c},\lambda_{k_d},\lambda_{k_j},[\lambda_k]_{j+1}^{\zeta}))^{([q]_1^{j-2}[0,0,0,q_j][q]_{j+1}^{\zeta})}}{ ([q!]_1^{j-2}[0!,0!,0!,q_j !][q!]_{j+1}^{\zeta})  } \nonumber \\
&& \times \left(\prod\limits_{p=1}^{j-2}S_{k_p,i_p}\bm{Y}_p\right)\bm{P}_{k_{j-1},i_{j-1}}\bm{Y}_{j-1}\bm{P}_{k_c,i_c}(\bm{N}_{k_c,i_c}^{q_c}-\bm{N}_{k_d,i_d}^{q_d})\bm{P}_{k_d,i_d}\bm{Y}_j\bm{N}_{k_{j},i_{j}}^{q_j}\left(\prod\limits_{p'=j+1}^{\zeta}\bm{Y}_{p'}S_{k_{p'},i_{p'}}\right),
\end{eqnarray}
\begin{eqnarray}\label{eq4:cross:thm:GMOI Perturbation Formula}
\lefteqn{\mathfrak{X}^{\bm{X}_{j-1,P},\bm{C}_N,\bm{D}_N,\bm{X}_{j,N}}=}\nonumber \\
&&\sum\limits_{\substack{S_{k_p,i_p}, p \in [1,j-2], \\
S_{k_{p'},i_{p'}}, p' \in [j+1,\zeta]}}\sum\limits_{\substack{k_{j-1},i_{j-1},\\k_c, i_c,q_c,k_d,i_d,q_d\\ k_{j},i_{j},q_j}}\frac{(\beta^{[\zeta+1]}([\lambda_k]_1^{j-2},\lambda_{k_{j-1}},\lambda_{k_c},\lambda_{k_d},\lambda_{k_j},[\lambda_k]_{j+1}^{\zeta}))^{([q]_1^{j-2}[0,0,q_d,q_j][q]_{j+1}^{\zeta})}}{ ([q!]_1^{j-2}[0!,0!,q_d!,q_j!][q!]_{j+1}^{\zeta})  } \nonumber \\
&& \times \left(\prod\limits_{p=1}^{j-2}S_{k_p,i_p}\bm{Y}_p\right)\bm{P}_{k_{j-1},i_{j-1}}\bm{Y}_{j-1}\bm{P}_{k_c,i_c}(\bm{N}_{k_c,i_c}^{q_c}-\bm{N}_{k_d,i_d}^{q_d})\bm{N}_{k_d,i_d}^{q_d}\bm{Y}_j\bm{N}_{k_{j},i_{j}}^{q_j}\left(\prod\limits_{p'=j+1}^{\zeta}\bm{Y}_{p'}S_{k_{p'},i_{p'}}\right) \nonumber \\
&+& \sum\limits_{\substack{S_{k_p,i_p}, p \in [1,j-2], \\
S_{k_{p'},i_{p'}}, p' \in [j+1,\zeta]}}\sum\limits_{\substack{k_{j-1},i_{j-1},\\k_c, i_c,q_c,k_d,i_d,q_d\\ k_{j},i_{j},q_j}}\frac{(\beta^{[\zeta+1]}([\lambda_k]_1^{j-2},\lambda_{k_{j-1}},\lambda_{k_c},\lambda_{k_d},\lambda_{k_j},[\lambda_k]_{j+1}^{\zeta}))^{([q]_1^{j-2}[0,q_c,0,q_j][q]_{j+1}^{\zeta})}}{ ([q!]_1^{j-2}[0!,q_c!,0!,q_j!][q!]_{j+1}^{\zeta})  } \nonumber \\
&& \times \left(\prod\limits_{p=1}^{j-2}S_{k_p,i_p}\bm{Y}_p\right)\bm{P}_{k_{j-1},i_{j-1}}\bm{Y}_{j-1}\bm{N}_{k_c,i_c}^{q_c}(\bm{N}_{k_c,i_c}^{q_c}-\bm{N}_{k_d,i_d}^{q_d})\bm{P}_{k_d,i_d}\bm{Y}_j\bm{N}_{k_{j},i_{j}}^{q_j}\left(\prod\limits_{p'=j+1}^{\zeta}\bm{Y}_{p'}S_{k_{p'},i_{p'}}\right) \nonumber \\
&+& \sum\limits_{\substack{S_{k_p,i_p}, p \in [1,j-2], \\
S_{k_{p'},i_{p'}}, p' \in [j+1,\zeta]}}\sum\limits_{\substack{k_{j-1},i_{j-1},\\k_c,i_c, k_d,i_d,q_d\\ k_{j},i_{j},q_j}} \frac{(\beta^{[\zeta]}([\lambda_k]_1^{j-2},\lambda_{k_{j-1}},\lambda_{k_c},\lambda_{k_j},[\lambda_k]_{j+1}^{\zeta}))^{([q]_1^{j-2}[0,q_d,q_j][q]_{j+1}^{\zeta})}}{ ([q!]_1^{j-2}[0!,q_d!,q_j!][q!]_{j+1}^{\zeta})  } \nonumber \\
&& \times \left(\prod\limits_{p=1}^{j-2}S_{k_p,i_p}\bm{Y}_p\right)\bm{P}_{k_{j-1},i_{j-1}}\bm{Y}_{j-1}\bm{P}_{k_c,i_c}\bm{N}_{k_d,i_d}^{q_d}\bm{Y}_j\bm{N}_{k_{j},i_{j}}^{q_j}\left(\prod\limits_{p'=j+1}^{\zeta}\bm{Y}_{p'}S_{k_{p'},i_{p'}}\right) \nonumber \\  
&-& \sum\limits_{\substack{S_{k_p,i_p}, p \in [1,j-2], \\
S_{k_{p'},i_{p'}}, p' \in [j+1,\zeta]}}\sum\limits_{\substack{k_{j-1},i_{j-1},\\k_c, i_c, q_c, k_d,i_d\\ k_{j},i_{j},q_j}} \frac{(\beta^{[\zeta]}([\lambda_k]_1^{j-2},\lambda_{k_{j-1}},\lambda_{k_d},\lambda_{k_j},[\lambda_k]_{j+1}^{\zeta}))^{([q]_1^{j-2}[0,q_c,q_j][q]_{j+1}^{\zeta})}}{ ([q!]_1^{j-2}[0!,q_c!,q_j!][q!]_{j+1}^{\zeta})  } \nonumber \\
&& \times \left(\prod\limits_{p=1}^{j-2}S_{k_p,i_p}\bm{Y}_p\right)\bm{P}_{k_{j-1},i_{j-1}}\bm{Y}_{j-1}\bm{N}_{k_c,i_c}^{q_c}\bm{P}_{k_d,i_d}\bm{Y}_j\bm{N}_{k_{j},i_{j}}^{q_j}\left(\prod\limits_{p'=j+1}^{\zeta}\bm{Y}_{p'}S_{k_{p'},i_{p'}}\right),
\end{eqnarray}
\begin{eqnarray}\label{eq5:cross:thm:GMOI Perturbation Formula}
\lefteqn{\mathfrak{X}^{\bm{X}_{j-1,N},\bm{C}_P,\bm{D}_P,\bm{X}_{j,P}}=}\nonumber \\
&& \sum\limits_{\substack{S_{k_p,i_p}, p \in [1,j-2], \\
S_{k_{p'},i_{p'}}, p' \in [j+1,\zeta]}}\sum\limits_{\substack{k_{j-1},i_{j-1}, q_{j-1}\\k_c, i_c,q_c,k_d,i_d,q_d\\ k_{j},i_{j}}}\frac{(\beta^{[\zeta+1]}([\lambda_k]_1^{j-2},\lambda_{k_{j-1}},\lambda_{k_c},\lambda_{k_d},\lambda_{k_j},[\lambda_k]_{j+1}^{\zeta}))^{([q]_1^{j-2}[q_{j-1},0,0,0][q]_{j+1}^{\zeta})}}{ ([q!]_1^{j-2}[q_{j-1}!,0!,0!,0!][q!]_{j+1}^{\zeta})  } \nonumber \\
&& \times \left(\prod\limits_{p=1}^{j-2}S_{k_p,i_p}\bm{Y}_p\right)\bm{N}_{k_{j-1},i_{j-1}}^{q_{j-1}}\bm{Y}_{j-1}\bm{P}_{k_c,i_c}(\bm{N}_{k_c,i_c}^{q_c}-\bm{N}_{k_d,i_d}^{q_d})\bm{P}_{k_d,i_d}\bm{Y}_j\bm{P}_{k_{j},i_{j}}\left(\prod\limits_{p'=j+1}^{\zeta}\bm{Y}_{p'}S_{k_{p'},i_{p'}}\right),
\end{eqnarray}
\begin{eqnarray}\label{eq6:cross:thm:GMOI Perturbation Formula}
\lefteqn{\mathfrak{X}^{\bm{X}_{j-1,N},\bm{C}_N,\bm{D}_N,\bm{X}_{j,P}}=}\nonumber \\
&& \sum\limits_{\substack{S_{k_p,i_p}, p \in [1,j-2], \\
S_{k_{p'},i_{p'}}, p' \in [j+1,\zeta]}}\sum\limits_{\substack{k_{j-1},i_{j-1},q_{j-1},\\k_c, i_c,q_c,k_d,i_d,q_d\\ k_{j},i_{j}}}\frac{(\beta^{[\zeta+1]}([\lambda_k]_1^{j-2},\lambda_{k_{j-1}},\lambda_{k_c},\lambda_{k_d},\lambda_{k_j},[\lambda_k]_{j+1}^{\zeta}))^{([q]_1^{j-2}[q_{j-1},0,q_d,0][q]_{j+1}^{\zeta})}}{ ([q!]_1^{j-2}[q_{j-1}!,0!,q_d!,0!][q!]_{j+1}^{\zeta})  } \nonumber \\
&& \times \left(\prod\limits_{p=1}^{j-2}S_{k_p,i_p}\bm{Y}_p\right)\bm{N}_{k_{j-1},i_{j-1}}^{q_{j-1}}\bm{Y}_{j-1}\bm{P}_{k_c,i_c}(\bm{N}_{k_c,i_c}^{q_c}-\bm{N}_{k_d,i_d}^{q_d})\bm{N}_{k_d,i_d}^{q_d}\bm{Y}_j\bm{P}_{k_{j},i_{j}}\left(\prod\limits_{p'=j+1}^{\zeta}\bm{Y}_{p'}S_{k_{p'},i_{p'}}\right) \nonumber \\
&+& \sum\limits_{\substack{S_{k_p,i_p}, p \in [1,j-2], \\
S_{k_{p'},i_{p'}}, p' \in [j+1,\zeta]}}\sum\limits_{\substack{k_{j-1},i_{j-1},q_{j-1},\\k_c, i_c,q_c,k_d,i_d,q_d\\ k_{j},i_{j}}}\frac{(\beta^{[\zeta+1]}([\lambda_k]_1^{j-2},\lambda_{k_{j-1}},\lambda_{k_c},\lambda_{k_d},\lambda_{k_j},[\lambda_k]_{j+1}^{\zeta}))^{([q]_1^{j-2}[q_{j-1},q_c,0,0][q]_{j+1}^{\zeta})}}{ ([q!]_1^{j-2}[q_{j-1}!,q_c!,0!,0!][q!]_{j+1}^{\zeta})  } \nonumber \\
&& \times \left(\prod\limits_{p=1}^{j-2}S_{k_p,i_p}\bm{Y}_p\right)\bm{N}_{k_{j-1},i_{j-1}}^{q_{j-1}}\bm{Y}_{j-1}\bm{N}_{k_c,i_c}^{q_c}(\bm{N}_{k_c,i_c}^{q_c}-\bm{N}_{k_d,i_d}^{q_d})\bm{P}_{k_d,i_d}\bm{Y}_j\bm{P}_{k_{j},i_{j}}\left(\prod\limits_{p'=j+1}^{\zeta}\bm{Y}_{p'}S_{k_{p'},i_{p'}}\right) \nonumber \\
&+& \sum\limits_{\substack{S_{k_p,i_p}, p \in [1,j-2], \\
S_{k_{p'},i_{p'}}, p' \in [j+1,\zeta]}}\sum\limits_{\substack{k_{j-1},i_{j-1},q_{j-1},\\k_c,i_c, k_d,i_d,q_d\\ k_{j},i_{j}}} \frac{(\beta^{[\zeta]}([\lambda_k]_1^{j-2},\lambda_{k_{j-1}},\lambda_{k_c},\lambda_{k_j},[\lambda_k]_{j+1}^{\zeta}))^{([q]_1^{j-2}[q_{j-1},q_d,0][q]_{j+1}^{\zeta})}}{ ([q!]_1^{j-2}[q_{j-1}!,q_d!,0!][q!]_{j+1}^{\zeta})  } \nonumber \\
&& \times \left(\prod\limits_{p=1}^{j-2}S_{k_p,i_p}\bm{Y}_p\right)\bm{N}_{k_{j-1},i_{j-1}}^{q_{j-1}}\bm{Y}_{j-1}\bm{P}_{k_c,i_c}\bm{N}_{k_d,i_d}^{q_d}\bm{Y}_j\bm{P}_{k_{j},i_{j}}\left(\prod\limits_{p'=j+1}^{\zeta}\bm{Y}_{p'}S_{k_{p'},i_{p'}}\right) \nonumber \\  
&-& \sum\limits_{\substack{S_{k_p,i_p}, p \in [1,j-2], \\
S_{k_{p'},i_{p'}}, p' \in [j+1,\zeta]}}\sum\limits_{\substack{k_{j-1},i_{j-1},q_{j-1},\\k_c, i_c, q_c, k_d,i_d\\ k_{j},i_{j}}} \frac{(\beta^{[\zeta]}([\lambda_k]_1^{j-2},\lambda_{k_{j-1}},\lambda_{k_d},\lambda_{k_j},[\lambda_k]_{j+1}^{\zeta}))^{([q]_1^{j-2}[q_{j-1},q_c,0][q]_{j+1}^{\zeta})}}{ ([q!]_1^{j-2}[q_{j-1}!,q_c!,0!][q!]_{j+1}^{\zeta})  } \nonumber \\
&& \times \left(\prod\limits_{p=1}^{j-2}S_{k_p,i_p}\bm{Y}_p\right)\bm{N}_{k_{j-1},i_{j-1}}^{q_{j-1}}\bm{Y}_{j-1}\bm{N}_{k_c,i_c}^{q_c}\bm{P}_{k_d,i_d}\bm{Y}_j\bm{P}_{k_{j},i_{j}}\left(\prod\limits_{p'=j+1}^{\zeta}\bm{Y}_{p'}S_{k_{p'},i_{p'}}\right),
\end{eqnarray}
\begin{eqnarray}\label{eq7:cross:thm:GMOI Perturbation Formula}
\lefteqn{\mathfrak{X}^{\bm{X}_{j-1,N},\bm{C}_P,\bm{D}_P,\bm{X}_{j,N}}=}\nonumber \\
&& \sum\limits_{\substack{S_{k_p,i_p}, p \in [1,j-2], \\
S_{k_{p'},i_{p'}}, p' \in [j+1,\zeta]}}\sum\limits_{\substack{k_{j-1},i_{j-1}, q_{j-1}\\k_c, i_c,q_c,k_d,i_d,q_d\\ k_{j},i_{j}, q_j}}\frac{(\beta^{[\zeta+1]}([\lambda_k]_1^{j-2},\lambda_{k_{j-1}},\lambda_{k_c},\lambda_{k_d},\lambda_{k_j},[\lambda_k]_{j+1}^{\zeta}))^{([q]_1^{j-2}[q_{j-1},0,0,q_j][q]_{j+1}^{\zeta})}}{ ([q!]_1^{j-2}[q_{j-1}!,0!,0!,q_j!][q!]_{j+1}^{\zeta})  } \nonumber \\
&& \times \left(\prod\limits_{p=1}^{j-2}S_{k_p,i_p}\bm{Y}_p\right)\bm{N}_{k_{j-1},i_{j-1}}^{q_{j-1}}\bm{Y}_{j-1}\bm{P}_{k_c,i_c}(\bm{N}_{k_c,i_c}^{q_c}-\bm{N}_{k_d,i_d}^{q_d})\bm{P}_{k_d,i_d}\bm{Y}_j\bm{N}_{k_{j},i_{j}}^{q_j}\left(\prod\limits_{p'=j+1}^{\zeta}\bm{Y}_{p'}S_{k_{p'},i_{p'}}\right),
\end{eqnarray}
\begin{eqnarray}\label{eq8:cross:thm:GMOI Perturbation Formula}
\lefteqn{\mathfrak{X}^{\bm{X}_{j-1,N},\bm{C}_N,\bm{D}_N,\bm{X}_{j,N}}=}\nonumber \\
&&  \sum\limits_{\substack{S_{k_p,i_p}, p \in [1,j-2], \\
S_{k_{p'},i_{p'}}, p' \in [j+1,\zeta]}}\sum\limits_{\substack{k_{j-1},i_{j-1},q_{j-1},\\k_c, i_c,q_c,k_d,i_d,q_d\\ k_{j},i_{j}, q_j}}\frac{(\beta^{[\zeta+1]}([\lambda_k]_1^{j-2},\lambda_{k_{j-1}},\lambda_{k_c},\lambda_{k_d},\lambda_{k_j},[\lambda_k]_{j+1}^{\zeta}))^{([q]_1^{j-2}[q_{j-1},0,q_d,q_j][q]_{j+1}^{\zeta})}}{ ([q!]_1^{j-2}[q_{j-1}!,0!,q_d!,q_j!][q!]_{j+1}^{\zeta})  } \nonumber \\
&& \times \left(\prod\limits_{p=1}^{j-2}S_{k_p,i_p}\bm{Y}_p\right)\bm{N}_{k_{j-1},i_{j-1}}^{q_{j-1}}\bm{Y}_{j-1}\bm{P}_{k_c,i_c}(\bm{N}_{k_c,i_c}^{q_c}-\bm{N}_{k_d,i_d}^{q_d})\bm{N}_{k_d,i_d}^{q_d}\bm{Y}_j\bm{N}_{k_{j},i_{j}}^{q_j}\left(\prod\limits_{p'=j+1}^{\zeta}\bm{Y}_{p'}S_{k_{p'},i_{p'}}\right) \nonumber \\
&+& \sum\limits_{\substack{S_{k_p,i_p}, p \in [1,j-2], \\
S_{k_{p'},i_{p'}}, p' \in [j+1,\zeta]}}\sum\limits_{\substack{k_{j-1},i_{j-1},q_{j-1},\\k_c, i_c,q_c,k_d,i_d,q_d\\ k_{j},i_{j},q_j}}\frac{(\beta^{[\zeta+1]}([\lambda_k]_1^{j-2},\lambda_{k_{j-1}},\lambda_{k_c},\lambda_{k_d},\lambda_{k_j},[\lambda_k]_{j+1}^{\zeta}))^{([q]_1^{j-2}[q_{j-1},q_c,0,q_j][q]_{j+1}^{\zeta})}}{ ([q!]_1^{j-2}[q_{j-1}!,q_c!,0!,q_j!][q!]_{j+1}^{\zeta})  } \nonumber \\
&& \times \left(\prod\limits_{p=1}^{j-2}S_{k_p,i_p}\bm{Y}_p\right)\bm{N}_{k_{j-1},i_{j-1}}^{q_{j-1}}\bm{Y}_{j-1}\bm{N}_{k_c,i_c}^{q_c}(\bm{N}_{k_c,i_c}^{q_c}-\bm{N}_{k_d,i_d}^{q_d})\bm{P}_{k_d,i_d}\bm{Y}_j\bm{N}_{k_{j},i_{j}}^{q_j}\left(\prod\limits_{p'=j+1}^{\zeta}\bm{Y}_{p'}S_{k_{p'},i_{p'}}\right) \nonumber \\
&+& \sum\limits_{\substack{S_{k_p,i_p}, p \in [1,j-2], \\
S_{k_{p'},i_{p'}}, p' \in [j+1,\zeta]}}\sum\limits_{\substack{k_{j-1},i_{j-1},q_{j-1},\\k_c,i_c, k_d,i_d,q_d\\ k_{j},i_{j}, q_j}} \frac{(\beta^{[\zeta]}([\lambda_k]_1^{j-2},\lambda_{k_{j-1}},\lambda_{k_c},\lambda_{k_j},[\lambda_k]_{j+1}^{\zeta}))^{([q]_1^{j-2}[q_{j-1},q_d,q_j][q]_{j+1}^{\zeta})}}{ ([q!]_1^{j-2}[q_{j-1}!,q_d!,q_j!][q!]_{j+1}^{\zeta})  } \nonumber \\
&& \times \left(\prod\limits_{p=1}^{j-2}S_{k_p,i_p}\bm{Y}_p\right)\bm{N}_{k_{j-1},i_{j-1}}^{q_{j-1}}\bm{Y}_{j-1}\bm{P}_{k_c,i_c}\bm{N}_{k_d,i_d}^{q_d}\bm{Y}_j\bm{N}_{k_{j},i_{j}}^{q_j}\left(\prod\limits_{p'=j+1}^{\zeta}\bm{Y}_{p'}S_{k_{p'},i_{p'}}\right) \nonumber \\  
&-& \sum\limits_{\substack{S_{k_p,i_p}, p \in [1,j-2], \\
S_{k_{p'},i_{p'}}, p' \in [j+1,\zeta]}}\sum\limits_{\substack{k_{j-1},i_{j-1},q_{j-1},\\k_c, i_c, q_c, k_d,i_d\\ k_{j},i_{j},q_j}} \frac{(\beta^{[\zeta]}([\lambda_k]_1^{j-2},\lambda_{k_{j-1}},\lambda_{k_d},\lambda_{k_j},[\lambda_k]_{j+1}^{\zeta}))^{([q]_1^{j-2}[q_{j-1},q_c,q_j][q]_{j+1}^{\zeta})}}{ ([q!]_1^{j-2}[q_{j-1}!,q_c!,q_j!][q!]_{j+1}^{\zeta})  } \nonumber \\
&& \times \left(\prod\limits_{p=1}^{j-2}S_{k_p,i_p}\bm{Y}_p\right)\bm{N}_{k_{j-1},i_{j-1}}^{q_{j-1}}\bm{Y}_{j-1}\bm{N}_{k_c,i_c}^{q_c}\bm{P}_{k_d,i_d}\bm{Y}_j\bm{N}_{k_{j},i_{j}}^{q_j}\left(\prod\limits_{p'=j+1}^{\zeta}\bm{Y}_{p'}S_{k_{p'},i_{p'}}\right).
\end{eqnarray}
\end{theorem}
\textbf{Proof:}
From Proposition~\ref{prop:GMOI decomp by parameters X P X N}, we have 
{\small
\begin{eqnarray}\label{eq3:thm:GMOI Perturbation Formula}
\lefteqn{T_{\beta^{[\zeta+1]}}^{[\bm{X}]_1^{j-1},\bm{C},\bm{D},[\bm{X}]_j^{\zeta}}([\bm{Y}]_1^{j-1},\bm{C}-\bm{D},[\bm{Y}]_j^{\zeta})=\bm{O}}\nonumber \\
&+_1& T_{\beta^{[\zeta+1]}}^{[\bm{X}]_1^{j-2},\bm{X}_{j-1,P},\bm{C}_P,\bm{D}_P,\bm{X}_{j,P},[\bm{X}]_{j+1}^{\zeta}}([\bm{Y}]_1^{j-1},\bm{C}-\bm{D},[\bm{Y}]_j^{\zeta})\nonumber \\
&+_2& T_{\beta^{[\zeta+1]}}^{[\bm{X}]_1^{j-2},\bm{X}_{j-1,P},\bm{C}_N,\bm{D}_P,\bm{X}_{j,P},[\bm{X}]_{j+1}^{\zeta}}([\bm{Y}]_1^{j-1},\bm{C}-\bm{D},[\bm{Y}]_j^{\zeta})\nonumber \\
&+_3& T_{\beta^{[\zeta+1]}}^{[\bm{X}]_1^{j-2},\bm{X}_{j-1,P},\bm{C}_P,\bm{D}_N,\bm{X}_{j,P},[\bm{X}]_{j+1}^{\zeta}}([\bm{Y}]_1^{j-1},\bm{C}-\bm{D},[\bm{Y}]_j^{\zeta})\nonumber \\
&+_4& T_{\beta^{[\zeta+1]}}^{[\bm{X}]_1^{j-2},\bm{X}_{j-1,P},\bm{C}_P,\bm{D}_P,\bm{X}_{j,N},[\bm{X}]_{j+1}^{\zeta}}([\bm{Y}]_1^{j-1},\bm{C}-\bm{D},[\bm{Y}]_j^{\zeta}) \nonumber \\
&+_5& T_{\beta^{[\zeta+1]}}^{[\bm{X}]_1^{j-2},\bm{X}_{j-1,P},\bm{C}_N,\bm{D}_P,\bm{X}_{j,N},[\bm{X}]_{j+1}^{\zeta}}([\bm{Y}]_1^{j-1},\bm{C}-\bm{D},[\bm{Y}]_j^{\zeta})\nonumber \\
&+_6& T_{\beta^{[\zeta+1]}}^{[\bm{X}]_1^{j-2},\bm{X}_{j-1,P},\bm{C}_P,\bm{D}_N,\bm{X}_{j,N},[\bm{X}]_{j+1}^{\zeta}}([\bm{Y}]_1^{j-1},\bm{C}-\bm{D},[\bm{Y}]_j^{\zeta})\nonumber \\
&+_7& T_{\beta^{[\zeta+1]}}^{[\bm{X}]_1^{j-2},\bm{X}_{j-1,N},\bm{C}_P,\bm{D}_P,\bm{X}_{j,P},[\bm{X}]_{j+1}^{\zeta}}([\bm{Y}]_1^{j-1},\bm{C}-\bm{D},[\bm{Y}]_j^{\zeta})\nonumber \\
&+_8& T_{\beta^{[\zeta+1]}}^{[\bm{X}]_1^{j-2},\bm{X}_{j-1,N},\bm{C}_N,\bm{D}_P,\bm{X}_{j,P},[\bm{X}]_{j+1}^{\zeta}}([\bm{Y}]_1^{j-1},\bm{C}-\bm{D},[\bm{Y}]_j^{\zeta})\nonumber \\
&+_9& T_{\beta^{[\zeta+1]}}^{[\bm{X}]_1^{j-2},\bm{X}_{j-1,N},\bm{C}_P,\bm{D}_N,\bm{X}_{j,P},[\bm{X}]_{j+1}^{\zeta}}([\bm{Y}]_1^{j-1},\bm{C}-\bm{D},[\bm{Y}]_j^{\zeta})\nonumber \\
&+_{10}&T_{\beta^{[\zeta+1]}}^{[\bm{X}]_1^{j-2},\bm{X}_{j-1,N},\bm{C}_P,\bm{D}_P,\bm{X}_{j,N},[\bm{X}]_{j+1}^{\zeta}}([\bm{Y}]_1^{j-1},\bm{C}-\bm{D},[\bm{Y}]_j^{\zeta})\nonumber \\
&+_{11}& T_{\beta^{[\zeta+1]}}^{[\bm{X}]_1^{j-2},\bm{X}_{j-1,N},\bm{C}_N,\bm{D}_P,\bm{X}_{j,N},[\bm{X}]_{j+1}^{\zeta}}([\bm{Y}]_1^{j-1},\bm{C}-\bm{D},[\bm{Y}]_j^{\zeta})\nonumber \\
&+_{12}& T_{\beta^{[\zeta+1]}}^{[\bm{X}]_1^{j-2},\bm{X}_{j-1,N},\bm{C}_P,\bm{D}_N,\bm{X}_{j,N},[\bm{X}]_{j+1}^{\zeta}}([\bm{Y}]_1^{j-1},\bm{C}-\bm{D},[\bm{Y}]_j^{\zeta})\nonumber \\
&+_{}& T_{\beta^{[\zeta+1]}}^{[\bm{X}]_1^{j-2},\bm{X}_{j-1,P},\bm{C}_N,\bm{D}_N,\bm{X}_{j,P},[\bm{X}]_{j+1}^{\zeta}}([\bm{Y}]_1^{j-1},\bm{C}-\bm{D},[\bm{Y}]_j^{\zeta})\nonumber \\
&+& T_{\beta^{[\zeta+1]}}^{[\bm{X}]_1^{j-2},\bm{X}_{j-1,P},\bm{C}_N,\bm{D}_N,\bm{X}_{j,N},[\bm{X}]_{j+1}^{\zeta}}([\bm{Y}]_1^{j-1},\bm{C}-\bm{D},[\bm{Y}]_j^{\zeta})\nonumber \\
&+& T_{\beta^{[\zeta+1]}}^{[\bm{X}]_1^{j-2},\bm{X}_{j-1,N},\bm{C}_N,\bm{D}_N,\bm{X}_{j,P},[\bm{X}]_{j+1}^{\zeta}}([\bm{Y}]_1^{j-1},\bm{C}-\bm{D},[\bm{Y}]_j^{\zeta})\nonumber \\
&+& T_{\beta^{[\zeta+1]}}^{[\bm{X}]_1^{j-2},\bm{X}_{j-1,N},\bm{C}_N,\bm{D}_N,\bm{X}_{j,N},[\bm{X}]_{j+1}^{\zeta}}([\bm{Y}]_1^{j-1},\bm{C}-\bm{D},[\bm{Y}]_j^{\zeta}), 
\end{eqnarray}
}

By comparing Eq.~\eqref{eq1:thm:GMOI Perturbation Formula} and Eq.~\eqref{eq3:thm:GMOI Perturbation Formula}, we claim we have the following pairs of identitoes:
\begin{eqnarray}\label{eq4-1:thm:GMOI Perturbation Formula}
\mbox{GMOI after $+_1$} &=&I+\mathfrak{X}^{\bm{X}_{j-1,P},\bm{C}_P,\bm{D}_P,\bm{X}_{j,P}}, 
\end{eqnarray}
\begin{eqnarray}\label{eq4-2:thm:GMOI Perturbation Formula}
 \mbox{GMOIs after $+_2$ and $+_3$} &=&II+\mathfrak{X}^{\bm{X}_{j-1,P},\bm{C}_N,\bm{D}_N,\bm{X}_{j,P}}, 
\end{eqnarray}
\begin{eqnarray}\label{eq4-3:thm:GMOI Perturbation Formula}
\mbox{GMOI after $+_4$} &=&III+\mathfrak{X}^{\bm{X}_{j-1,P},\bm{C}_P,\bm{D}_P,\bm{X}_{j,N}},  
\end{eqnarray}
\begin{eqnarray}\label{eq4-4:thm:GMOI Perturbation Formula}
\mbox{GMOIs after $+_5$ and $+_6$} &=&IV+\mathfrak{X}^{\bm{X}_{j-1,P},\bm{C}_N,\bm{D}_N,\bm{X}_{j,N}},  
\end{eqnarray}
\begin{eqnarray}\label{eq4-5:thm:GMOI Perturbation Formula}
\mbox{GMOI after $+_7$} &=&V+\mathfrak{X}^{\bm{X}_{j-1,N},\bm{C}_P,\bm{D}_P,\bm{X}_{j,P}},  
\end{eqnarray}
\begin{eqnarray}\label{eq4-6:thm:GMOI Perturbation Formula}
\mbox{GMOIs after $+_8$ and $+_9$} &=&VI+\mathfrak{X}^{\bm{X}_{j-1,N},\bm{C}_N,\bm{D}_N,\bm{X}_{j,P}},  
\end{eqnarray}
\begin{eqnarray}\label{eq4-7:thm:GMOI Perturbation Formula}
\mbox{GMOI after $+_{10}$} &=&VII+\mathfrak{X}^{\bm{X}_{j-1,N},\bm{C}_P,\bm{D}_P,\bm{X}_{j,N}},  
\end{eqnarray}
\begin{eqnarray}\label{eq4-8:thm:GMOI Perturbation Formula}
\mbox{GMOIs after $+_{11}$ and $+_{12}$} &=&VIII+\mathfrak{X}^{\bm{X}_{j-1,N},\bm{C}_N,\bm{D}_N,\bm{X}_{j,N}},
\end{eqnarray}
We will provide detailed proof steps for identity given by Eq.~\eqref{eq4-1:thm:GMOI Perturbation Formula} and for identity given by Eq.~\eqref{eq4-2:thm:GMOI Perturbation Formula} as other identities from Eq.~\eqref{eq4-3:thm:GMOI Perturbation Formula} to Eq.~\eqref{eq4-8:thm:GMOI Perturbation Formula} can be proved similarly. 

\textbf{Identity  given by Eq.~\eqref{eq4-1:thm:GMOI Perturbation Formula} Proof}

Because we have 
{\tiny
\begin{eqnarray}\label{eq1:eq4-1:thm:GMOI Perturbation Formula Proof}
\lefteqn{T_{\beta^{[\zeta]}}^{[\bm{X}]_1^{j-2},\bm{X}_{j-1,P},\bm{C}_P,\bm{X}_{j,P},[\bm{X}]_{j+1}^{\zeta}}([\bm{Y}]_1^{j-1},[\bm{Y}]_j^{\zeta}) - T_{\beta^{[\zeta]}}^{[\bm{X}]_1^{j-2},\bm{X}_{j-1,P},\bm{D}_P,\bm{X}_{j,P},[\bm{X}]_{j+1}^{\zeta}}([\bm{Y}]_1^{j-1},[\bm{Y}]_j^{\zeta})=}\nonumber \\
&&\sum\limits_{\substack{S_{k_p,i_p}, p \in [1,j-2], \\
S_{k_{p'},i_{p'}}, p' \in [j+1,\zeta]}}\sum\limits_{\substack{k_{j-1},i_{j-1},\\k_c, i_c,\\ k_{j},i_{j}}} \frac{(\beta^{[\zeta]}([\lambda_k]_1^{j-2},\lambda_{k_{j-1}},\lambda_{k_c},\lambda_{k_j},[\lambda_k]_{j+1}^{\zeta}))^{([q]_1^{j-2}[0,0,0][q]_{j+1}^{\zeta})}}{ ([q!]_1^{j-2}[0!,0!,0!][q!]_{j+1}^{\zeta})  } \nonumber \\
&& \times \left(\prod\limits_{p=1}^{j-2}S_{k_p,i_p}\bm{Y}_p\right)\bm{P}_{k_{j-1},i_{j-1}}\bm{Y}_{j-1}\bm{P}_{k_c,i_c}\bm{Y}_j\bm{P}_{k_{j},i_{j}}\left(\prod\limits_{p'=j+1}^{\zeta}\bm{Y}_{p'}S_{k_{p'},i_{p'}}\right) \nonumber \\
&-& \sum\limits_{\substack{S_{k_p,i_p}, p \in [1,j-2], \\
S_{k_{p'},i_{p'}}, p' \in [j+1,\zeta]}}\sum\limits_{\substack{k_{j-1},i_{j-1},\\k_d,i_d,\\ k_{j},i_{j}}} \frac{(\beta^{[\zeta]}([\lambda_k]_1^{j-2},\lambda_{k_{j-1}},\lambda_{k_d},\lambda_{k_j},[\lambda_k]_{j+1}^{\zeta}))^{([q]_1^{j-2}[0,0,0][q]_{j+1}^{\zeta})}}{ ([q!]_1^{j-2}[0!,0!,0!][q!]_{j+1}^{\zeta})  } \nonumber \\
&& \times \left(\prod\limits_{p=1}^{j-2}S_{k_p,i_p}\bm{Y}_p\right)\bm{P}_{k_{j-1},i_{j-1}}\bm{Y}_{j-1}\bm{P}_{k_d,i_d}\bm{Y}_j\bm{P}_{k_{j},i_{j}}\left(\prod\limits_{p'=j+1}^{\zeta}\bm{Y}_{p'}S_{k_{p'},i_{p'}}\right) \nonumber \\
&=_1& \sum\limits_{\substack{S_{k_p,i_p}, p \in [1,j-2], \\
S_{k_{p'},i_{p'}}, p' \in [j+1,\zeta]}}\sum\limits_{\substack{k_{j-1},i_{j-1},\\k_c, i_c,k_d,i_d,\\ k_{j},i_{j}}} \frac{(\beta^{[\zeta]}([\lambda_k]_1^{j-2},\lambda_{k_{j-1}},\lambda_{k_c},\lambda_{k_j},[\lambda_k]_{j+1}^{\zeta}))^{([q]_1^{j-2}[0,0,0][q]_{j+1}^{\zeta})}}{ ([q!]_1^{j-2}[0!,0!,0!][q!]_{j+1}^{\zeta})  } \nonumber \\
&& \times \left(\prod\limits_{p=1}^{j-2}S_{k_p,i_p}\bm{Y}_p\right)\bm{P}_{k_{j-1},i_{j-1}}\bm{Y}_{j-1}\bm{P}_{k_c,i_c}\bm{P}_{k_d,i_d}\bm{Y}_j\bm{P}_{k_{j},i_{j}}\left(\prod\limits_{p'=j+1}^{\zeta}\bm{Y}_{p'}S_{k_{p'},i_{p'}}\right) \nonumber \\
&-& \sum\limits_{\substack{S_{k_p,i_p}, p \in [1,j-2], \\
S_{k_{p'},i_{p'}}, p' \in [j+1,\zeta]}}\sum\limits_{\substack{k_{j-1},i_{j-1},\\k_c,i_c, k_d,i_d,\\ k_{j},i_{j}}}\frac{(\beta^{[\zeta]}([\lambda_k]_1^{j-2},\lambda_{k_{j-1}},\lambda_{k_d},\lambda_{k_j},[\lambda_k]_{j+1}^{\zeta}))^{([q]_1^{j-2}[0,0,0][q]_{j+1}^{\zeta})}}{ ([q!]_1^{j-2}[0!,0!,0!][q!]_{j+1}^{\zeta})  } \nonumber \\
&& \times \left(\prod\limits_{p=1}^{j-2}S_{k_p,i_p}\bm{Y}_p\right)\bm{P}_{k_{j-1},i_{j-1}}\bm{Y}_{j-1}\bm{P}_{k_c,i_c}\bm{P}_{k_d,i_d}\bm{Y}_j\bm{P}_{k_{j},i_{j}}\left(\prod\limits_{p'=j+1}^{\zeta}\bm{Y}_{p'}S_{k_{p'},i_{p'}}\right) \nonumber \\
&=_2&   \sum\limits_{\substack{S_{k_p,i_p}, p \in [1,j-2], \\
S_{k_{p'},i_{p'}}, p' \in [j+1,\zeta]}}\sum\limits_{\substack{k_{j-1},i_{j-1},\\k_c, i_c,k_d,i_d,\\ k_{j},i_{j}}}(\lambda_{k_c} - \lambda_{k_d})\frac{(\beta^{[\zeta+1]}([\lambda_k]_1^{j-2},\lambda_{k_{j-1}},\lambda_{k_c},\lambda_{k_d},\lambda_{k_j},[\lambda_k]_{j+1}^{\zeta}))^{([q]_1^{j-2}[0,0,0,0][q]_{j+1}^{\zeta})}}{ ([q!]_1^{j-2}[0!,0!,0!,0!][q!]_{j+1}^{\zeta})  } \nonumber \\
&& \times \left(\prod\limits_{p=1}^{j-2}S_{k_p,i_p}\bm{Y}_p\right)\bm{P}_{k_{j-1},i_{j-1}}\bm{Y}_{j-1}\bm{P}_{k_c,i_c}(\bm{P}_{k_c,i_c}-\bm{P}_{k_d,i_d})\bm{P}_{k_d,i_d}\bm{Y}_j\bm{P}_{k_{j},i_{j}}\left(\prod\limits_{p'=j+1}^{\zeta}\bm{Y}_{p'}S_{k_{p'},i_{p'}}\right)\nonumber \\
&=&  T_{\beta^{[\zeta+1]}}^{[\bm{X}]_1^{j-2},\bm{X}_{j-1,P},\bm{C}_P,\bm{D}_P,\bm{X}_{j,P},[\bm{X}]_{j+1}^{\zeta}}([\bm{Y}]_1^{j-1},\bm{C}-\bm{D},[\bm{Y}]_j^{\zeta}) \nonumber \\
&-&   \sum\limits_{\substack{S_{k_p,i_p}, p \in [1,j-2], \\
S_{k_{p'},i_{p'}}, p' \in [j+1,\zeta]}}\sum\limits_{\substack{k_{j-1},i_{j-1},\\k_c, i_c,q_c,k_d,i_d,q_d\\ k_{j},i_{j}}}\frac{(\beta^{[\zeta+1]}([\lambda_k]_1^{j-2},\lambda_{k_{j-1}},\lambda_{k_c},\lambda_{k_d},\lambda_{k_j},[\lambda_k]_{j+1}^{\zeta}))^{([q]_1^{j-2}[0,0,0,0][q]_{j+1}^{\zeta})}}{ ([q!]_1^{j-2}[0!,0!,0!,0!][q!]_{j+1}^{\zeta})  } \nonumber \\
&& \times \left(\prod\limits_{p=1}^{j-2}S_{k_p,i_p}\bm{Y}_p\right)\bm{P}_{k_{j-1},i_{j-1}}\bm{Y}_{j-1}\bm{P}_{k_c,i_c}(\bm{N}_{k_c,i_c}^{q_c}-\bm{N}_{k_d,i_d}^{q_d})\bm{P}_{k_d,i_d}\bm{Y}_j\bm{P}_{k_{j},i_{j}}\left(\prod\limits_{p'=j+1}^{\zeta}\bm{Y}_{p'}S_{k_{p'},i_{p'}}\right),
\end{eqnarray}
}
where we apply $\sum\limits_{k_c,i_c}\bm{P}_{k_c,i_c} = \sum\limits_{k_d,i_d}\bm{P}_{k_d,i_d} = \bm{I}$ in $=_1$,  apply Lemma~\ref{lma:first-order divided difference identity} in $=_2$. Moreover, the perturbation correct term $\mathfrak{X}^{\bm{X}_{j-1,P},\bm{C}_P,\bm{D}_P,\bm{X}_{j,P}}$ can be expressed by
\begin{eqnarray}\label{eq1:eq4-1:thm:GMOI Perturbation Formula Proof}
\lefteqn{\mathfrak{X}^{\bm{X}_{j-1,P},\bm{C}_P,\bm{D}_P,\bm{X}_{j,P}}=}\nonumber \\
&& \sum\limits_{\substack{S_{k_p,i_p}, p \in [1,j-2], \\
S_{k_{p'},i_{p'}}, p' \in [j+1,\zeta]}}\sum\limits_{\substack{k_{j-1},i_{j-1},\\k_c, i_c,q_c,k_d,i_d,q_d\\ k_{j},i_{j}}}\frac{(\beta^{[\zeta+1]}([\lambda_k]_1^{j-2},\lambda_{k_{j-1}},\lambda_{k_c},\lambda_{k_d},\lambda_{k_j},[\lambda_k]_{j+1}^{\zeta}))^{([q]_1^{j-2}[0,0,0,0][q]_{j+1}^{\zeta})}}{ ([q!]_1^{j-2}[0!,0!,0!,0!][q!]_{j+1}^{\zeta})  } \nonumber \\
&& \times \left(\prod\limits_{p=1}^{j-2}S_{k_p,i_p}\bm{Y}_p\right)\bm{P}_{k_{j-1},i_{j-1}}\bm{Y}_{j-1}\bm{P}_{k_c,i_c}(\bm{N}_{k_c,i_c}^{q_c}-\bm{N}_{k_d,i_d}^{q_d})\bm{P}_{k_d,i_d}\bm{Y}_j\bm{P}_{k_{j},i_{j}}\left(\prod\limits_{p'=j+1}^{\zeta}\bm{Y}_{p'}S_{k_{p'},i_{p'}}\right).
\end{eqnarray}

\textbf{Identity  given by Eq.~\eqref{eq4-2:thm:GMOI Perturbation Formula} Proof}

Because we have 
{\tiny
\begin{eqnarray}\label{eq1:eq4-2:thm:GMOI Perturbation Formula Proof}
\lefteqn{T_{\beta^{[\zeta]}}^{[\bm{X}]_1^{j-2},\bm{X}_{j-1,P},\bm{C}_N,\bm{X}_{j,P},[\bm{X}]_{j+1}^{\zeta}}([\bm{Y}]_1^{j-1},[\bm{Y}]_j^{\zeta}) - T_{\beta^{[\zeta]}}^{[\bm{X}]_1^{j-2},\bm{X}_{j-1,P},\bm{D}_N,\bm{X}_{j,P},[\bm{X}]_{j+1}^{\zeta}}([\bm{Y}]_1^{j-1},[\bm{Y}]_j^{\zeta})=}\nonumber \\
&&\sum\limits_{\substack{S_{k_p,i_p}, p \in [1,j-2], \\
S_{k_{p'},i_{p'}}, p' \in [j+1,\zeta]}}\sum\limits_{\substack{k_{j-1},i_{j-1},\\k_c, i_c, q_c\\ k_{j},i_{j}}} \frac{(\beta^{[\zeta]}([\lambda_k]_1^{j-2},\lambda_{k_{j-1}},\lambda_{k_c},\lambda_{k_j},[\lambda_k]_{j+1}^{\zeta}))^{([q]_1^{j-2}[0,q_c,0][q]_{j+1}^{\zeta})}}{ ([q!]_1^{j-2}[0!,q_c!,0!][q!]_{j+1}^{\zeta})  } \nonumber \\
&& \times \left(\prod\limits_{p=1}^{j-2}S_{k_p,i_p}\bm{Y}_p\right)\bm{P}_{k_{j-1},i_{j-1}}\bm{Y}_{j-1}\bm{N}_{k_c,i_c}^{q_c}\bm{Y}_j\bm{P}_{k_{j},i_{j}}\left(\prod\limits_{p'=j+1}^{\zeta}\bm{Y}_{p'}S_{k_{p'},i_{p'}}\right) \nonumber \\
&-& \sum\limits_{\substack{S_{k_p,i_p}, p \in [1,j-2], \\
S_{k_{p'},i_{p'}}, p' \in [j+1,\zeta]}}\sum\limits_{\substack{k_{j-1},i_{j-1},\\k_d,i_d,q_d\\ k_{j},i_{j}}} \frac{(\beta^{[\zeta]}([\lambda_k]_1^{j-2},\lambda_{k_{j-1}},\lambda_{k_d},\lambda_{k_j},[\lambda_k]_{j+1}^{\zeta}))^{([q]_1^{j-2}[0,q_d,0][q]_{j+1}^{\zeta})}}{ ([q!]_1^{j-2}[0!,q_d!,0!][q!]_{j+1}^{\zeta})  } \nonumber \\
&& \times \left(\prod\limits_{p=1}^{j-2}S_{k_p,i_p}\bm{Y}_p\right)\bm{P}_{k_{j-1},i_{j-1}}\bm{Y}_{j-1}\bm{N}_{k_d,i_d}^{q_d}\bm{Y}_j\bm{P}_{k_{j},i_{j}}\left(\prod\limits_{p'=j+1}^{\zeta}\bm{Y}_{p'}S_{k_{p'},i_{p'}}\right) \nonumber \\
&=_1&\sum\limits_{\substack{S_{k_p,i_p}, p \in [1,j-2], \\
S_{k_{p'},i_{p'}}, p' \in [j+1,\zeta]}}\sum\limits_{\substack{k_{j-1},i_{j-1},\\k_c, i_c, q_c, k_d,i_d\\ k_{j},i_{j}}} \frac{(\beta^{[\zeta]}([\lambda_k]_1^{j-2},\lambda_{k_{j-1}},\lambda_{k_c},\lambda_{k_j},[\lambda_k]_{j+1}^{\zeta}))^{([q]_1^{j-2}[0,q_c,0][q]_{j+1}^{\zeta})}}{ ([q!]_1^{j-2}[0!,q_c!,0!][q!]_{j+1}^{\zeta})  } \nonumber \\
&& \times \left(\prod\limits_{p=1}^{j-2}S_{k_p,i_p}\bm{Y}_p\right)\bm{P}_{k_{j-1},i_{j-1}}\bm{Y}_{j-1}\bm{N}_{k_c,i_c}^{q_c}\bm{P}_{k_d,i_d}\bm{Y}_j\bm{P}_{k_{j},i_{j}}\left(\prod\limits_{p'=j+1}^{\zeta}\bm{Y}_{p'}S_{k_{p'},i_{p'}}\right) \nonumber \\
&-& \sum\limits_{\substack{S_{k_p,i_p}, p \in [1,j-2], \\
S_{k_{p'},i_{p'}}, p' \in [j+1,\zeta]}}\sum\limits_{\substack{k_{j-1},i_{j-1},\\k_c,i_c, k_d,i_d,q_d\\ k_{j},i_{j}}} \frac{(\beta^{[\zeta]}([\lambda_k]_1^{j-2},\lambda_{k_{j-1}},\lambda_{k_d},\lambda_{k_j},[\lambda_k]_{j+1}^{\zeta}))^{([q]_1^{j-2}[0,q_d,0][q]_{j+1}^{\zeta})}}{ ([q!]_1^{j-2}[0!,q_d!,0!][q!]_{j+1}^{\zeta})  } \nonumber \\
&& \times \left(\prod\limits_{p=1}^{j-2}S_{k_p,i_p}\bm{Y}_p\right)\bm{P}_{k_{j-1},i_{j-1}}\bm{Y}_{j-1}\bm{P}_{k_c,i_c}\bm{N}_{k_d,i_d}^{q_d}\bm{Y}_j\bm{P}_{k_{j},i_{j}}\left(\prod\limits_{p'=j+1}^{\zeta}\bm{Y}_{p'}S_{k_{p'},i_{p'}}\right),
\end{eqnarray}
}
where we apply $\sum\limits_{k_c,i_c}\bm{P}_{k_c,i_c} = \sum\limits_{k_d,i_d}\bm{P}_{k_d,i_d} = \bm{I}$ in $=_1$.

Different proof in Eq.~\eqref{eq4-1:thm:GMOI Perturbation Formula}, the matrices production part of the first and the second terms are different, i.e., $\bm{N}_{k_c,i_c}^{q_c}\bm{P}_{k_d,i_d} \neq \bm{P}_{k_c,i_c}\bm{N}_{k_d,i_d}^{q_d}$.  Therefore, we will add and subtract the same auxillary terms to make two GMOI terems agree with GMOIs terms after $+_2$ and $+_3$ gievn by Eq.~\eqref{eq3:thm:GMOI Perturbation Formula}. Continue from Eq.~\eqref{eq1:eq4-2:thm:GMOI Perturbation Formula Proof}, we have 
{\tiny
\begin{eqnarray}\label{eq2:eq4-2:thm:GMOI Perturbation Formula Proof}
\lefteqn{T_{\beta^{[\zeta]}}^{[\bm{X}]_1^{j-2},\bm{X}_{j-1,P},\bm{C}_N,\bm{X}_{j,P},[\bm{X}]_{j+1}^{\zeta}}([\bm{Y}]_1^{j-1},[\bm{Y}]_j^{\zeta}) - T_{\beta^{[\zeta]}}^{[\bm{X}]_1^{j-2},\bm{X}_{j-1,P},\bm{D}_N,\bm{X}_{j,P},[\bm{X}]_{j+1}^{\zeta}}([\bm{Y}]_1^{j-1},[\bm{Y}]_j^{\zeta})=}\nonumber \\
&&\Bigg[\sum\limits_{\substack{S_{k_p,i_p}, p \in [1,j-2], \\
S_{k_{p'},i_{p'}}, p' \in [j+1,\zeta]}}\sum\limits_{\substack{k_{j-1},i_{j-1},\\k_c, i_c, q_c, k_d,i_d\\ k_{j},i_{j}}} \frac{(\beta^{[\zeta]}([\lambda_k]_1^{j-2},\lambda_{k_{j-1}},\lambda_{k_c},\lambda_{k_j},[\lambda_k]_{j+1}^{\zeta}))^{([q]_1^{j-2}[0,q_c,0][q]_{j+1}^{\zeta})}}{ ([q!]_1^{j-2}[0!,q_c!,0!][q!]_{j+1}^{\zeta})  } \nonumber \\
&& \times \left(\prod\limits_{p=1}^{j-2}S_{k_p,i_p}\bm{Y}_p\right)\bm{P}_{k_{j-1},i_{j-1}}\bm{Y}_{j-1}\bm{N}_{k_c,i_c}^{q_c}\bm{P}_{k_d,i_d}\bm{Y}_j\bm{P}_{k_{j},i_{j}}\left(\prod\limits_{p'=j+1}^{\zeta}\bm{Y}_{p'}S_{k_{p'},i_{p'}}\right) \nonumber \\
&-& \sum\limits_{\substack{S_{k_p,i_p}, p \in [1,j-2], \\
S_{k_{p'},i_{p'}}, p' \in [j+1,\zeta]}}\sum\limits_{\substack{k_{j-1},i_{j-1},\\k_c, i_c, q_c, k_d,i_d\\ k_{j},i_{j}}} \frac{(\beta^{[\zeta]}([\lambda_k]_1^{j-2},\lambda_{k_{j-1}},\lambda_{k_d},\lambda_{k_j},[\lambda_k]_{j+1}^{\zeta}))^{([q]_1^{j-2}[0,q_c,0][q]_{j+1}^{\zeta})}}{ ([q!]_1^{j-2}[0!,q_c!,0!][q!]_{j+1}^{\zeta})  } \nonumber \\
&& \times \left(\prod\limits_{p=1}^{j-2}S_{k_p,i_p}\bm{Y}_p\right)\bm{P}_{k_{j-1},i_{j-1}}\bm{Y}_{j-1}\bm{N}_{k_c,i_c}^{q_c}\bm{P}_{k_d,i_d}\bm{Y}_j\bm{P}_{k_{j},i_{j}}\left(\prod\limits_{p'=j+1}^{\zeta}\bm{Y}_{p'}S_{k_{p'},i_{p'}}\right) \Bigg]\nonumber \\
&+& \sum\limits_{\substack{S_{k_p,i_p}, p \in [1,j-2], \\
S_{k_{p'},i_{p'}}, p' \in [j+1,\zeta]}}\sum\limits_{\substack{k_{j-1},i_{j-1},\\k_c, i_c, q_c, k_d,i_d\\ k_{j},i_{j}}} \frac{(\beta^{[\zeta]}([\lambda_k]_1^{j-2},\lambda_{k_{j-1}},\lambda_{k_d},\lambda_{k_j},[\lambda_k]_{j+1}^{\zeta}))^{([q]_1^{j-2}[0,q_c,0][q]_{j+1}^{\zeta})}}{ ([q!]_1^{j-2}[0!,q_c!,0!][q!]_{j+1}^{\zeta})  } \nonumber \\
&& \times \left(\prod\limits_{p=1}^{j-2}S_{k_p,i_p}\bm{Y}_p\right)\bm{P}_{k_{j-1},i_{j-1}}\bm{Y}_{j-1}\bm{N}_{k_c,i_c}^{q_c}\bm{P}_{k_d,i_d}\bm{Y}_j\bm{P}_{k_{j},i_{j}}\left(\prod\limits_{p'=j+1}^{\zeta}\bm{Y}_{p'}S_{k_{p'},i_{p'}}\right) \nonumber \\
&-&\bigg[\sum\limits_{\substack{S_{k_p,i_p}, p \in [1,j-2], \\
S_{k_{p'},i_{p'}}, p' \in [j+1,\zeta]}}\sum\limits_{\substack{k_{j-1},i_{j-1},\\k_c,i_c, k_d,i_d,q_d\\ k_{j},i_{j}}} \frac{(\beta^{[\zeta]}([\lambda_k]_1^{j-2},\lambda_{k_{j-1}},\lambda_{k_d},\lambda_{k_j},[\lambda_k]_{j+1}^{\zeta}))^{([q]_1^{j-2}[0,q_d,0][q]_{j+1}^{\zeta})}}{ ([q!]_1^{j-2}[0!,q_d!,0!][q!]_{j+1}^{\zeta})  } \nonumber \\
&& \times \left(\prod\limits_{p=1}^{j-2}S_{k_p,i_p}\bm{Y}_p\right)\bm{P}_{k_{j-1},i_{j-1}}\bm{Y}_{j-1}\bm{P}_{k_c,i_c}\bm{N}_{k_d,i_d}^{q_d}\bm{Y}_j\bm{P}_{k_{j},i_{j}}\left(\prod\limits_{p'=j+1}^{\zeta}\bm{Y}_{p'}S_{k_{p'},i_{p'}}\right)\nonumber \\
&+&\sum\limits_{\substack{S_{k_p,i_p}, p \in [1,j-2], \\
S_{k_{p'},i_{p'}}, p' \in [j+1,\zeta]}}\sum\limits_{\substack{k_{j-1},i_{j-1},\\k_c,i_c, k_d,i_d,q_d\\ k_{j},i_{j}}} \frac{(\beta^{[\zeta]}([\lambda_k]_1^{j-2},\lambda_{k_{j-1}},\lambda_{k_c},\lambda_{k_j},[\lambda_k]_{j+1}^{\zeta}))^{([q]_1^{j-2}[0,q_d,0][q]_{j+1}^{\zeta})}}{ ([q!]_1^{j-2}[0!,q_d!,0!][q!]_{j+1}^{\zeta})  } \nonumber \\
&& \times \left(\prod\limits_{p=1}^{j-2}S_{k_p,i_p}\bm{Y}_p\right)\bm{P}_{k_{j-1},i_{j-1}}\bm{Y}_{j-1}\bm{P}_{k_c,i_c}\bm{N}_{k_d,i_d}^{q_d}\bm{Y}_j\bm{P}_{k_{j},i_{j}}\left(\prod\limits_{p'=j+1}^{\zeta}\bm{Y}_{p'}S_{k_{p'},i_{p'}}\right)\Bigg]\nonumber \\
&-&\sum\limits_{\substack{S_{k_p,i_p}, p \in [1,j-2], \\
S_{k_{p'},i_{p'}}, p' \in [j+1,\zeta]}}\sum\limits_{\substack{k_{j-1},i_{j-1},\\k_c,i_c, k_d,i_d,q_d\\ k_{j},i_{j}}} \frac{(\beta^{[\zeta]}([\lambda_k]_1^{j-2},\lambda_{k_{j-1}},\lambda_{k_c},\lambda_{k_j},[\lambda_k]_{j+1}^{\zeta}))^{([q]_1^{j-2}[0,q_d,0][q]_{j+1}^{\zeta})}}{ ([q!]_1^{j-2}[0!,q_d!,0!][q!]_{j+1}^{\zeta})  } \nonumber \\
&& \times \left(\prod\limits_{p=1}^{j-2}S_{k_p,i_p}\bm{Y}_p\right)\bm{P}_{k_{j-1},i_{j-1}}\bm{Y}_{j-1}\bm{P}_{k_c,i_c}\bm{N}_{k_d,i_d}^{q_d}\bm{Y}_j\bm{P}_{k_{j},i_{j}}\left(\prod\limits_{p'=j+1}^{\zeta}\bm{Y}_{p'}S_{k_{p'},i_{p'}}\right) \nonumber \\
&=_1& T_{\beta^{[\zeta+1]}}^{[\bm{X}]_1^{j-2},\bm{X}_{j-1,P},\bm{C}_N,\bm{D}_P,\bm{X}_{j,P},[\bm{X}]_{j+1}^{\zeta}}([\bm{Y}]_1^{j-1},\bm{C}-\bm{D},[\bm{Y}]_j^{\zeta})+ T_{\beta^{[\zeta+1]}}^{[\bm{X}]_1^{j-2},\bm{X}_{j-1,P},\bm{C}_P,\bm{D}_N,\bm{X}_{j,P},[\bm{X}]_{j+1}^{\zeta}}([\bm{Y}]_1^{j-1},\bm{C}-\bm{D},[\bm{Y}]_j^{\zeta}) \nonumber \\
&-&\mathfrak{X}^{\bm{X}_{j-1,P},\bm{C}_N,\bm{D}_N,\bm{X}_{j,P}},
\end{eqnarray}
}
where we apply Lemma~\ref{lma:first-order divided difference identity} in $=_1$ twice to get two GMOI terms with $\zeta+2$ parameter matrices. Moreover, the perturbation correct term $\mathfrak{X}^{\bm{X}_{j-1,P},\bm{C}_N,\bm{D}_N,\bm{X}_{j,P}}$ can be expressed by
\begin{eqnarray}\label{eq2:eq4-2:thm:GMOI Perturbation Formula Proof}
\lefteqn{\mathfrak{X}^{\bm{X}_{j-1,P},\bm{C}_N,\bm{D}_N,\bm{X}_{j,P}}=}\nonumber \\
&&\sum\limits_{\substack{S_{k_p,i_p}, p \in [1,j-2], \\
S_{k_{p'},i_{p'}}, p' \in [j+1,\zeta]}}\sum\limits_{\substack{k_{j-1},i_{j-1},\\k_c, i_c,q_c,k_d,i_d,q_d\\ k_{j},i_{j}}}\frac{(\beta^{[\zeta+1]}([\lambda_k]_1^{j-2},\lambda_{k_{j-1}},\lambda_{k_c},\lambda_{k_d},\lambda_{k_j},[\lambda_k]_{j+1}^{\zeta}))^{([q]_1^{j-2}[0,0,q_d,0][q]_{j+1}^{\zeta})}}{ ([q!]_1^{j-2}[0!,0!,q_d!,0!][q!]_{j+1}^{\zeta})  } \nonumber \\
&& \times \left(\prod\limits_{p=1}^{j-2}S_{k_p,i_p}\bm{Y}_p\right)\bm{P}_{k_{j-1},i_{j-1}}\bm{Y}_{j-1}\bm{P}_{k_c,i_c}(\bm{N}_{k_c,i_c}^{q_c}-\bm{N}_{k_d,i_d}^{q_d})\bm{N}_{k_d,i_d}^{q_d}\bm{Y}_j\bm{P}_{k_{j},i_{j}}\left(\prod\limits_{p'=j+1}^{\zeta}\bm{Y}_{p'}S_{k_{p'},i_{p'}}\right) \nonumber \\
&+& \sum\limits_{\substack{S_{k_p,i_p}, p \in [1,j-2], \\
S_{k_{p'},i_{p'}}, p' \in [j+1,\zeta]}}\sum\limits_{\substack{k_{j-1},i_{j-1},\\k_c, i_c,q_c,k_d,i_d,q_d\\ k_{j},i_{j}}}\frac{(\beta^{[\zeta+1]}([\lambda_k]_1^{j-2},\lambda_{k_{j-1}},\lambda_{k_c},\lambda_{k_d},\lambda_{k_j},[\lambda_k]_{j+1}^{\zeta}))^{([q]_1^{j-2}[0,q_c,0,0][q]_{j+1}^{\zeta})}}{ ([q!]_1^{j-2}[0!,q_c!,0!,0!][q!]_{j+1}^{\zeta})  } \nonumber \\
&& \times \left(\prod\limits_{p=1}^{j-2}S_{k_p,i_p}\bm{Y}_p\right)\bm{P}_{k_{j-1},i_{j-1}}\bm{Y}_{j-1}\bm{N}_{k_c,i_c}^{q_c}(\bm{N}_{k_c,i_c}^{q_c}-\bm{N}_{k_d,i_d}^{q_d})\bm{P}_{k_d,i_d}\bm{Y}_j\bm{P}_{k_{j},i_{j}}\left(\prod\limits_{p'=j+1}^{\zeta}\bm{Y}_{p'}S_{k_{p'},i_{p'}}\right) \nonumber \\
&+& \sum\limits_{\substack{S_{k_p,i_p}, p \in [1,j-2], \\
S_{k_{p'},i_{p'}}, p' \in [j+1,\zeta]}}\sum\limits_{\substack{k_{j-1},i_{j-1},\\k_c,i_c, k_d,i_d,q_d\\ k_{j},i_{j}}} \frac{(\beta^{[\zeta]}([\lambda_k]_1^{j-2},\lambda_{k_{j-1}},\lambda_{k_c},\lambda_{k_j},[\lambda_k]_{j+1}^{\zeta}))^{([q]_1^{j-2}[0,q_d,0][q]_{j+1}^{\zeta})}}{ ([q!]_1^{j-2}[0!,q_d!,0!][q!]_{j+1}^{\zeta})  } \nonumber \\
&& \times \left(\prod\limits_{p=1}^{j-2}S_{k_p,i_p}\bm{Y}_p\right)\bm{P}_{k_{j-1},i_{j-1}}\bm{Y}_{j-1}\bm{P}_{k_c,i_c}\bm{N}_{k_d,i_d}^{q_d}\bm{Y}_j\bm{P}_{k_{j},i_{j}}\left(\prod\limits_{p'=j+1}^{\zeta}\bm{Y}_{p'}S_{k_{p'},i_{p'}}\right) \nonumber \\  
&-& \sum\limits_{\substack{S_{k_p,i_p}, p \in [1,j-2], \\
S_{k_{p'},i_{p'}}, p' \in [j+1,\zeta]}}\sum\limits_{\substack{k_{j-1},i_{j-1},\\k_c, i_c, q_c, k_d,i_d\\ k_{j},i_{j}}} \frac{(\beta^{[\zeta]}([\lambda_k]_1^{j-2},\lambda_{k_{j-1}},\lambda_{k_d},\lambda_{k_j},[\lambda_k]_{j+1}^{\zeta}))^{([q]_1^{j-2}[0,q_c,0][q]_{j+1}^{\zeta})}}{ ([q!]_1^{j-2}[0!,q_c!,0!][q!]_{j+1}^{\zeta})  } \nonumber \\
&& \times \left(\prod\limits_{p=1}^{j-2}S_{k_p,i_p}\bm{Y}_p\right)\bm{P}_{k_{j-1},i_{j-1}}\bm{Y}_{j-1}\bm{N}_{k_c,i_c}^{q_c}\bm{P}_{k_d,i_d}\bm{Y}_j\bm{P}_{k_{j},i_{j}}\left(\prod\limits_{p'=j+1}^{\zeta}\bm{Y}_{p'}S_{k_{p'},i_{p'}}\right).
\end{eqnarray}
$\hfill\Box$

Following Example~\ref{exp:GMOI Perturbation Formula} will be provided by applying Theorem~\ref{thm:GMOI Perturbation Formula} to the GTOI. This example is an extension of Lemma 5 in~\cite{chang2025GDOIMatrix} as we need to require all parameter matrices without any nilpotent parts to have Lemma 5 in~\cite{chang2025GDOIMatrix} valid. 

\begin{example}\label{exp:GMOI Perturbation Formula}
We have
{\small
\begin{eqnarray}\label{eq1:exp:GMOI Perturbation Formula}
\lefteqn{T_{\beta^{[2]}}^{\bm{C},\bm{D},\bm{X}_1}(\bm{C}-\bm{D},\bm{Y})=}\nonumber \\
&& \Bigg(T_{\beta^{[1]}}^{\bm{C}_P,\bm{X}_{1,P}}(\bm{Y}) - T_{\beta^{[1]}}^{\bm{D}_P,\bm{X}_{1,P}}(\bm{Y})\Bigg) + \mathfrak{X}^{\bm{C}_P,\bm{D}_P,\bm{X}_{1,P}} \nonumber \\
&+&\Bigg(T_{\beta^{[1]}}^{\bm{C}_P,\bm{X}_{1,N}}(\bm{Y}) - T_{\beta^{[1]}}^{\bm{D}_P,\bm{X}_{1,N}}(\bm{Y})\Bigg) + \mathfrak{X}^{\bm{C}_P,\bm{D}_P,\bm{X}_{1,N}} \nonumber \\
&+& \Bigg(T_{\beta^{[1]}}^{\bm{C}_N,\bm{X}_{1,P}}(\bm{Y}) - T_{\beta^{[1]}}^{\bm{D}_N,\bm{X}_{1,P}}(\bm{Y})\Bigg) + \mathfrak{X}^{\bm{C}_N,\bm{D}_N,\bm{X}_{1,P}} \nonumber \\
&+&\Bigg(T_{\beta^{[1]}}^{\bm{C}_N,\bm{X}_{1,N}}(\bm{Y}) - T_{\beta^{[1]}}^{\bm{D}_N,\bm{X}_{1,N}}(\bm{Y})\Bigg) + \mathfrak{X}^{\bm{C}_N,\bm{D}_N,\bm{X}_{1,N}} \nonumber \\
&+&
T_{\beta^{[2]}}^{\bm{C}_N,\bm{D}_N,\bm{X}_{1,P}}(\bm{C}-\bm{D},\bm{Y})+ T_{\beta^{[2]}}^{\bm{C}_N,\bm{D}_N,\bm{X}_{1,N}}(\bm{C}-\bm{D},\bm{Y})\nonumber \\
&=&T_{\beta^{[1]}}^{\bm{C},\bm{X}_1}(\bm{Y}) - T_{\beta^{[1]}}^{\bm{D},\bm{X}_1}(\bm{Y})  + \bar{\mathfrak{X}}(\bm{C},\bm{D},\bm{X}_1)\nonumber \\
&+&
T_{\beta^{[2]}}^{\bm{C}_N,\bm{D}_N,\bm{X}_{1,P}}(\bm{C}-\bm{D},\bm{Y})+ T_{\beta^{[2]}}^{\bm{C}_N,\bm{D}_N,\bm{X}_{1,N}}(\bm{C}-\bm{D},\bm{Y}),
\end{eqnarray}
}
where $\beta^{[1]}$ and $\beta^{[2]}$ are first and second divide differences, and the term $\bar{\mathfrak{X}}$ is expressed by
\begin{eqnarray}\label{eq2:exp:GMOI Perturbation Formula}
\bar{\mathfrak{X}}(\bm{C},\bm{D},\bm{X}_1)&=&  \mathfrak{X}^{\bm{C}_P,\bm{D}_P,\bm{X}_{1,P}} +  \mathfrak{X}^{\bm{C}_P,\bm{D}_P,\bm{X}_{1,N}} \nonumber \\
&&+ \mathfrak{X}^{\bm{C}_N,\bm{D}_N,\bm{X}_{1,P}} + \mathfrak{X}^{\bm{C}_N,\bm{D}_N,\bm{X}_{1,N}}.
\end{eqnarray}

We also have the following expressions for terms $\mathfrak{X}^{\bm{C}_P,\bm{D}_P,\bm{X}_{1,P}}$, $\mathfrak{X}^{\bm{C}_P,\bm{D}_P,\bm{X}_{1,N}}$, $\mathfrak{X}^{\bm{C}_N,\bm{D}_N,\bm{X}_{1,P}}$, and $\mathfrak{X}^{\bm{C}_N,\bm{D}_N,\bm{X}_{1,N}}$. 
\begin{eqnarray}\label{e3-1:exp:GMOI Perturbation Formula}
\lefteqn{\mathfrak{X}^{\bm{C}_P,\bm{D}_P,\bm{X}_{1,P}}= }\nonumber \\
&&\sum\limits_{\substack{k_c, i_c,q_c,k_d,i_d,q_d\\ k_{1},i_{1}}}\frac{(\beta^{[2]}(\lambda_{k_c},\lambda_{k_d},\lambda_{k_1}))^{([0,0,0])}}{ ([0!,0!,0!])} \times \bm{P}_{k_c,i_c}(\bm{N}_{k_c,i_c}^{q_c}-\bm{N}_{k_d,i_d}^{q_d})\bm{P}_{k_d,i_d}\bm{Y}\bm{P}_{k_{1},i_{1}},
\end{eqnarray}
\begin{eqnarray}\label{e3-2:exp:GMOI Perturbation Formula}
\lefteqn{ \mathfrak{X}^{\bm{C}_P,\bm{D}_P,\bm{X}_{1,N}} = }\nonumber \\
&&\sum\limits_{\substack{k_c, i_c,q_c,k_d,i_d,q_d,\\ k_{1},i_{1},q_1}}\frac{(\beta^{[2]}(\lambda_{k_c},\lambda_{k_d},\lambda_{k_1}))^{([0,0,q_1])}}{ ([0!,0!,q_1!])} \times \bm{P}_{k_c,i_c}(\bm{N}_{k_c,i_c}^{q_c}-\bm{N}_{k_d,i_d}^{q_d})\bm{P}_{k_d,i_d}\bm{Y}\bm{N}_{k_{1},i_{1}}^{q_1},
\end{eqnarray}
\begin{eqnarray}\label{e3-3:exp:GMOI Perturbation Formula}
\lefteqn{ \mathfrak{X}^{\bm{C}_N,\bm{D}_N,\bm{X}_{1,P}} = }\nonumber \\
&&\sum\limits_{\substack{k_c, i_c,q_c,k_d,i_d,q_d,\\ k_{1},i_{1}}}\frac{(\beta^{[2]}(\lambda_{k_c},\lambda_{k_d},\lambda_{k_1}))^{([0,q_d,0])}}{ ([0!,q_d!,0!])} \times \bm{P}_{k_c,i_c}(\bm{N}_{k_c,i_c}^{q_c}-\bm{N}_{k_d,i_d}^{q_d})\bm{N}_{k_d,i_d}^{q_d}\bm{Y}\bm{P}_{k_{1},i_{1}}\nonumber \\
&+&\sum\limits_{\substack{k_c, i_c,q_c,k_d,i_d,q_d,\\ k_{1},i_{1}}}\frac{(\beta^{[2]}(\lambda_{k_c},\lambda_{k_d},\lambda_{k_1}))^{([q_c,0,0])}}{ ([q_c!,0!,0!])} \times \bm{N}_{k_c,i_c}^{q_c}(\bm{N}_{k_c,i_c}^{q_c}-\bm{N}_{k_d,i_d}^{q_d})\bm{P}_{k_d,i_d}\bm{Y}\bm{P}_{k_{1},i_{1}}\nonumber \\
&+&\sum\limits_{\substack{k_c, i_c,k_d,i_d,q_d,\\ k_{1},i_{1}}}\frac{(\beta^{[1]}(\lambda_{k_c},\lambda_{k_1}))^{([q_d,0])}}{ ([q_d!,0!])} \times \bm{P}_{k_c,i_c}\bm{N}_{k_d,i_d}^{q_d}\bm{Y}\bm{P}_{k_{1},i_{1}}\nonumber \\
&-&\sum\limits_{\substack{k_c, i_c,q_c,k_d,i_d,\\ k_{1},i_{1}}}\frac{(\beta^{[1]}(\lambda_{k_d},\lambda_{k_1}))^{([q_c,0])}}{ ([q_c!,0!])} \times \bm{N}_{k_c,i_c}^{q_c}\bm{P}_{k_d,i_d}\bm{Y}\bm{P}_{k_{1},i_{1}},
\end{eqnarray}
\begin{eqnarray}\label{e3-3:exp:GMOI Perturbation Formula}
\lefteqn{ \mathfrak{X}^{\bm{C}_N,\bm{D}_N,\bm{X}_{1,N}} = }\nonumber \\
&&\sum\limits_{\substack{k_c, i_c,q_c,k_d,i_d,q_d,\\ k_{1},i_{1},q_1}}\frac{(\beta^{[2]}(\lambda_{k_c},\lambda_{k_d},\lambda_{k_1}))^{([0,q_d,q_1])}}{ ([0!,q_d!,q_1!])} \times \bm{P}_{k_c,i_c}(\bm{N}_{k_c,i_c}^{q_c}-\bm{N}_{k_d,i_d}^{q_d})\bm{N}_{k_d,i_d}^{q_d}\bm{Y}\bm{N}_{k_{1},i_{1}}^{q_1}\nonumber \\
&+&\sum\limits_{\substack{k_c, i_c,q_c,k_d,i_d,q_d,\\ k_{1},i_{1},q_1}}\frac{(\beta^{[2]}(\lambda_{k_c},\lambda_{k_d},\lambda_{k_1}))^{([q_c,0,q_1])}}{ ([q_c!,0!,q_1!])} \times \bm{N}_{k_c,i_c}^{q_c}(\bm{N}_{k_c,i_c}^{q_c}-\bm{N}_{k_d,i_d}^{q_d})\bm{P}_{k_d,i_d}\bm{Y}\bm{N}_{k_{1},i_{1}}^{q_1}\nonumber \\
&+&\sum\limits_{\substack{k_c, i_c,k_d,i_d,q_d,\\ k_{1},i_{1},q_1}}\frac{(\beta^{[1]}(\lambda_{k_c},\lambda_{k_1}))^{([q_d,q_1])}}{ ([q_d!,q_1!])} \times \bm{P}_{k_c,i_c}\bm{N}_{k_d,i_d}^{q_d}\bm{Y}\bm{N}_{k_{1},i_{1}}^{q_1}\nonumber \\
&-&\sum\limits_{\substack{k_c, i_c,q_c,k_d,i_d,\\ k_{1},i_{1},q_1}}\frac{(\beta^{[1]}(\lambda_{k_d},\lambda_{k_1}))^{([q_c,q_1])}}{ ([q_c!,q_1!])} \times \bm{N}_{k_c,i_c}^{q_c}\bm{P}_{k_d,i_d}\bm{Y}\bm{N}_{k_{1},i_{1}}^{q_1}.
\end{eqnarray}
\end{example}

\subsection{GMOI Lipschitz Estimations}\label{sec: GMOI Lipschitz Estimations}

In this section, we will apply Theorem~\ref{thm: GMOI norm Est} to determine GMOI Lipschitz estimations. The Lipschitz estimation for GMOI is to determine the norm for the difference between GMOI $T_{\beta}^{\bm{X}_1,\ldots,\bm{X}_{\zeta+1}}(\bm{Y}_1,\ldots,\bm{Y}_{\zeta})$ and  GMOI $T_{\beta}^{\bm{X}_1,\ldots,\bm{X}_{\zeta+1}}(\bm{Y}'_1,\ldots,\bm{Y}'_{\zeta})$ is terms of norms for the difference $\bm{Y}_i$ and $\bm{Y}'_i$.  Theorem~\ref{thm: GMOI Lip Est} below is provided to characterize GMOI Lipschitz estimations.

\begin{theorem}\label{thm: GMOI Lip Est}
We have Lipschitz estimation for the difference between  GMOI $T_{\beta}^{\bm{X}_1,\ldots,\bm{X}_{\zeta+1}}(\bm{Y}_1,\ldots,\bm{Y}_{\zeta})$ and GMOI $T_{\beta}^{\bm{X}_1,\ldots,\bm{X}_{\zeta+1}}(\bm{Y}'_1,\ldots,\bm{Y}'_{\zeta})$:
\begin{eqnarray}\label{eq0:  thm: GMOI Lip Est}
\lefteqn{\left\Vert T_{\beta}^{\bm{X}_1,\ldots,\bm{X}_{\zeta+1}}(\bm{Y}_1,\ldots,\bm{Y}_{\zeta}) - T_{\beta}^{\bm{X}_1,\ldots,\bm{X}_{\zeta+1}}(\bm{Y}'_1,\ldots,\bm{Y}'_{\zeta}) \right\Vert}\nonumber \\
&\leq& 
\sum\limits_{i=1}^{\zeta}\Upsilon_i (\bm{X}_1,\ldots,\bm{X}_{\zeta+1},\beta) \left(\prod\limits_{j=1}^{i-1}\left\Vert\bm{Y}'_{j}\right\Vert\right)\left\Vert \bm{Y}_i - \bm{Y}'_i \right\Vert \left(\prod\limits_{j=i+1}^{\zeta}\left\Vert\bm{Y}_{j}\right\Vert\right),
\end{eqnarray}
where $\Upsilon_i (\bm{X}_1,\ldots,\bm{X}_{\zeta+1},\beta)$ is the scalar component of the upper bound for the following GMOI expressed by 
\begin{eqnarray}\label{eq0-1:  thm: GMOI Lip Est}
\lefteqn{\left\Vert T_{\beta}^{\bm{X}_1,\ldots,\bm{X}_{\zeta+1}}([\bm{Y}']_1^{i-1},\bm{Y}_i - \bm{Y}'_i,[\bm{Y}]_{i+1}^\zeta)  \right\Vert}\nonumber \\
&\leq_1&\sum\limits_{i'=0}^{2^{\zeta+1}-1} \left\Vert \bm{A}_{i'} \right\Vert_{up} \nonumber \\
&=& 
\Upsilon_i (\bm{X}_1,\ldots,\bm{X}_{\zeta+1},\beta)\left(\prod\limits_{j=1}^{i-1}\left\Vert\bm{Y}'_{j}\right\Vert\right)\left\Vert \bm{Y}_i - \bm{Y}'_i \right\Vert \left(\prod\limits_{j=i+1}^{\zeta}\left\Vert\bm{Y}_{j}\right\Vert\right).
\end{eqnarray}
The inequality $\leq_1$ comes from Theorem~\ref{thm: GMOI norm Est}. 
\end{theorem}
\textbf{Proof:}
We can perform a telescoping decomposition for the L.H.S. of Eq.~\eqref{eq0:  thm: GMOI Lip Est} as
\begin{eqnarray}\label{eq2:  thm: GMOI Lip Est}
\lefteqn{\left\Vert T_{\beta}^{\bm{X}_1,\ldots,\bm{X}_{\zeta+1}}(\bm{Y}_1,\ldots,\bm{Y}_{\zeta}) - T_{\beta}^{\bm{X}_1,\ldots,\bm{X}_{\zeta+1}}(\bm{Y}'_1,\ldots,\bm{Y}'_{\zeta}) \right\Vert}\nonumber \\
&=& \left\Vert \sum\limits_{i=1}^{\zeta}\left(T_{\beta}^{\bm{X}_1,\ldots,\bm{X}_{\zeta+1}}([\bm{Y}']_1^{i-1},[\bm{Y}]_i^\zeta) - T_{\beta}^{\bm{X}_1,\ldots,\bm{X}_{\zeta+1}}([\bm{Y}']_1^{i},[\bm{Y}]_{i+1}^\zeta)\right)\right\Vert
\nonumber \\
&=& \left\Vert \sum\limits_{i=1}^{\zeta}T_{\beta}^{\bm{X}_1,\ldots,\bm{X}_{\zeta+1}}([\bm{Y}']_1^{i-1},\bm{Y}_i - \bm{Y}'_i,[\bm{Y}]_{i+1}^\zeta)\right\Vert.
\end{eqnarray}
Then, this theorem is proved by applying Theorem~\ref{thm: GMOI norm Est} to obtain $\Upsilon_i (\bm{X}_1,\ldots,\bm{X}_{\zeta+1},\beta)$ given by Eq.~\eqref{eq0-1:  thm: GMOI Lip Est}.
$\hfill\Box$

\section{GMOI Continuity}\label{sec: GMOI Continuity}

The purpose of this section is to establish the continuity property of GMOI, i.e., we will have
\begin{eqnarray}\label{eq1:sec: GMOI Continuity}
 T_{\beta}^{\bm{X}_{1, \ell_1},\ldots,\bm{X}_{\zeta+1,\ell_{\zeta+1}}}(\bm{Y}_1,\ldots,\bm{Y}_{\zeta})\rightarrow T_{\beta}^{\bm{X}_1,\ldots,\bm{X}_{\zeta+1}}(\bm{Y}_1,\ldots,\bm{Y}_{\zeta}),
\end{eqnarray}
given that 
\begin{eqnarray}\label{eq2:sec: GMOI Continuity}
\bm{X}_{i, \ell_i} \rightarrow \bm{X}_{i}.
\end{eqnarray}
Note that $\rightarrow$ is in the sense of Frobenius norm. 

The following Lemma~\ref{lma: bar mathfrak X zero} has to be built before GMOI continuity property establishment. 
\begin{lemma}\label{lma: bar mathfrak X zero}
Given $\bm{C} \rightarrow \bm{D}$, we have
\begin{eqnarray}
\bar{\mathfrak{X}}([\bm{X}]_1^{j-1},\bm{C}, \bm{D}, [\bm{X}]_j^{\zeta})&\rightarrow& \bm{O}, 
\end{eqnarray}
where $\bm{O}$ is a matrix with all entries as zeros. 
\end{lemma} 
\textbf{Proof:}
From Eq.~\eqref{eq2:thm:GMOI Perturbation Formula}, we have
\begin{eqnarray}
\bar{\mathfrak{X}}([\bm{X}]_1^{j-1},\bm{C},\bm{D},[\bm{X}]_j^{\zeta})&=& \mathfrak{X}^{\bm{X}_{j-1,P},\bm{C}_P,\bm{D}_P,\bm{X}_{j,P}} + \mathfrak{X}^{\bm{X}_{j-1,P},\bm{C}_N,\bm{D}_N,\bm{X}_{j,P}} \nonumber \\
&+& \mathfrak{X}^{\bm{X}_{j-1,P},\bm{C}_P,\bm{D}_P,\bm{X}_{j,N}} + \mathfrak{X}^{\bm{X}_{j-1,P},\bm{C}_N,\bm{D}_N,\bm{X}_{j,N}} \nonumber \\
&+&  \mathfrak{X}^{\bm{X}_{j-1,N},\bm{C}_P,\bm{D}_P,\bm{X}_{j,P}} + \mathfrak{X}^{\bm{X}_{j-1,N},\bm{C}_N,\bm{D}_N,\bm{X}_{j,P}} \nonumber \\
&+& \mathfrak{X}^{\bm{X}_{j-1,N},\bm{C}_P,\bm{D}_P,\bm{X}_{j,N}} + \mathfrak{X}^{\bm{X}_{j-1,N},\bm{C}_N,\bm{D}_N,\bm{X}_{j,N}}. 
\end{eqnarray}

\textbf{Claim 1: $\mathfrak{X}^{\bm{X}_{j-1,P},\bm{C}_P,\bm{D}_P,\bm{X}_{j,P}}  \rightarrow \bm{O}$}

Let $\bm{C} = \bm{D} + t\bm{V}$, where $t \in \mathbb{R}$ and $\bm{V}$ is a fixed perturbation matrix. Then under smooth variation of the Jordan decomposition (i.e., assuming stable geometric multiplicities), we have:
\[
\mathfrak{X}^{\bm{X}_{j-1,P},\bm{C}_P,\bm{D}_P,\bm{X}_{j,P}} = O(t)
\quad \text{as} \quad t \to 0.
\]

We recall the expression:
{\small
\begin{align}
\lefteqn{\mathfrak{X}^{\bm{X}_{j-1,P},\bm{C}_P,\bm{D}_P,\bm{X}_{j,P}}=}\nonumber \\
&& 
\sum_{\substack{S_{k_p,i_p}, p \in [1,j-2], \\
S_{k_{p'},i_{p'}}, p' \in [j+1,\zeta]}} \sum_{\substack{k_{j-1},i_{j-1},\\k_c, i_c,q_c,k_d,i_d,q_d\\ k_{j},i_{j}}}
\frac{
\left( \beta^{[\zeta+1]}([\lambda_k]_1^{j-2},\lambda_{k_{j-1}},\lambda_{k_c},\lambda_{k_d},\lambda_{k_j},[\lambda_k]_{j+1}^{\zeta}) \right)^{([q]_1^{j-2}[0,0,0,0][q]_{j+1}^{\zeta})}
}{
([q!]_1^{j-2}[0!,0!,0!,0!][q!]_{j+1}^{\zeta})
} \nonumber \\
&& \times \left( \prod_{p=1}^{j-2} S_{k_p,i_p} \bm{Y}_p \right)
\bm{P}_{k_{j-1},i_{j-1}} \bm{Y}_{j-1}
\bm{P}_{k_c,i_c}(\bm{N}_{k_c,i_c}^{q_c} - \bm{N}_{k_d,i_d}^{q_d}) \bm{P}_{k_d,i_d} \bm{Y}_j \bm{P}_{k_j,i_j}
\left( \prod_{p'=j+1}^{\zeta} \bm{Y}_{p'} S_{k_{p'},i_{p'}} \right). 
\end{align}
}

We now study the behavior of the central difference:
\[
\bm{N}_{k_c,i_c}^{q_c}(\bm{C}) - \bm{N}_{k_d,i_d}^{q_d}(\bm{D}),
\]
where $\bm{N}_{k_c,i_c}^{q_c}(\bm{C})$ denotes the nilpotent part of the matrix $\bm{C}$ corresponding to the $k_c$-th eigenvalue and the $i_c$-th geometric component, with exponent $q_c$. Similarly, $\bm{N}_{k_d,i_d}^{q_d}(\bm{D})$ denotes the nilpotent part of the matrix $\bm{D}$ corresponding to the $k_d$-th eigenvalue and the $i_d$-th geometric component, with exponent $q_d$.

Let $\bm{C} = \bm{D} + t \bm{V}$ and suppose the eigenvalues of $\bm{C}$ and $\bm{D}$ satisfy
\[
\lambda_{k_c} = \lambda_{k_d} + t\mu + o(t),
\]
with corresponding Jordan subspaces continuously varying. Then by perturbation theory for Jordan blocks:
\[
\bm{N}_{k_c,i_c}^{q_c}(\bm{C}) = \bm{N}_{k_d,i_d}^{q_d}(\bm{D}) + t \cdot \dot{\bm{N}}_{k_d,i_d}^{q_c} + o(t),
\]
and
\[
\bm{P}_{k_c,i_c}(\bm{C}) = \bm{P}_{k_d,i_d}(\bm{D}) + t \cdot \dot{\bm{P}}_{k_d,i_d} + o(t),
\]
where $\dot{\bm{N}}_{k_d,i_d}$ and $\dot{\bm{P}}_{k_d,i_d}$ are Fréchet derivatives with respect to $t$.

Thus, the difference term becomes:
\[
\bm{P}_{k_c,i_c}(\bm{N}_{k_c,i_c}^{q_c} - \bm{N}_{k_d,i_d}^{q_d})\bm{P}_{k_d,i_d}
= t \cdot \left( \bm{P}_{k_d,i_d} \dot{\bm{N}}_{k_d,i_d}^{q_c} \bm{P}_{k_d,i_d} \right) + o(t).
\]

Substituting this into the full sum, each term is of the form:
\[
t \cdot \left( \text{bounded coefficient involving } \beta^{[\zeta+1]} \text{ and } \bm{Y}_p, \bm{P}_{k,i}, \bm{S}_{k,i} \right) + o(t),
\]
which implies:
\[
\mathfrak{X}^{\bm{X}_{j-1,P},\bm{C}_P,\bm{D}_P,\bm{X}_{j,P}} = O(t) \quad \text{as } t \to 0.
\]

\textbf{Claim 2: $\mathfrak{X}^{\bm{X}_{j-1,P},\bm{C}_N,\bm{D}_N,\bm{X}_{j,P}}  \rightarrow \bm{O}$}

Let $\bm{C} = \bm{D} + t\bm{V}$ for small $t > 0$ and fixed matrix $\bm{V}$, and suppose the Jordan decomposition of $\bm{C}$ depends smoothly on $t$. Then:
\[
\mathfrak{X}^{\bm{X}_{j-1,P},\bm{C}_N,\bm{D}_N,\bm{X}_{j,P}} = O(t)
\quad \text{as} \quad t \to 0.
\]

The expression for $\mathfrak{X}$ consists of four multilinear terms in Eq.~\eqref{eq2:cross:thm:GMOI Perturbation Formula}, each involving differences between nilpotent components of $\bm{C}$ and $\bm{D}$. 

Again, by perturbation theory for Jordan blocks:
\[
\bm{N}_{k_c,i_c}^{q_c}(\bm{C}) = \bm{N}_{k_d,i_d}^{q_d}(\bm{D}) + t \cdot \dot{\bm{N}}_{k_d,i_d}^{q_c} + o(t),
\]
and
\[
\bm{P}_{k_c,i_c}(\bm{C}) = \bm{P}_{k_d,i_d}(\bm{D}) + t \cdot \dot{\bm{P}}_{k_d,i_d} + o(t),
\]
where $\dot{\bm{N}}_{k_d,i_d}$ and $\dot{\bm{P}}_{k_d,i_d}$ are Fréchet derivatives with respect to $t$.

Now analyze each of the four terms in Eq.~\eqref{eq2:cross:thm:GMOI Perturbation Formula}:

Term 1:
\[
\bm{P}_{k_c,i_c}(\bm{N}_{k_c,i_c}^{q_c} - \bm{N}_{k_d,i_d}^{q_d})\bm{N}_{k_d,i_d}^{q_d}
= O(t),
\]
since $(\bm{N}_{k_c,i_c}^{q_c} - \bm{N}_{k_d,i_d}^{q_d}) = O(t)$ and all other matrices (e.g., $\bm{Y}_j$) are fixed.

Term 2:
\[
\bm{N}_{k_c,i_c}^{q_c}(\bm{N}_{k_c,i_c}^{q_c} - \bm{N}_{k_d,i_d}^{q_d})\bm{P}_{k_d,i_d}
= O(t),
\]
because $(\bm{N}_{k_c,i_c}^{q_c} - \bm{N}_{k_d,i_d}^{q_d}) = O(t)$ is again an $O(t)$ difference.

We prove that the third and fourth terms in Eq.~\eqref{eq2:cross:thm:GMOI Perturbation Formula} cancel in the limit $\bm{C} \to \bm{D}$. Denote the third term as $T_3$ and the fourth term as $T_4$. Then, we have:
\begin{eqnarray}
T_3 &=& \sum\limits_{\substack{S_{k_p,i_p}, p \in [1,j-2], \\
S_{k_{p'},i_{p'}}, p' \in [j+1,\zeta]}}\sum\limits_{\substack{k_{j-1},i_{j-1},\\k_c,i_c, k_d,i_d,q_d\\ k_{j},i_{j}}} \frac{(\beta^{[\zeta]}([\lambda_k]_1^{j-2},\lambda_{k_{j-1}},\lambda_{k_c},\lambda_{k_j},[\lambda_k]_{j+1}^{\zeta}))^{([q]_1^{j-2}[0,q_d,0][q]_{j+1}^{\zeta})}}{ ([q!]_1^{j-2}[0!,q_d!,0!][q!]_{j+1}^{\zeta})  } \nonumber \\
&& \times \left(\prod\limits_{p=1}^{j-2}S_{k_p,i_p}\bm{Y}_p\right)\bm{P}_{k_{j-1},i_{j-1}}\bm{Y}_{j-1}\bm{P}_{k_c,i_c}\bm{N}_{k_d,i_d}^{q_d}\bm{Y}_j\bm{P}_{k_{j},i_{j}}\left(\prod\limits_{p'=j+1}^{\zeta}\bm{Y}_{p'}S_{k_{p'},i_{p'}}\right),
\end{eqnarray}

\begin{eqnarray}
T_4 &=& \sum\limits_{\substack{S_{k_p,i_p}, p \in [1,j-2], \\
S_{k_{p'},i_{p'}}, p' \in [j+1,\zeta]}}\sum\limits_{\substack{k_{j-1},i_{j-1},\\k_c, i_c, q_c, k_d,i_d\\ k_{j},i_{j}}} \frac{(\beta^{[\zeta]}([\lambda_k]_1^{j-2},\lambda_{k_{j-1}},\lambda_{k_d},\lambda_{k_j},[\lambda_k]_{j+1}^{\zeta}))^{([q]_1^{j-2}[0,q_c,0][q]_{j+1}^{\zeta})}}{ ([q!]_1^{j-2}[0!,q_c!,0!][q!]_{j+1}^{\zeta})  } \nonumber \\
&& \times \left(\prod\limits_{p=1}^{j-2}S_{k_p,i_p}\bm{Y}_p\right)\bm{P}_{k_{j-1},i_{j-1}}\bm{Y}_{j-1}\bm{N}_{k_c,i_c}^{q_c}\bm{P}_{k_d,i_d}\bm{Y}_j\bm{P}_{k_{j},i_{j}}\left(\prod\limits_{p'=j+1}^{\zeta}\bm{Y}_{p'}S_{k_{p'},i_{p'}}\right)
\end{eqnarray}

Now, observe the following:

\begin{enumerate}
    \item As $\bm{C} \to \bm{D}$, we have $\lambda_{k_c} \to \lambda_{k_d}$ and the generalized divided differences $\beta^{[\zeta]}$ converge:
    \[
    \beta^{[\zeta]}(\cdots, \lambda_{k_c}, \cdots) \to \beta^{[\zeta]}(\cdots, \lambda_{k_d}, \cdots).
    \]
    
    \item The projectors and nilpotents satisfy:
    \[
    \bm{P}_{k_c,i_c} \to \bm{P}_{k_d,i_d}, \quad \bm{N}_{k_c,i_c}^{q_c} \to \bm{N}_{k_d,i_d}^{q_d}.
    \]
    
    \item Furthermore, in the Jordan canonical basis, the nilpotent and projector commute:
    \[
    \bm{P}_{k_d,i_d} \bm{N}_{k_d,i_d}^{q_d} = \bm{N}_{k_d,i_d}^{q_d} \bm{P}_{k_d,i_d}.
    \]
\end{enumerate}

Combining these, we find that in the limit $\bm{C} \to \bm{D}$:

\[
T_3 \to \sum \beta(\cdots) \bm{P}_{k_d,i_d} \bm{N}_{k_d,i_d}^{q_d}, \qquad 
T_4 \to \sum \beta(\cdots) \bm{N}_{k_d,i_d}^{q_d} \bm{P}_{k_d,i_d}.
\]

Since these matrices commute in the limit and the coefficients match, we conclude:
\[
T_3 - T_4 \to \bm{0}.
\]

Therefore, combining the asymptotic bounds and noting all differences are \(O(t)\), we conclude:
\[
\mathfrak{X}^{\bm{X}_{j-1,P},\bm{C}_N,\bm{D}_N,\bm{X}_{j,P}} = O(t) \to \bm{O} \quad \text{as } \bm{C} \to \bm{D}.
\]

This lemma can be proved by using the same methods from proving Claim 1 and Claim 2 to show the following
\begin{eqnarray}
\mathfrak{X}^{\bm{X}_{j-1,P},\bm{C}_P,\bm{D}_P,\bm{X}_{j,N}}  \rightarrow \bm{O}, \nonumber \\
\mathfrak{X}^{\bm{X}_{j-1,P},\bm{C}_N,\bm{D}_N,\bm{X}_{j,N}} \rightarrow \bm{O}, \nonumber \\\mathfrak{X}^{\bm{X}_{j-1,N},\bm{C}_P,\bm{D}_P,\bm{X}_{j,P}} \rightarrow \bm{O}, \nonumber \\
\mathfrak{X}^{\bm{X}_{j-1,N},\bm{C}_N,\bm{D}_N,\bm{X}_{j,P}} \rightarrow \bm{O}, \nonumber \\\mathfrak{X}^{\bm{X}_{j-1,N},\bm{C}_P,\bm{D}_P,\bm{X}_{j,N}} \rightarrow \bm{O}, \nonumber \\ \mathfrak{X}^{\bm{X}_{j-1,N},\bm{C}_N,\bm{D}_N,\bm{X}_{j,N}} \rightarrow \bm{O}.
\end{eqnarray}
$\hfill\Box$

We are ready to present Theorem~\ref{thm:GMOI continuity}, which is used to establish the GMOI continuity property. 
\begin{theorem}\label{thm:GMOI continuity}
Suppose that
\begin{equation}\label{eq1:thm:GMOI continuity}
\bm{X}_{j, \ell_j} \rightarrow \bm{X}_{j},
\end{equation}
for all $j = 1, 2, \ldots, \zeta+1$, and assume that all partial derivatives of the function $\beta^{[\zeta+1]}$ up to any order are bounded. Then, it follows that
\begin{eqnarray}\label{eq2:thm:GMOI continuity}
 T_{\beta^{[\zeta]}}^{\bm{X}_{1, \ell_1},\ldots,\bm{X}_{\zeta+1,\ell_{\zeta+1}}}(\bm{Y}_1,\ldots,\bm{Y}_{\zeta})\rightarrow T_{\beta^{[\zeta]}}^{\bm{X}_1,\ldots,\bm{X}_{\zeta+1}}(\bm{Y}_1,\ldots,\bm{Y}_{\zeta}).
\end{eqnarray}
\end{theorem}
\textbf{Proof:}
By telescoping summation formula, we have
\begin{eqnarray}\label{eq3:thm:GMOI continuity}
\lefteqn{\left\Vert T_{\beta^{[\zeta]}}^{\bm{X}_{1, \ell_1},\ldots,\bm{X}_{\zeta+1,\ell_{\zeta+1}}}(\bm{Y}_1,\ldots,\bm{Y}_{\zeta}) - T_{\beta^{[\zeta]}}^{\bm{X}_1,\ldots,\bm{X}_{\zeta+1}}(\bm{Y}_1,\ldots,\bm{Y}_{\zeta}) \right\Vert}\nonumber \\ &=&  \left\Vert \sum\limits_{j=1}^{\zeta+1}\left(T_{\beta^{[\zeta]}}^{[\bm{X}]_1^{j-1},[\bm{X}_\ell]_j^{\zeta+1}}(\bm{Y}_1,\ldots,\bm{Y}_{\zeta}) - T_{\beta^{[\zeta]}}^{ [\bm{X}]_1^{j},[\bm{X}_\ell]_{j+1}^{\zeta+1} }(\bm{Y}_1,\ldots,\bm{Y}_{\zeta})\right)\right\Vert.
\end{eqnarray}
For each index $j$, by applying GMOI perturbation formula given by Theorem~\ref{thm:GMOI Perturbation Formula}, we have 
\begin{eqnarray}\label{eq4:thm:GMOI continuity}
\lefteqn{\left\Vert T_{\beta^{[\zeta]}}^{[\bm{X}]_1^{j-1},[\bm{X}_\ell]_j^{\zeta+1}}(\bm{Y}_1,\ldots,\bm{Y}_{\zeta}) - T_{\beta^{[\zeta]}}^{ [\bm{X}]_1^{j},[\bm{X}_\ell]_{j+1}^{\zeta+1} }(\bm{Y}_1,\ldots,\bm{Y}_{\zeta} \right\Vert}\nonumber \\ &=&  \left\Vert T_{\beta^{[\zeta+1]}}^{[\bm{X}]_1^{j-1},\bm{X}_{j, \ell_j},\bm{X}_{j},[\bm{X}_\ell]_{j+1}^{\zeta+1}}([\bm{Y}]_1^{j-1},\bm{X}_{j, \ell_j}-\bm{X}_{j},[\bm{Y}]_j^{\zeta}) -
\bar{\mathfrak{X}}([\bm{X}]_1^{j-1},\bm{X}_{j, \ell_j},\bm{X}_{j},[\bm{X}_\ell]_{j+1}^{\zeta+1}) \right. \nonumber \\
&& - \left. T_{\beta^{[\zeta+1]}}^{[\bm{X}]_1^{j-2},\bm{X}_{j-1,P},\bm{X}_{j, \ell_{j},N},\bm{X}_{j,N},\bm{X}_{j+1,\ell_{j+1}, P},[\bm{X}_\ell]_{j+2}^{\zeta+1}}([\bm{Y}]_1^{j-1},\bm{X}_{j, \ell_j}-\bm{X}_{j},[\bm{Y}]_j^{\zeta}) \right. \nonumber \\
&& - \left. T_{\beta^{[\zeta+1]}}^{[\bm{X}]_1^{j-2},\bm{X}_{j-1,P},\bm{X}_{j, \ell_{j}, N},\bm{X}_{j,N},\bm{X}_{j+1,\ell_{j+1},N},[\bm{X}]_{j+2}^{\zeta+1}}([\bm{Y}]_1^{j-1},\bm{X}_{j, \ell_j}-\bm{X}_{j},[\bm{Y}]_j^{\zeta}) \right. \nonumber \\
&& - \left. T_{\beta^{[\zeta+1]}}^{[\bm{X}]_1^{j-2},\bm{X}_{j-1,N},\bm{X}_{j, \ell_{j},N},\bm{X}_{j,N},\bm{X}_{j+1,\ell_{j+1}, P},[\bm{X}_\ell]_{j+2}^{\zeta+1}}([\bm{Y}]_1^{j-1},\bm{X}_{j, \ell_j}-\bm{X}_{j},[\bm{Y}]_j^{\zeta}) \right. \nonumber \\
&& - \left. T_{\beta^{[\zeta+1]}}^{[\bm{X}]_1^{j-2},\bm{X}_{j-1,N},\bm{X}_{j, \ell_{j}, N},\bm{X}_{j,N},\bm{X}_{j+1,\ell_{j+1},N},[\bm{X}]_{j+2}^{\zeta+1}}([\bm{Y}]_1^{j-1},\bm{X}_{j, \ell_j}-\bm{X}_{j},[\bm{Y}]_j^{\zeta})   \right\Vert  \nonumber \\
&\leq& \left\Vert T_{\beta^{[\zeta+1]}}^{[\bm{X}]_1^{j-1},\bm{X}_{j, \ell_j},\bm{X}_{j},[\bm{X}_\ell]_{j+1}^{\zeta}}([\bm{Y}]_1^{j-1},\bm{X}_{j, \ell_j}-\bm{X}_{j},[\bm{Y}]_j^{\zeta}) \right\Vert + \nonumber \\ 
&& + \left\Vert \bar{\mathfrak{X}}([\bm{X}]_1^{j-1},\bm{X}_{j, \ell_j},\bm{X}_{j},[\bm{X}_\ell]_{j+1}^{\zeta+1}) \right\Vert  \nonumber \\
&& + \left\Vert T_{\beta^{[\zeta+1]}}^{[\bm{X}]_1^{j-2},\bm{X}_{j-1,P},\bm{X}_{j, \ell_{j},N},\bm{X}_{j,N},\bm{X}_{j+1,\ell_{j+1}, P},[\bm{X}_\ell]_{j+2}^{\zeta}}([\bm{Y}]_1^{j-1},\bm{X}_{j, \ell_j}-\bm{X}_{j},[\bm{Y}]_j^{\zeta}) \right\Vert \nonumber \\
&& + \left\Vert T_{\beta^{[\zeta+1]}}^{[\bm{X}]_1^{j-2},\bm{X}_{j-1,P},\bm{X}_{j, \ell_{j}, N},\bm{X}_{j,N},\bm{X}_{j+1,\ell_{j+1},N},[\bm{X}]_{j+2}^{\zeta}}([\bm{Y}]_1^{j-1},\bm{X}_{j, \ell_j}-\bm{X}_{j},[\bm{Y}]_j^{\zeta}) \right\Vert \nonumber \\ 
&& + \left\Vert T_{\beta^{[\zeta+1]}}^{[\bm{X}]_1^{j-2},\bm{X}_{j-1,N},\bm{X}_{j, \ell_{j}, N},\bm{X}_{j,N},\bm{X}_{j+1,\ell_{j+1},P},[\bm{X}]_{j+2}^{\zeta}}([\bm{Y}]_1^{j-1},\bm{X}_{j, \ell_j}-\bm{X}_{j},[\bm{Y}]_j^{\zeta}) \right\Vert \nonumber \\ 
&& + \left\Vert T_{\beta^{[\zeta+1]}}^{[\bm{X}]_1^{j-2},\bm{X}_{j-1,N},\bm{X}_{j, \ell_{j}, N},\bm{X}_{j,N},\bm{X}_{j+1,\ell_{j+1},N},[\bm{X}]_{j+2}^{\zeta}}([\bm{Y}]_1^{j-1},\bm{X}_{j, \ell_j}-\bm{X}_{j},[\bm{Y}]_j^{\zeta}) \right\Vert \nonumber \\ 
&\leq_1& \frac{\epsilon}{6}+ \left(\frac{\epsilon}{6}\right)_\star +  \frac{\epsilon}{6}+ \frac{\epsilon}{6} +  \frac{\epsilon}{6}+ \frac{\epsilon}{6}=\epsilon,
\end{eqnarray}
where each $\frac{\epsilon}{6}$ without $\star$ comes from Theorem~\ref{thm: GMOI norm Est} and the assumption that all partial derivatives of the function $\beta^{[\zeta+1]}$ up to any order are bounded. Besides, the term $ (\frac{\epsilon}{6})_\star$ comes from Lemma~\ref{lma: bar mathfrak X zero}.

Then, this theorem is proved because each summand in Eq.~\eqref{eq3:thm:GMOI continuity} approaches to  zero from Eq.~\eqref{eq4:thm:GMOI continuity}.
$\hfill\Box$

\section{Applications: Differentiation of Matrix Functions}\label{sec: Applications: Differentiation of Matrix Functions}

The purpose of this section is to apply previous GMOI theory to study the differentiation of a matrix function.  The goal is to obtain 
\begin{eqnarray}\label{eq:n-th derivative}
\frac{ d^n f (\bm{X}+t \bm{Y})}{d t^n},
\end{eqnarray}
where $n \in \mathbb{N}$. Recall Theorem 11~\cite{chang2025GDOICont}, we have
\begin{theorem}\label{thm:d f X dt}
Let $f$ be a complex-valued and differentiable function, we have 
\begin{eqnarray}\label{eq1:thm:d f X dt}
\frac{df(\bm{X} + t \bm{Y})}{dt}&=&T_{f^{[1]}}^{\bm{X} + t \bm{Y},\bm{X} + t \bm{Y}}\left(\bm{Y}\right),
\end{eqnarray}
where $T_{f^{[1]}}^{\bm{X} + t \bm{Y},\bm{X} + t \bm{Y}}\left(\bm{Y}\right)$ is the GDOI with matrices $\bm{X}$ and $\bm{Y}$.   
\end{theorem}

Before presenting the main theorem in this section, we have to introduce the following lemmas. Lemma~\ref{lma:diff of GMOI with all para matrices are nilpotent}  shows tha the differentiation of a GMOI will be a zero matrix if all parameter matrices are nilpotent parts of the original parameter matrices. 

\begin{lemma}\label{lma:diff of GMOI with all para matrices are nilpotent} 
Given a GMOI with the following format:
\begin{eqnarray}\label{eq1:lma:diff of GMOI with all para matrices are nilpotent} 
T_{\beta^{[\zeta]}}^{[\bm{X}_N]_{1}^{\zeta+1}}([\bm{Y}]_1^{\zeta}), 
\end{eqnarray}
we have
\begin{eqnarray}\label{eq2:lma:diff of GMOI with all para matrices are nilpotent} 
\left.\frac{d T_{\beta^{[\zeta]}}^{[(\bm{X}+t\bm{Y})_N]_{1}^{\zeta+1}}([\bm{Y}]_1^{\zeta})}{dt}\right\vert_{t=0} = \bm{O}. 
\end{eqnarray}
\end{lemma}
\textbf{Proof:}
By derivative definition and telescoping summation relation, we have
\begin{eqnarray}\label{eq3:lma:diff of GMOI with all para matrices are nilpotent} 
\lefteqn{\lim\limits_{t \rightarrow 0} \frac{ T_{\beta^{[\zeta]}}^{[(\bm{X}+t\bm{Y})_N]_{1}^{\zeta+1}}([\bm{Y}]_1^{\zeta}) - T_{\beta^{[\zeta]}}^{[\bm{X}_N]_{1}^{\zeta+1}}([\bm{Y}]_1^{\zeta}) }{t}=}\nonumber \\
&& \lim\limits_{t \rightarrow 0} \frac{\sum\limits_{j=1}^{\zeta+1}\left(T_{\beta^{[\zeta]}}^{[\bm{X}_N]_{1}^{j-1},[(\bm{X}+t\bm{Y})_N]_{j}^{\zeta+1}}([\bm{Y}]_1^{\zeta}) - T_{\beta^{[\zeta]}}^{[\bm{X}_N]_{1}^{j},[(\bm{X}+t\bm{Y})_N]_{j+1}^{\zeta+1}}([\bm{Y}]_1^{\zeta})\right)}{t}.
\end{eqnarray}
Then, for each $j$ in the summand of the numerator of Eq.~\eqref{eq3:lma:diff of GMOI with all para matrices are nilpotent}, we have
\begin{eqnarray}\label{eq4:lma:diff of GMOI with all para matrices are nilpotent} 
\lefteqn{T_{\beta^{[\zeta]}}^{[\bm{X}_N]_{1}^{j-1},[(\bm{X}+t\bm{Y})_N]_{j}^{\zeta+1}}([\bm{Y}]_1^{\zeta}) - T_{\beta^{[\zeta]}}^{[\bm{X}_N]_{1}^{j},[(\bm{X}+t\bm{Y})_N]_{j+1}^{\zeta+1}}([\bm{Y}]_1^{\zeta})}\nonumber \\
&=_1& T_{\beta^{[\zeta+1]}}^{[\bm{X}_N]_{1}^{j-1},(\bm{X}+t\bm{Y})_N,\bm{X}_N,[(\bm{X}+t\bm{Y})_N]_{j}^{\zeta+1}}([\bm{Y}]_1^{j-1},t\bm{Y},[\bm{Y}]_{j}^{\zeta}) \nonumber \\
&& - T_{\beta^{[\zeta+1]}}^{[\bm{X}_N]_{1}^{j-1},(\bm{X}+t\bm{Y})_N,\bm{X}_N,[(\bm{X}+t\bm{Y})_N]_{j}^{\zeta+1}}([\bm{Y}]_1^{j-1},t\bm{Y},[\bm{Y}]_{j}^{\zeta}) \nonumber \\
&& -  \overline{\mathfrak{X}}([\bm{X}_N]_{1}^{j-1},(\bm{X}+t\bm{Y})_N,\bm{X}_N,[(\bm{X}+t\bm{Y})_N]_{j}^{\zeta+1})\nonumber \\
&=_2&\bm{O},
\end{eqnarray}
where the perturbation formula given by Theorem~\ref{thm:GMOI Perturbation Formula} is applied to get $=_1$, and Eq.~\eqref{eq2:thm:GMOI Perturbation Formula} is used to get $=_2$ since $ \overline{\mathfrak{X}}([\bm{X}_N]_{1}^{j-1},(\bm{X}+t\bm{Y})_N,\bm{X}_N,[(\bm{X}+t\bm{Y})_N]_{j}^{\zeta+1}) = \bm{O}$ by setting $\bm{C} = (\bm{X}+t\bm{Y})_N$ and $\bm{D}=\bm{X}_N$, respectively. 
$\hfill\Box$

Lemma~\ref{lma:diff of GMOI with para matrices having some org matrices} will study the the differentiation of a GMOI given that not all parameter matrices are nilpotent parts of the original parameter matrices. 

\begin{lemma}\label{lma:diff of GMOI with para matrices having some org matrices}
Given a GMOI with the following format:
\begin{eqnarray}\label{eq1:lma:diff of GMOI with para matrices having some org matrices}
T_{\beta^{[\zeta]}}^{\underline{[\bm{X},\bm{X}_N]}_\ell}([\bm{Y}]_1^{\zeta}), 
\end{eqnarray}
where $\underline{[\bm{X},\bm{X}_N]}_\ell$ is an array of parameter matrices with $\zeta+1$ matrices composed by at least one matrix $\bm{X}$ and at least one matrix $\bm{X}_N$, and the subscript $\ell$ is used to index those many possible arrays of parameter matrices with $\zeta+1$ matrices composed by at least one matrix $\bm{X}$ and at least one matrix $\bm{X}_N$.  We define the mapping function $\vartheta_\ell(j)$ as the position of the $j$-th paramter matrix $\bm{X}$ in the array $\underline{[\bm{X},\bm{X}_N]}_\ell$,  where $j=1,2,\ldots, \left\vert\underline{[\bm{X},\bm{X}_N]}_\ell\right\vert_{\bm{X}}$ and $\left\vert\underline{[\bm{X},\bm{X}_N]}_\ell\right\vert_{\bm{X}}$ represents the number of $\bm{X}$ in the array of parameter matrices $\underline{[\bm{X},\bm{X}_N]}_\ell$. 
Let $[\underline{[\bm{X},\bm{X}_N]}_\ell]_1^{\vartheta_\ell(j) -1}$ represents the first $\vartheta_\ell(j) -1$ parameter matrices of the array $\underline{[\bm{X},\bm{X}_N]}_\ell$, and $[\underline{[\bm{X},\bm{X}_N]}_\ell]_{\vartheta_\ell(j)+1}^{\zeta+1}$ represents the last $\zeta-\vartheta_\ell(j) +1$ parameter matrices of the array $\underline{[\bm{X},\bm{X}_N]}_\ell$. We also use $\tilde{\bm{X}} = \bm{X}+ t \bm{Y}$ for later notations simplification.

Then, we have
\begin{eqnarray}\label{eq2:lma:diff of GMOI with para matrices having some org matrices}
\lefteqn{\left.\frac{d T_{\beta^{[\zeta]}}^{\underline{[\bm{X}+t\bm{Y},(\bm{X}+t\bm{Y})_N]}_\ell}([\bm{Y}]_1^{\zeta})}{dt}\right\vert_{t=0}}\nonumber \\
&=& \sum\limits_{j=1}^{\left\vert\underline{[\bm{X},\bm{X}_N]}_\ell\right\vert_{\bm{X}}}\Bigg[T_{\beta^{[\zeta+1]}}^{
[\underline{[\bm{X},\bm{X}_N]}_\ell]_1^{\vartheta_\ell(j) -1},\bm{X},\bm{X},
[\underline{[\bm{X},\bm{X}_N]}_\ell]_{\vartheta_\ell(j)+1}^{\zeta+1}
}([\bm{Y}]_1^{\vartheta_\ell(j)-1}, \bm{Y}, [\bm{Y}]_{\vartheta_\ell(j)}^{\zeta}) \nonumber \\
&& -  T_{\beta^{[\zeta+1]}}^{
[\underline{[\bm{X},\bm{X}_N]}_\ell]_1^{\vartheta_\ell(j) -1},\bm{X}_N,\bm{X}_N,
[\underline{[\bm{X},\bm{X}_N]}_\ell]_{\vartheta_\ell(j)+1}^{\zeta+1}
}([\bm{Y}]_1^{\vartheta_\ell(j)-1}, \bm{Y}, [\bm{Y}]_{\vartheta_\ell(j)}^{\zeta}) \nonumber \\
&& - \left.\mathfrak{X}^{(1)}\left(
[\underline{[\bm{X},\bm{X}_N]}_\ell]_1^{\vartheta_\ell(j) -1},\tilde{\bm{X}},\bm{X},
[\underline{[\tilde{\bm{X}},\tilde{\bm{X}}_N]}_\ell]_{\vartheta_\ell(j)+1}^{\zeta+1}
\right)\right\vert_{t=0} \Bigg]
\end{eqnarray}
where $ \left.\mathfrak{X}^{(1)}\left(
[\underline{[\bm{X},\bm{X}_N]}_\ell]_1^{\vartheta_\ell(j) -1},\tilde{\bm{X}},\bm{X},
[\underline{[\tilde{\bm{X}},\tilde{\bm{X}}_N]}_\ell]_{\vartheta_\ell(j)+1}^{\zeta+1}
\right)\right\vert_{t=0}$ is obtained by
\begin{eqnarray}\label{eq3:lma:diff of GMOI with para matrices having some org matrices}
\lefteqn{\left.\mathfrak{X}^{(1)}\left([\underline{[\bm{X},\bm{X}_N]}_\ell]_1^{\vartheta_\ell(j) -1},\tilde{\bm{X}},\bm{X},
[\underline{[\tilde{\bm{X}},\tilde{\bm{X}}_N]}_\ell]_{\vartheta_\ell(j)+1}^{\zeta+1}\right)\right\vert_{t=0}}\nonumber \\
&=&
\left.\frac{d \mathfrak{X}\left( [\underline{[\bm{X},\bm{X}_N]}_\ell]_1^{\vartheta_\ell(j) -1},\tilde{\bm{X}},\bm{X},
[\underline{[\tilde{\bm{X}},\tilde{\bm{X}}_N]}_\ell]_{\vartheta_\ell(j)+1}^{\zeta+1} \right)}{dt} \right\vert_{t=0}
\end{eqnarray}
\end{lemma}
\textbf{Proof:}
From Lemma~\ref{lma:diff of GMOI with all para matrices are nilpotent}, we only need to consider those terms involving $\bm{X}$ in the telescoping summation decomposition for $ T_{\beta^{[\zeta]}}^{\underline{[\bm{X}+t\bm{Y},(\bm{X}+t\bm{Y})_N]}_\ell}([\bm{Y}]_1^{\zeta}) -  T_{\beta^{[\zeta]}}^{\underline{[\bm{X},\bm{X}_N]}_\ell}([\bm{Y}]_1^{\zeta})$. As derivative for GMOIs involving $\bm{X}_N$ will be zero matrices. 

Therefore, we have
\begin{eqnarray}\label{eq4:lma:diff of GMOI with para matrices having some org matrices}
\lefteqn{\lim\limits_{t \rightarrow 0} \frac{T_{\beta^{[\zeta]}}^{\underline{[\bm{X}+t\bm{Y},(\bm{X}+t\bm{Y})_N]}_\ell}([\bm{Y}]_1^{\zeta}) -  T_{\beta^{[\zeta]}}^{\underline{[\bm{X},\bm{X}_N]}_\ell}([\bm{Y}]_1^{\zeta})}{t}=}\nonumber \\
&& \lim\limits_{t \rightarrow 0} \frac{1}{t}
\sum\limits_{j=1}^{\left\vert\underline{[\bm{X},\bm{X}_N]}_\ell\right\vert_{\bm{X}}}
\Bigg[T_{\beta^{[\zeta]}}^{
[\underline{[\bm{X},\bm{X}_N]}_\ell]_1^{\vartheta_\ell(j) -1},\tilde{\bm{X}},
[\underline{[\tilde{\bm{X}},\tilde{\bm{X}}_N]}_\ell]_{\vartheta_\ell(j)+1}^{\zeta+1}
}([\bm{Y}]_1^{\vartheta_\ell(j)-1}, [\bm{Y}]_{\vartheta_\ell(j)}^{\zeta}) \nonumber \\
&& -  T_{\beta^{[\zeta]}}^{
[\underline{[\bm{X},\bm{X}_N]}_\ell]_1^{\vartheta_\ell(j) -1},\bm{X},
[\underline{[\tilde{\bm{X}},\tilde{\bm{X}}_N]}_\ell]_{\vartheta_\ell(j)+1}^{\zeta+1}
}([\bm{Y}]_1^{\vartheta_\ell(j)-1}, [\bm{Y}]_{\vartheta_\ell(j)}^{\zeta})\Bigg].
\end{eqnarray}

Again, by the perturbation formula given by Theorem~\ref{thm:GMOI Perturbation Formula}, we obtain
\begin{eqnarray}\label{eq5:lma:diff of GMOI with para matrices having some org matrices}
\lefteqn{\mbox{R.H.S. of Eq.~\eqref{eq4:lma:diff of GMOI with para matrices having some org matrices}}}\nonumber \\
&=&\lim\limits_{t \rightarrow 0} \frac{1}{t}  \sum\limits_{j=1}^{\left\vert\underline{[\bm{X},\bm{X}_N]}_\ell\right\vert_{\bm{X}}}\Bigg[
T_{\beta^{[\zeta+1]}}^{
[\underline{[\bm{X},\bm{X}_N]}_\ell]_1^{\vartheta_\ell(j) -1},\tilde{\bm{X}},\bm{X},
[\underline{[\tilde{\bm{X}},\tilde{\bm{X}}_N]}_\ell]_{\vartheta_\ell(j)+1}^{\zeta+1}
}([\bm{Y}]_1^{\vartheta_\ell(j)-1}, t\bm{Y}, [\bm{Y}]_{\vartheta_\ell(j)}^{\zeta}) \nonumber \\
&& -  T_{\beta^{[\zeta+1]}}^{
[\underline{[\bm{X},\bm{X}_N]}_\ell]_1^{\vartheta_\ell(j) -1},\tilde{\bm{X}}_N,\bm{X}_N,
[\underline{[\tilde{\bm{X}},\tilde{\bm{X}}_N]}_\ell]_{\vartheta_\ell(j)+1}^{\zeta+1}
}([\bm{Y}]_1^{\vartheta_\ell(j)-1}, t\bm{Y}, [\bm{Y}]_{\vartheta_\ell(j)}^{\zeta}) \nonumber \\
&& - \mathfrak{X}\left(
[\underline{[\bm{X},\bm{X}_N]}_\ell]_1^{\vartheta_\ell(j) -1},\tilde{\bm{X}},\bm{X},
[\underline{[\tilde{\bm{X}},\tilde{\bm{X}}_N]}_\ell]_{\vartheta_\ell(j)+1}^{\zeta+1}
\right)\Bigg].
\end{eqnarray}
This lemma is proved by setting $t=0$ in Eq.~\eqref{eq5:lma:diff of GMOI with para matrices having some org matrices}.
$\hfill\Box$

The exact computation of $ \left.\mathfrak{X}^{(1)}\left(
[\underline{[\bm{X},\bm{X}_N]}_\ell]_1^{\vartheta_\ell(j) -1},\tilde{\bm{X}},\bm{X},
[\underline{[\tilde{\bm{X}},\tilde{\bm{X}}_N]}_\ell]_{\vartheta_\ell(j)+1}^{\zeta+1}
\right)\right\vert_{t=0}$ is tedious and there is no simple formula to express this exactly due to various combination of entries in the array $[\bm{X},\bm{X}_N]_\ell$. But, we will provide an example to show how to calculate this. 

We are ready to present the main Theorem~\ref{thm:n-th der by GMOI} in this section about the representing $n$-th derivative of a matrix-valued function in terms of GMOIs.
\begin{theorem}\label{thm:n-th der by GMOI}
Given any matrices $\bm{X}$ and $\bm{Y}$, we have
\begin{eqnarray}\label{eq1:thm:n-th der by GMOI}
\left.\frac{d^n f (\bm{X}+t \bm{Y})}{d t^n}\right\vert_{t=0}&=&n! T^{[\bm{X}]_1^{n+1}}_{f^{[n]}}([\bm{Y}]_1^n)-
(n-1)!\sum\limits_{j=1}^n T_{f^{[n]}}^{[\bm{X}]_1^{j-1},\bm{X}_N,\bm{X}_N,[\bm{X}]_{j+1}^{n}}([\bm{Y}]_1^n) \nonumber \\
&&-\sum\limits_{i=1}^{n-1}(n-i)!\left.\sum\limits_{j=1}^{n-i+1} \mathfrak{X}^{(i)}([\bm{X}]_1^{j-1},\tilde{\bm{X}}, \bm{X},[\tilde{\bm{X}}]_{j+1}^{n-i+1})\right\vert_{t=0} \nonumber \\
&&- \sum\limits_{\ell_{n-1}=1}^{\left(\sum\limits_{k=1}^{\lceil\frac{n-1}{2}\rceil}\frac{(n-k)!}{k!(n-2k)!}\right)- 1(\mbox{~if $n$ is even})} a_{\ell_{n-1}}\left.\frac{d T_{f^{[n-1]}}^{\underline{[\tilde{\bm{X}},\tilde{\bm{X}}_N]}_{\ell_{n-1}}}([\bm{Y}]_1^{n-1})}{dt}\right\vert_{t=0} \nonumber \\
&&+ \left[ \sum\limits_{\rho=4}^n\sum\limits_{i_\rho=1}^{\gamma(\rho)} b_{i_\rho}\left.\mathfrak{X}^{(n+2-\rho)}([\tilde{\bm{X}},\bm{X},\bm{X}_N]_{i_\rho})\right\vert_{t=0}\right]\bm{1}(\mbox{if $n \geq 4$}),
\end{eqnarray}
where $n \geq 2$, $a_{\ell_{n-1}}$ are coefficients for those GMOIs with parameters matrices containing both $\bm{X}$ and $\bm{X}_N$ from $\frac{d^{n-1} f (\bm{X}+t \bm{Y})}{d t^{n-1}}$, $b_{i_\rho}$ are coefficients for those $\mathfrak{X}^{(k)}([\tilde{\bm{X}},\bm{X},\bm{X}_N]_{i_\rho})$ from  $\frac{d^{n-1} f (\bm{X}+t \bm{Y})}{d t^{n-1}}$ with $k=1,2,\ldots,n-3$, and $[\tilde{\bm{X}},\bm{X},\bm{X}_N]_{i_\rho}$ is the array indexed by $i_\rho$ composed by $\tilde{\bm{X}}, \bm{X}$ and $\bm{X}_N$, respectively with length $\rho$. The term $\gamma(\rho)$ is determined by
\begin{eqnarray}\label{eq2:thm:n-th der by GMOI}
\gamma(\rho)&=&  \left[\frac{(\rho-2)!}{1!(\rho-3)!},  \frac{(\rho-3)!}{2!(\rho-5)!}, \ldots, 1 \right] \cdot [\rho-3, \rho-5, \ldots,0],
\end{eqnarray}
where $\cdot$ is the inner product operation. Note that $\bm{1}(\mbox{if $n \geq 4$})$ is the indicator condition function.
\end{theorem}
\textbf{Proof:}
Becuase we have
\begin{eqnarray}\label{eq2:thm:n-th der by GMOI}
\left.\frac{d f (\bm{X}+t \bm{Y})}{d t}\right\vert_{t=0}&=&T^{[\bm{X}]_1^{2}}_{f^{[1]}}([\bm{Y}]), 
\end{eqnarray}
by applying Theorem~\ref{thm:GMOI Perturbation Formula} to Eq.~\ref{eq2:thm:n-th der by GMOI}, we obtain
\begin{eqnarray}\label{eq3:thm:n-th der by GMOI}
\left.\frac{d^2 f (\bm{X}+t \bm{Y})}{d^2 t}\right\vert_{t=0}&=&\lim\limits_{t\rightarrow 0}\frac{ T^{[\bm{X}+t\bm{Y}]_1^{2}}_{f^{[1]}}([\bm{Y}])  - T^{[\bm{X}]_1^{2}}_{f^{[1]}}([\bm{Y}])}{t}\nonumber \\
&=&2 T^{[\bm{X}]_1^{3}}_{f^{[2]}}([\bm{Y}]_1^2) - T^{\bm{X},\bm{X}_N,\bm{X}_N}_{f^{[2]}}([\bm{Y}]_1^2) - T^{\bm{X}_N,\bm{X}_N,\bm{X}}_{f^{[2]}}([\bm{Y}]_1^2) \nonumber \\
&&-\left.\mathfrak{X}^{(1)}(\tilde{\bm{X}}, \bm{X}, \tilde{\bm{X}})\right\vert_{t=0} -  \left.\mathfrak{X}^{(1)}(\bm{X}, \tilde{\bm{X}}, \bm{X})\right\vert_{t=0}.
\end{eqnarray}

By taking the derivative of $f(\bm{X}+t\bm{Y})$ repeatdely and applying Theorem~\ref{thm:GMOI Perturbation Formula} to every terms except those GMOIs with parameter matrices conposed by $\bm{X}$ and $\bm{X}_N$ jointly, we have
\begin{eqnarray}\label{eq4:thm:n-th der by GMOI}
\lefteqn{\left.\frac{d^n f (\bm{X}+t \bm{Y})}{d t^n}\right\vert_{t=0}, \mbox{without $T^{\underline{[\bm{X}, \bm{X}_N]}}$ derivatives}}\nonumber \\
&=& n! T^{[\bm{X}]_1^{n+1}}_{f^{[n]}}([\bm{Y}]_1^n)-
(n-1)!\sum\limits_{j=1}^n T_{f^{[n]}}^{[\bm{X}]_1^{j-1},\bm{X}_N,\bm{X}_N,[\bm{X}]_{j+1}^{n}}([\bm{Y}]_1^n) \nonumber \\
&&-\sum\limits_{i=1}^{n-1}(n-i)!\left.\sum\limits_{j=1}^{n-i+1} \mathfrak{X}^{(i)}([\bm{X}]_1^{j-1},\tilde{\bm{X}}, \bm{X},[\tilde{\bm{X}}]_{j+1}^{n-i+1})\right\vert_{t=0}.
\end{eqnarray}

For those derivatives of GMOIs involving  $T^{[\bm{X}, \bm{X}_N]}$, there are $\left(\sum\limits_{k=1}^{\lceil\frac{n-1}{2}\rceil}\frac{(n-k)!}{k!(n-2k)!}\right)$ terms of GMOIs with parameter matrices having $\bm{X}_N$. However, from Lemma~\ref{lma:diff of GMOI with all para matrices are nilpotent}, we only need to consider GMOIs with parameter matrices composed by $\bm{X}$ and $\bm{X}_N$ jointly. Therefore, we need to reduce one in the total number of $\ell_{n-1}$. By Lemma~\ref{lma:diff of GMOI with para matrices having some org matrices}, we can evaluate each term of $\left.\frac{d T_{f^{[n-1]}}^{\underline{[\tilde{\bm{X}},\tilde{\bm{X}}_N]}_{\ell_{n-1}}}([\bm{Y}]_1^{n-1})}{dt}\right\vert_{t=0}$. Because there are $\gamma(\rho)$ terms for each $i_\rho$ to form $\mathfrak{X}^{(n+2-\rho)}([\tilde{\bm{X}},\bm{X},\bm{X}_N]_{i_\rho})$ when we apply perturbation formula from Theorem~\ref{thm:GMOI Perturbation Formula} in evaluating $\frac{d T_{f^{[n-1]}}^{\underline{[\tilde{\bm{X}},\tilde{\bm{X}}_N]}_{\ell_{n-1}}}([\bm{Y}]_1^{n-1})}{dt}$. 

Let us explain why the value $\gamma(\rho)$ can be expressed by
\begin{eqnarray}\label{eq4-1:thm:n-th der by GMOI}
\gamma(\rho)&=&  \left[\frac{(\rho-2)!}{1!(\rho-3)!},  \frac{(\rho-3)!}{2!(\rho-5)!}, \ldots, 1 \right] \cdot [\rho-3, \rho-5, \ldots,0].
\end{eqnarray}

We obtain
\[
\frac{(\rho - k - 1)!}{k! (\rho - 2k - 1)!}
\]
GMOIs of parameter matrices with the length \( \rho - 1 \) from the derivative
\[
\frac{d^{n-2}}{dt^{n-2}} f(\bm{X} + t \bm{Y}),
\]
since we begin considering the operator derivative structure of
\(\mathfrak{X}([\tilde{\bm{X}}, \bm{X}, \bm{X}_N]_{i_\rho})\)
from order 2. Moreover, each such GMOI includes \( \rho - 2k - 1 \) instances of \( \tilde{\bm{X}} \). As a result, we have
\[
\frac{(\rho - k - 1)!}{k! (\rho - 2k - 1)!} \times (\rho - 2k - 1)
\]
terms of the form
\(\mathfrak{X}([\tilde{\bm{X}}, \bm{X}, \bm{X}_N])\),
where the total array length is \( \rho \), and there are exactly \( k \) pairs of \( \tilde{\bm{X}}_N \). The total count
\(\gamma(\rho)\) is then obtained by summing over all admissible values of \( k \).

Therefore, we have
\begin{eqnarray}\label{eq5:thm:n-th der by GMOI}
\lefteqn{\left.\frac{d^n f (\bm{X}+t \bm{Y})}{d t^n}\right\vert_{t=0}, \mbox{with $T^{\underline{[\bm{X}, \bm{X}_N]}}$ derivatives}}\nonumber \\
&=& - \sum\limits_{\ell_{n-1}=1}^{\left(\sum\limits_{k=1}^{\lceil\frac{n-1}{2}\rceil}\frac{(n-k)!}{k!(n-2k)!}\right)- 1(\mbox{~if $n$ is even})} a_{\ell_{n-1}}\left.\frac{d T_{f^{[n-1]}}^{\underline{[\tilde{\bm{X}},\tilde{\bm{X}}_N]}_{\ell_{n-1}}}([\bm{Y}]_1^{n-1})}{dt}\right\vert_{t=0} \nonumber \\
&&+  \sum\limits_{\rho=4}^n\sum\limits_{i_\rho=1}^{\gamma(\rho)} b_{i_\rho}\left.\mathfrak{X}^{(n+2-\rho)}([\tilde{\bm{X}},\bm{X},\bm{X}_N]_{i_\rho})\right\vert_{t=0}.
\end{eqnarray}
This theorem is proved by combing Eq.~\eqref{eq4:thm:n-th der by GMOI} and Eq.~\eqref{eq5:thm:n-th der by GMOI}
$\hfill\Box$

Below, we will present examples to illustrate $\left.\frac{d^3 f (\bm{X}+t \bm{Y})}{d^3 t}\right\vert_{t=0}$ and $\left.\frac{d^4 f (\bm{X}+t \bm{Y})}{d^4 t}\right\vert_{t=0}$. 
\begin{example}\label{exp: der order 3 and 4}
Since we have
\begin{eqnarray}\label{eq1:exp: der order 3 and 4}
\left.\frac{d^2 f (\bm{X}+t \bm{Y})}{d^2 t}\right\vert_{t=0}&=&\lim\limits_{t\rightarrow 0}\frac{ T^{[\bm{X}+t\bm{Y}]_1^{2}}_{f^{[1]}}([\bm{Y}])  - T^{[\bm{X}]_1^{2}}_{f^{[1]}}([\bm{Y}])}{t}\nonumber \\
&=&2 T^{[\bm{X}]_1^{3}}_{f^{[2]}}([\bm{Y}]_1^2) - T^{\bm{X},\bm{X}_N,\bm{X}_N}_{f^{[2]}}([\bm{Y}]_1^2) - T^{\bm{X}_N,\bm{X}_N,\bm{X}}_{f^{[2]}}([\bm{Y}]_1^2) \nonumber \\
&&-\left.\mathfrak{X}^{(1)}(\tilde{\bm{X}}, \bm{X}, \tilde{\bm{X}})\right\vert_{t=0} -  \left.\mathfrak{X}^{(1)}(\bm{X}, \tilde{\bm{X}}, \bm{X})\right\vert_{t=0},
\end{eqnarray}
by the telescoping summation and the perturbation formula given by Theorem~\ref{thm:GMOI Perturbation Formula}, we have 
\begin{eqnarray}\label{eq2:exp: der order 3 and 4}
\left.\frac{d 2 T^{[\bm{X}]_1^{3}}_{f^{[2]}}([\bm{Y}]_1^2)}{d t}\right\vert_{t=0}&=& 3! T^{[\bm{X}]_1^{4}}_{f^{[3]}}([\bm{Y}]_1^3) - 2 T^{\bm{X}_N,\bm{X}_N,\bm{X},\bm{X}}_{f^{[3]}}([\bm{Y}]_1^3) \nonumber \\
&& - 2 T^{\bm{X}, \bm{X}_N,\bm{X}_N,\bm{X}}_{f^{[3]}}([\bm{Y}]_1^3) - 2 T^{\bm{X},\bm{X},\bm{X}_N,\bm{X}_N}_{f^{[3]}}([\bm{Y}]_1^3)\nonumber \\
&& - 2 \mathfrak{X}^{(1)}(\tilde{\bm{X}},\bm{X},\tilde{\bm{X}},\tilde{\bm{X}}) - 2 \mathfrak{X}^{(1)}(\bm{X}, \tilde{\bm{X}},\bm{X},\tilde{\bm{X}}) \nonumber \\
&& - 2 \mathfrak{X}^{(1)}(\bm{X},\bm{X},\tilde{\bm{X}},\bm{X}).
\end{eqnarray}
By applying Lemma~\ref{lma:diff of GMOI with para matrices having some org matrices} to $T^{\bm{X},\bm{X}_N,\bm{X}_N}_{f^{[2]}}([\bm{Y}]_1^2)$ and $T^{\bm{X}_N,\bm{X}_N,\bm{X}}_{f^{[2]}}([\bm{Y}]_1^2)$, we have 
\begin{eqnarray}\label{eq3:exp: der order 3 and 4}
\left.\frac{d T^{\bm{X},\bm{X}_N,\bm{X}_N}_{f^{[2]}}([\bm{Y}]_1^2)}{d t}\right\vert_{t=0}&=& T^{\bm{X},\bm{X},\bm{X}_N,\bm{X}_N}_{f^{[3]}}([\bm{Y}]_1^3) - T^{\bm{X}_N,\bm{X}_N,\bm{X}_N,\bm{X}_N}_{f^{[3]}}([\bm{Y}]_1^3) \nonumber \\
&&-\mathfrak{X}^{(1)}(\tilde{\bm{X}},\bm{X},\tilde{\bm{X}}_N,\tilde{\bm{X}}_N),
\end{eqnarray}
and
\begin{eqnarray}\label{eq4:exp: der order 3 and 4}
\left.\frac{d T^{\bm{X}_N,\bm{X}_N,\bm{X}}_{f^{[2]}}([\bm{Y}]_1^2)}{d t}\right\vert_{t=0}&=& T^{\bm{X}_N,\bm{X}_N,\bm{X},\bm{X}}_{f^{[3]}}([\bm{Y}]_1^3) - T^{\bm{X}_N,\bm{X}_N,\bm{X}_N,\bm{X}_N}_{f^{[3]}}([\bm{Y}]_1^3) \nonumber \\
&&-\mathfrak{X}^{(1)}(\bm{X}_N,\bm{X}_N,\tilde{\bm{X}},\bm{X}).
\end{eqnarray}
By combining Eq.~\eqref{eq2:exp: der order 3 and 4} to Eq.~\eqref{eq4:exp: der order 3 and 4}, we have
\begin{eqnarray}\label{eq5:exp: der order 3 and 4}
\left.\frac{d^3 f (\bm{X}+t \bm{Y})}{d^3 t}\right\vert_{t=0}&=& 3! T^{[\bm{X}]_1^{4}}_{f^{[3]}}([\bm{Y}]_1^3) - 3 T^{\bm{X}_N,\bm{X}_N,\bm{X},\bm{X}}_{f^{[3]}}([\bm{Y}]_1^3) \nonumber \\
&& - 2 T^{\bm{X}, \bm{X}_N,\bm{X}_N,\bm{X}}_{f^{[3]}}([\bm{Y}]_1^3) - 3 T^{\bm{X},\bm{X},\bm{X}_N,\bm{X}_N}_{f^{[3]}}([\bm{Y}]_1^3)\nonumber \\
&& - 2 \left.\mathfrak{X}^{(1)}(\tilde{\bm{X}},\bm{X},\tilde{\bm{X}},\tilde{\bm{X}})\right\vert_{t=0} - 2\left. \mathfrak{X}^{(1)}(\bm{X}, \tilde{\bm{X}},\bm{X},\tilde{\bm{X}})\right\vert_{t=0} \nonumber \\
&& - 2 \left.\mathfrak{X}^{(1)}(\bm{X},\bm{X},\tilde{\bm{X}},\bm{X})\right\vert_{t=0} + 2T^{\bm{X}_N,\bm{X}_N,\bm{X}_N,\bm{X}_N}_{f^{[3]}}([\bm{Y}]_1^3) \nonumber \\
&&-\left.\mathfrak{X}^{(2)}(\tilde{\bm{X}}, \bm{X}, \tilde{\bm{X}})\right\vert_{t=0} -  \left.\mathfrak{X}^{(2)}(\bm{X}, \tilde{\bm{X}}, \bm{X})\right\vert_{t=0}\nonumber \\
&&+\left.\mathfrak{X}^{(1)}(\tilde{\bm{X}},\bm{X},\tilde{\bm{X}}_N,\tilde{\bm{X}}_N)\right\vert_{t=0} +\left.\mathfrak{X}^{(1)}(\bm{X}_N,\bm{X}_N,\tilde{\bm{X}},\bm{X})\right\vert_{t=0} \nonumber \\
\end{eqnarray}

For the evaluation of $\left.\frac{d^4 f (\bm{X}+t \bm{Y})}{d^4 t}\right\vert_{t=0}$, we first derive
\begin{eqnarray}\label{eq5:exp: der order 3 and 4}
 3! \frac{d T^{[\bm{X}]_1^{4}}_{f^{[3]}}([\bm{Y}]_1^3)}{d t}&=&  4! T^{[\bm{X}]_1^{5}}_{f^{[4]}}([\bm{Y}]_1^4) - 3! T^{\bm{X}_N,\bm{X}_N,\bm{X},\bm{X},\bm{X}}_{f^{[4]}}([\bm{Y}]_1^4) \nonumber \\
&& - 3! T^{\bm{X}, \bm{X}_N,\bm{X}_N,\bm{X},\bm{X}}_{f^{[4]}}([\bm{Y}]_1^4) - 3! T^{\bm{X},\bm{X},\bm{X}_N,\bm{X}_N, \bm{X}}_{f^{[4]}}([\bm{Y}]_1^4) - 3! T^{\bm{X},\bm{X},\bm{X},\bm{X}_N,\bm{X}_N}_{f^{[4]}}([\bm{Y}]_1^4)  \nonumber \\
&& - 3!  \mathfrak{X}^{(1)}(\tilde{\bm{X}},\bm{X},\tilde{\bm{X}},\tilde{\bm{X}},\tilde{\bm{X}}) - 3! \mathfrak{X}^{(1)}(\bm{X}, \tilde{\bm{X}},\bm{X},\tilde{\bm{X}},\tilde{\bm{X}}) \nonumber \\
&& - 3!  \mathfrak{X}^{(1)}(\bm{X},\bm{X},\tilde{\bm{X}},\bm{X},\tilde{\bm{X}})- 3!  \mathfrak{X}^{(1)}(\bm{X},\bm{X},\bm{X},\tilde{\bm{X}},\bm{X}).
\end{eqnarray}
by the telescoping summation and the perturbation formula given by Theorem~\ref{thm:GMOI Perturbation Formula}.

We also have
\begin{eqnarray}\label{eq6-1:exp: der order 3 and 4}
\left.\frac{d T^{\bm{X}_N,\bm{X}_N,\bm{X},\bm{X}}_{f^{[3]}}([\bm{Y}]_1^3)}{d t}\right\vert_{t=0}&=& T^{\bm{X}_N,\bm{X}_N,\bm{X},\bm{X},\bm{X}}_{f^{[4]}}([\bm{Y}]_1^4) - T^{\bm{X}_N,\bm{X}_N,\bm{X}_N,\bm{X}_N,\bm{X}}_{f^{[4]}}([\bm{Y}]_1^4) \nonumber \\
&&-\left.\mathfrak{X}^{(1)}(\bm{X}_N,\bm{X}_N,\tilde{\bm{X}},\bm{X},\tilde{\bm{X}})\right\vert_{t=0}\nonumber \\
&& +  T^{\bm{X}_N,\bm{X}_N,\bm{X},\bm{X},\bm{X}}_{f^{[4]}}([\bm{Y}]_1^4) - T^{\bm{X}_N,\bm{X}_N,\bm{X},\bm{X}_N,\bm{X}_N}_{f^{[4]}}([\bm{Y}]_1^4) \nonumber \\
&&-\left.\mathfrak{X}^{(1)}(\bm{X}_N,\bm{X}_N,\bm{X},\tilde{\bm{X}},\bm{X})\right\vert_{t=0},
\end{eqnarray}
\begin{eqnarray}\label{eq6-2:exp: der order 3 and 4}
\left.\frac{d T^{\bm{X},\bm{X}_N,\bm{X}_N,\bm{X}}_{f^{[3]}}([\bm{Y}]_1^3)}{d t}\right\vert_{t=0}&=& T^{\bm{X},\bm{X},\bm{X}_N,\bm{X}_N,\bm{X}}_{f^{[4]}}([\bm{Y}]_1^4) - T^{\bm{X}_N,\bm{X}_N,\bm{X}_N,\bm{X}_N,\bm{X}}_{f^{[4]}}([\bm{Y}]_1^4) \nonumber \\
&&-\left.\mathfrak{X}^{(1)}(\tilde{\bm{X}},\bm{X},\tilde{\bm{X}}_N,\tilde{\bm{X}}_N,\tilde{\bm{X}})\right\vert_{t=0} \nonumber \\
&& +  T^{\bm{X},\bm{X}_N,\bm{X}_N,\bm{X},\bm{X}}_{f^{[4]}}([\bm{Y}]_1^4) - T^{\bm{X},\bm{X}_N,\bm{X}_N,\bm{X}_N,\bm{X}_N}_{f^{[4]}}([\bm{Y}]_1^4) \nonumber \\
&&-\left.\mathfrak{X}^{(1)}(\bm{X},\bm{X}_N,\bm{X}_N,\tilde{\bm{X}},\bm{X})\right\vert_{t=0},
\end{eqnarray}
and
\begin{eqnarray}\label{eq6-3:exp: der order 3 and 4}
\left.\frac{d T^{\bm{X},\bm{X},\bm{X}_N,\bm{X}_N}_{f^{[3]}}([\bm{Y}]_1^3)}{d t}\right\vert_{t=0}&=& T^{\bm{X},\bm{X},\bm{X},\bm{X}_N,\bm{X}_N}_{f^{[4]}}([\bm{Y}]_1^4) - T^{\bm{X}_N,\bm{X}_N,\bm{X},\bm{X}_N,\bm{X}_N}_{f^{[4]}}([\bm{Y}]_1^4) \nonumber \\
&&-\left.\mathfrak{X}^{(1)}(\tilde{\bm{X}},\bm{X},\tilde{\bm{X}},\tilde{\bm{X}}_N,\tilde{\bm{X}}_N)\right\vert_{t=0} \nonumber \\
&& +  T^{\bm{X},\bm{X},\bm{X},\bm{X}_N,\bm{X}_N}_{f^{[4]}}([\bm{Y}]_1^4) - T^{\bm{X},\bm{X}_N,\bm{X}_N,\bm{X}_N,\bm{X}_N}_{f^{[4]}}([\bm{Y}]_1^4) \nonumber \\
&&-\left.\mathfrak{X}^{(1)}(\bm{X},\tilde{\bm{X}},\bm{X},\tilde{\bm{X}}_N,\tilde{\bm{X}}_N)\right\vert_{}.
\end{eqnarray}

By combining Eq.~\eqref{eq5:exp: der order 3 and 4} to Eq.~\eqref{eq6-3:exp: der order 3 and 4}, we have
\begin{eqnarray}\label{eq7:exp: der order 3 and 4}
\lefteqn{\left.\frac{d^4 f (\bm{X}+t \bm{Y})}{d^4 t}\right\vert_{t=0}}\nonumber\\
&=& 4! T^{[\bm{X}]_1^{5}}_{f^{[4]}}([\bm{Y}]_1^4) - 3! T^{\bm{X}_N,\bm{X}_N,\bm{X},\bm{X},\bm{X}}_{f^{[4]}}([\bm{Y}]_1^4) \nonumber \\
&& - 3! T^{\bm{X}, \bm{X}_N,\bm{X}_N,\bm{X},\bm{X}}_{f^{[4]}}([\bm{Y}]_1^4) - 3! T^{\bm{X},\bm{X},\bm{X}_N,\bm{X}_N, \bm{X}}_{f^{[4]}}([\bm{Y}]_1^4) - 3! T^{\bm{X},\bm{X},\bm{X},\bm{X}_N,\bm{X}_N}_{f^{[4]}}([\bm{Y}]_1^4)  \nonumber \\
&& - 3!  \left.\mathfrak{X}^{(1)}(\tilde{\bm{X}},\bm{X},\tilde{\bm{X}},\tilde{\bm{X}},\tilde{\bm{X}})\right\vert_{t=0} - 3! \left.\mathfrak{X}^{(1)}(\bm{X}, \tilde{\bm{X}},\bm{X},\tilde{\bm{X}},\tilde{\bm{X}})\right\vert_{t=0} \nonumber \\
&& - 3!  \left.\mathfrak{X}^{(1)}(\bm{X},\bm{X},\tilde{\bm{X}},\bm{X},\tilde{\bm{X}})\right\vert_{t=0} - 3!  \left.\mathfrak{X}^{(1)}(\bm{X},\bm{X},\bm{X},\tilde{\bm{X}},\bm{X})\right\vert_{t=0} \nonumber \\
&&-3\Bigg(T^{\bm{X}_N,\bm{X}_N,\bm{X},\bm{X},\bm{X}}_{f^{[4]}}([\bm{Y}]_1^4) - T^{\bm{X}_N,\bm{X}_N,\bm{X}_N,\bm{X}_N,\bm{X}}_{f^{[4]}}([\bm{Y}]_1^4) \nonumber \\
&&-\left.\mathfrak{X}^{(1)}(\bm{X}_N,\bm{X}_N,\tilde{\bm{X}},\bm{X},\tilde{\bm{X}})\right\vert_{t=0} \nonumber \\
&& +  T^{\bm{X}_N,\bm{X}_N,\bm{X},\bm{X},\bm{X}}_{f^{[4]}}([\bm{Y}]_1^4) - T^{\bm{X}_N,\bm{X}_N,\bm{X},\bm{X}_N,\bm{X}_N}_{f^{[4]}}([\bm{Y}]_1^4) \nonumber \\
&&-\left.\mathfrak{X}^{(1)}(\bm{X}_N,\bm{X}_N,\bm{X},\tilde{\bm{X}},\bm{X})\right\vert_{t=0} \Bigg)\nonumber \\
&&-2\Bigg(T^{\bm{X},\bm{X},\bm{X}_N,\bm{X}_N,\bm{X}}_{f^{[4]}}([\bm{Y}]_1^4) - T^{\bm{X}_N,\bm{X}_N,\bm{X}_N,\bm{X}_N,\bm{X}}_{f^{[4]}}([\bm{Y}]_1^4) \nonumber \\
&&-\left.\mathfrak{X}^{(1)}(\tilde{\bm{X}},\bm{X},\tilde{\bm{X}}_N,\tilde{\bm{X}}_N,\tilde{\bm{X}})\right\vert_{t=0} \nonumber \\
&& +  T^{\bm{X},\bm{X}_N,\bm{X}_N,\bm{X},\bm{X}}_{f^{[4]}}([\bm{Y}]_1^4) - T^{\bm{X},\bm{X}_N,\bm{X}_N,\bm{X}_N,\bm{X}_N}_{f^{[4]}}([\bm{Y}]_1^4) \nonumber \\
&&-\left.\mathfrak{X}^{(1)}(\bm{X},\bm{X}_N,\bm{X}_N,\tilde{\bm{X}},\bm{X})\right\vert_{t=0} \Bigg)\nonumber \\
&&-3\Bigg(T^{\bm{X},\bm{X},\bm{X},\bm{X}_N,\bm{X}_N}_{f^{[4]}}([\bm{Y}]_1^4) - T^{\bm{X}_N,\bm{X}_N,\bm{X},\bm{X}_N,\bm{X}_N}_{f^{[4]}}([\bm{Y}]_1^4) \nonumber \\
&&- \left.\mathfrak{X}^{(1)}(\tilde{\bm{X}},\bm{X},\tilde{\bm{X}}_N,\tilde{\bm{X}}_N,\tilde{\bm{X}})\right\vert_{t=0} \nonumber \\
&& +  T^{\bm{X},\bm{X},\bm{X},\bm{X}_N,\bm{X}_N}_{f^{[4]}}([\bm{Y}]_1^4) - T^{\bm{X},\bm{X}_N,\bm{X}_N,\bm{X}_N,\bm{X}_N}_{f^{[4]}}([\bm{Y}]_1^4) \nonumber \\
&&- \left.\mathfrak{X}^{(1)}(\bm{X},\tilde{\bm{X}},\bm{X},\tilde{\bm{X}}_N,\tilde{\bm{X}}_N)\right\vert_{t=0} \Bigg)\nonumber \\
&&  -2 \left.\mathfrak{X}^{(2)}(\tilde{\bm{X}},\bm{X},\tilde{\bm{X}},\tilde{\bm{X}})\right\vert_{t=0} - 2\left. \mathfrak{X}^{(2)}(\bm{X}, \tilde{\bm{X}},\bm{X},\tilde{\bm{X}})\right\vert_{t=0} \nonumber \\
&& - 2 \left.\mathfrak{X}^{(2)}(\bm{X},\bm{X},\tilde{\bm{X}},\bm{X})\right\vert_{t=0} \nonumber \\
&&-\left.\mathfrak{X}^{(3)}(\tilde{\bm{X}}, \bm{X}, \tilde{\bm{X}})\right\vert_{t=0} -  \left.\mathfrak{X}^{(3)}(\bm{X}, \tilde{\bm{X}}, \bm{X})\right\vert_{t=0}\nonumber \\
&&+\left.\mathfrak{X}^{(2)}(\tilde{\bm{X}},\bm{X},\tilde{\bm{X}}_N,\tilde{\bm{X}}_N)\right\vert_{t=0} +\left.\mathfrak{X}^{(2)}(\bm{X}_N,\bm{X}_N,\tilde{\bm{X}},\bm{X})\right\vert_{t=0}.
\end{eqnarray}
After arranging terms, we have
\begin{eqnarray}\label{eq7:exp: der order 3 and 4}
\lefteqn{\left.\frac{d^4 f (\bm{X}+t \bm{Y})}{d^4 t}\right\vert_{t=0}}\nonumber\\
&=& 4! T^{[\bm{X}]_1^{5}}_{f^{[4]}}([\bm{Y}]_1^4) - 12 T^{\bm{X}_N,\bm{X}_N,\bm{X},\bm{X},\bm{X}}_{f^{[4]}}([\bm{Y}]_1^4) \nonumber \\
&& - 8 T^{\bm{X}, \bm{X}_N,\bm{X}_N,\bm{X},\bm{X}}_{f^{[4]}}([\bm{Y}]_1^4) - 8 T^{\bm{X},\bm{X},\bm{X}_N,\bm{X}_N, \bm{X}}_{f^{[4]}}([\bm{Y}]_1^4) - 12 T^{\bm{X},\bm{X},\bm{X},\bm{X}_N,\bm{X}_N}_{f^{[4]}}([\bm{Y}]_1^4)  \nonumber \\
&& - 3! \left.\mathfrak{X}^{(1)}(\tilde{\bm{X}},\bm{X},\tilde{\bm{X}},\tilde{\bm{X}},\tilde{\bm{X}})\right\vert_{t=0}  - 3! \left.\mathfrak{X}^{(1)}(\bm{X}, \tilde{\bm{X}},\bm{X},\tilde{\bm{X}},\tilde{\bm{X}})\right\vert_{t=0} \nonumber \\
&& - 3!  \left.\mathfrak{X}^{(1)}(\bm{X},\bm{X},\tilde{\bm{X}},\bm{X},\tilde{\bm{X}})\right\vert_{t=0} - 3!  \left.\mathfrak{X}^{(1)}(\bm{X},\bm{X},\bm{X},\tilde{\bm{X}},\bm{X})\right\vert_{t=0}   \nonumber \\
&&  -2 \left.\mathfrak{X}^{(2)}(\tilde{\bm{X}},\bm{X},\tilde{\bm{X}},\tilde{\bm{X}})\right\vert_{t=0} - 2\left. \mathfrak{X}^{(2)}(\bm{X}, \tilde{\bm{X}},\bm{X},\tilde{\bm{X}})\right\vert_{t=0} \nonumber \\
&& - 2 \left.\mathfrak{X}^{(2)}(\bm{X},\bm{X},\tilde{\bm{X}},\bm{X})\right\vert_{t=0} -\left.\mathfrak{X}^{(3)}(\tilde{\bm{X}}, \bm{X}, \tilde{\bm{X}})\right\vert_{t=0} -  \left.\mathfrak{X}^{(3)}(\bm{X}, \tilde{\bm{X}}, \bm{X})\right\vert_{t=0}\nonumber \\
&&  +3 \left.\mathfrak{X}^{(1)}(\bm{X}_N,\bm{X}_N,\tilde{\bm{X}},\bm{X},\tilde{\bm{X}})\right\vert_{t=0}   +3 \left.\mathfrak{X}^{(1)}(\bm{X}_N,\bm{X}_N,\bm{X},\tilde{\bm{X}},\bm{X})\right\vert_{t=0} + 5 T^{\bm{X}_N,\bm{X}_N,\bm{X}_N,\bm{X}_N,\bm{X}}_{f^{[4]}}([\bm{Y}]_1^4) \nonumber \\
&&  +  2 \left.\mathfrak{X}^{(1)}(\bm{X},\bm{X}_N,\bm{X}_N,\tilde{\bm{X}},\bm{X})\right\vert_{t=0} + 6T^{\bm{X}_N,\bm{X}_N,\bm{X},\bm{X}_N,\bm{X}_N}_{f^{[4]}}([\bm{Y}]_1^4)+ 5\left.\mathfrak{X}^{(1)}(\tilde{\bm{X}},\bm{X},\tilde{\bm{X}}_N,\tilde{\bm{X}}_N,\tilde{\bm{X}})\right\vert_{t=0}  \nonumber \\
&&  +5 T^{\bm{X},\bm{X}_N,\bm{X}_N,\bm{X}_N,\bm{X}_N}_{f^{[4]}}([\bm{Y}]_1^4)  +3 \left.\mathfrak{X}^{(1)}(\bm{X},\tilde{\bm{X}},\bm{X},\tilde{\bm{X}}_N,\tilde{\bm{X}}_N)\right\vert_{t=0}  \nonumber \\
&&+\left.\mathfrak{X}^{(2)}(\tilde{\bm{X}},\bm{X},\tilde{\bm{X}}_N,\tilde{\bm{X}}_N)\right\vert_{t=0} +\left.\mathfrak{X}^{(2)}(\bm{X}_N,\bm{X}_N,\tilde{\bm{X}},\bm{X})\right\vert_{t=0}.
\end{eqnarray}
\end{example}

From Example~\ref{exp: der order 3 and 4}, we need to evaluate high-order derivatives of terms involving $\mathfrak{X}$. In following Example~\ref{exp: higher order derivative of mathfrak X}, we will try to evaluate any $n$-th order derivative of $\mathfrak{X}(\bm{X}, \tilde{\bm{X}}, \bm{X})$. For other evaluation of high-order derivatives of terms involving $\mathfrak{X}$, they can be calculated similarly although more complicated.

\begin{example}\label{exp: higher order derivative of mathfrak X}

By Eq.~\eqref{eq2:exp:GMOI Perturbation Formula}, the term $\mathfrak{X}(\bm{X}, \tilde{\bm{X}}, \bm{X})$ can be expressed by
\begin{eqnarray}\label{eq1:exp: higher order derivative of mathfrak X}
\bar{\mathfrak{X}}(\bm{C},\bm{D},\bm{X}_1)&=&  \mathfrak{X}^{\bm{C}_P,\bm{D}_P,\bm{X}_{1,P}} +  \mathfrak{X}^{\bm{C}_P,\bm{D}_P,\bm{X}_{1,N}} \nonumber \\
&&+ \mathfrak{X}^{\bm{C}_N,\bm{D}_N,\bm{X}_{1,P}} + \mathfrak{X}^{\bm{C}_N,\bm{D}_N,\bm{X}_{1,N}},
\end{eqnarray}
where $\bm{C}=\bm{X}, \bm{D}=\bm{X}+t\bm{Y}$ and $\bm{X}_1=\bm{X}$. Then, the $n$-th derivative of  $\mathfrak{X}(\bm{X}, \tilde{\bm{X}}, \bm{X})$ can be expressed by
\begin{eqnarray}\label{eq1:exp: higher order derivative of mathfrak X}
\frac{d^n}{dt^n}\bar{\mathfrak{X}}(\bm{C},\bm{D},\bm{X}_1)&=&\frac{d^n}{dt^n}\mathfrak{X}^{\bm{C}_P,\bm{D}_P,\bm{X}_{1,P}} + \frac{d^n}{dt^n}\mathfrak{X}^{\bm{C}_P,\bm{D}_P,\bm{X}_{1,N}} \nonumber \\
&&+\frac{d^n}{dt^n}\mathfrak{X}^{\bm{C}_N,\bm{D}_N,\bm{X}_{1,P}} + \frac{d^n}{dt^n}\mathfrak{X}^{\bm{C}_N,\bm{D}_N,\bm{X}_{1,N}}.
\end{eqnarray}
In the following, we will evaluate $\frac{d^n}{dt^n}\mathfrak{X}^{\bm{C}_P,\bm{D}_P,\bm{X}_{1,P}}$, $\frac{d^n}{dt^n}\mathfrak{X}^{\bm{C}_P,\bm{D}_P,\bm{X}_{1,N}}$, $\frac{d^n}{dt^n}\mathfrak{X}^{\bm{C}_N,\bm{D}_N,\bm{X}_{1,P}}$ and $\frac{d^n}{dt^n}\mathfrak{X}^{\bm{C}_N,\bm{D}_N,\bm{X}_{1,N}}$, respectively.

\textbf{The Evaluation of $\frac{d^n}{dt^n} \mathfrak{X}^{\bm{C}_P,\bm{D}_P,\bm{X}_{1,P}}$}

Because we have
\begin{eqnarray}\label{e3-1:exp: higher order derivative of mathfrak X}
\lefteqn{\mathfrak{X}^{\bm{C}_P,\bm{D}_P,\bm{X}_{1,P}}= }\nonumber \\
&&\sum\limits_{\substack{k_c, i_c,q_c,k_d,i_d,q_d\\ k_{1},i_{1}}}\frac{(\beta^{[2]}(\lambda_{k_c},\lambda_{k_d},\lambda_{k_1}))^{([0,0,0])}}{ ([0!,0!,0!])} \times \bm{P}_{k_c,i_c}(\bm{N}_{k_c,i_c}^{q_c}-\bm{N}_{k_d,i_d}^{q_d})\bm{P}_{k_d,i_d}\bm{Y}\bm{P}_{k_{1},i_{1}},
\end{eqnarray}
to derive the $n$-th derivative of the operator-valued expression
\[
\mathfrak{X}^{\bm{C}_P,\bm{D}_P,\bm{X}_{1,P}}
\]
with respect to $t$, from Eq.~\eqref{e3-1:exp: higher order derivative of mathfrak X}, we need to analyze how each term behaves under differentiation when the matrix $\bm{X}+t\bm{Y}$ appears implicitly in the spectral data ($\lambda_{k_d}$, $\bm{P}_{k_d,i_d}$, and $\bm{N}_{k_d,i_d}$). 

Differentiating $\mathfrak{X}^{\bm{C}_P,\bm{D}_P,\bm{X}_{1,P}}(t)$ with respect to $t$, the only $t$-dependence comes from:
\begin{itemize}
  \item $\lambda_{k_d}(t)$,
  \item $\bm{P}_{k_d,i_d}(t)$,
  \item $\bm{N}_{k_d,i_d}^{q_d}(t)$.
\end{itemize}

Then, the $n$-th derivative becomes:
\begin{eqnarray}
\lefteqn{\frac{d^n}{dt^n}\mathfrak{X}^{\bm{C}_P,\bm{D}_P,\bm{X}_{1,P}}(t)\bigg|_{t=0}}\nonumber \\
&=& \sum_{\substack{k_c, i_c,q_c,k_d,i_d,q_d\\ k_{1},i_{1}}}
\left.\frac{\partial^n}{\partial t^n}
\left[
\beta^{[2]}(\lambda_{k_c},\lambda_{k_d}(t),\lambda_{k_1})
\cdot 
\bm{P}_{k_c,i_c}(\bm{N}_{k_c,i_c}^{q_c} - \bm{N}_{k_d,i_d}^{q_d}(t))\bm{P}_{k_d,i_d}(t)\bm{Y}\bm{P}_{k_1,i_1}
\right]\right|_{t=0}
\end{eqnarray}

Since $\lambda_{k_c}, \lambda_{k_1}$ are independent of $t$, by applying multivariate Faà di Bruno formula for $\beta^{[2]}$, we have
\[
\frac{\partial^n}{\partial t^n}\beta^{[2]}(\lambda_{k_c},\lambda_{k_d}(t),\lambda_{k_1}) = \sum_{\pi \in \mathcal{P}(n)}
\frac{n!}{\prod_{j=1}^n b_j! \cdot j!^{b_j}}
\cdot \partial_{\lambda_{k_d}}^{|\pi|} \beta^{[2]}(\lambda_{k_c}, \lambda_{k_d}(t), \lambda_{k_1})
\cdot \prod_{j=1}^n \left( \lambda_{k_d}^{(j)}(t) \right)^{b_j},
\]
where $\mathcal{P}(n)$ denote the set of all integer partitions $\pi = (b_1, b_2, \dots, b_m)$ such that $\sum_{j=1}^n j \cdot b_j = n$.

Similarly,  by applying chain rule of derivative to projectors and nilpotents, we also have
\[
\frac{d^n}{dt^n} \left[
(\bm{N}_{k_c,i_c}^{q_c} - \bm{N}_{k_d,i_d}^{q_d}(t)) \cdot \bm{P}_{k_d,i_d}(t)
\right]
= \sum_{j=0}^n \binom{n}{j} 
\left( - \frac{d^j}{dt^j} \bm{N}_{k_d,i_d}^{q_d}(t) \right) 
\cdot \frac{d^{n-j}}{dt^{n-j}} \bm{P}_{k_d,i_d}(t).
\]
By defining 
\[
\bm{N}_{k_d,i_d}^{q_d,(j)} \define \left. \frac{d^j}{dt^j} \bm{N}_{k_d,i_d}^{q_d}(t) \right|_{t=0}, \quad
\bm{P}_{k_d,i_d}^{(n-j)} \define  \left.\frac{d^{n-j}}{dt^{n-j}} \bm{P}_{k_d,i_d}(t) \right|_{t=0},
\]
finally, we have
\begin{eqnarray}
\frac{d^n}{dt^n}\mathfrak{X}^{\bm{C}_P,\bm{D}_P,\bm{X}_{1,P}}(t)\bigg|_{t=0}
&=&
\sum_{\substack{k_c, i_c, q_c, k_d, i_d, q_d \\ k_1, i_1}}
\sum_{\substack{m + \ell = n}}
\sum_{\pi \in \mathcal{P}(m)}
\frac{1}{\ell!}
\cdot \frac{m!}{\prod_{j=1}^m b_j! j!^{b_j}} \nonumber \\
&& \cdot \left(\partial_{\lambda_{k_d}}^{|\pi|} \beta^{[2]}\right)(\lambda_{k_c}, \lambda_{k_d}, \lambda_{k_1})
\cdot \prod_{j=1}^m \left( \lambda_{k_d}^{(j)} \right)^{b_j} \nonumber \\
&& \cdot \bm{P}_{k_c,i_c}
\cdot \left[
\sum_{j=0}^{\ell} \binom{\ell}{j}
\left( - \bm{N}_{k_d,i_d}^{q_d,(j)} \right)
\cdot \bm{P}_{k_d,i_d}^{(\ell - j)}
\right]
\cdot \bm{Y}
\cdot \bm{P}_{k_1,i_1},
\end{eqnarray}
where $\lambda_{k_d}^{(j)} := \left. \frac{d^j}{dt^j} \lambda_{k_d}(t) \right|_{t=0}$. 

This formula computes how the generalized interpolation kernel operator $\mathfrak{X}$ changes under perturbations of the form $\bm{X} + t\bm{Y}$, taking into account the full Jordan structure and the role of divided differences in the operator function calculus.

\textbf{The evaluation of $\frac{d^n}{dt^n}  \mathfrak{X}^{\bm{C}_P,\bm{D}_P,\bm{X}_{1,N}}$}

Since we have
\begin{eqnarray}\label{e3-2:exp: higher order derivative of mathfrak X}
\lefteqn{ \mathfrak{X}^{\bm{C}_P,\bm{D}_P,\bm{X}_{1,N}} = }\nonumber \\
&&\sum\limits_{\substack{k_c, i_c,q_c,k_d,i_d,q_d,\\ k_{1},i_{1},q_1}}\frac{(\beta^{[2]}(\lambda_{k_c},\lambda_{k_d},\lambda_{k_1}))^{([0,0,q_1])}}{ ([0!,0!,q_1!])} \times \bm{P}_{k_c,i_c}(\bm{N}_{k_c,i_c}^{q_c}-\bm{N}_{k_d,i_d}^{q_d})\bm{P}_{k_d,i_d}\bm{Y}\bm{N}_{k_{1},i_{1}}^{q_1},
\end{eqnarray}
to derive the $n$-th derivative of the operator-valued expression

\[
\mathfrak{X}^{\bm{C}_P,\bm{D}_P,\bm{X}_{1,N}},
\]
with respect to $t$ based on Eq.~\eqref{e3-2:exp: higher order derivative of mathfrak X}, we must again analyze the impact of differentiating the terms where $\bm{X} + t\bm{Y}$ enters implicitly via the spectral data: $\lambda_{k_d}(t)$, $\bm{P}_{k_d,i_d}(t)$, and $\bm{N}_{k_d,i_d}(t)$. Finally, we have
\begin{eqnarray}
\frac{d^n}{dt^n} \mathfrak{X}^{\bm{C}_P,\bm{D}_P,\bm{X}_{1,N}}(t)\Big|_{t=0}
&= &\sum_{\substack{k_c, i_c,q_c,k_d,i_d,q_d,\\ k_{1},i_{1},q_1}}
\sum_{\substack{m + \ell = n}} 
\sum_{\pi \in \mathcal{P}(m)} 
\frac{1}{q_1! \cdot \ell!} 
\cdot \frac{m!}{\prod_{j=1}^m b_j! \cdot j!^{b_j}} \nonumber \\
&& \cdot \left(\partial_{\lambda_{k_d}}^{|\pi|} (\beta^{[2]})^{([0,0,q_1])}\right)(\lambda_{k_c}, \lambda_{k_d}, \lambda_{k_1})
\cdot \prod_{j=1}^m \left( \lambda_{k_d}^{(j)} \right)^{b_j} \nonumber \\
&& \cdot \bm{P}_{k_c,i_c}
\cdot \left[
\sum_{j=0}^{\ell} \binom{\ell}{j}
\left( - \bm{N}_{k_d,i_d}^{q_d,(j)} \right)
\cdot \bm{P}_{k_d,i_d}^{(\ell - j)}
\right]
\cdot \bm{Y}
\cdot \bm{N}_{k_1,i_1}^{q_1}.
\end{eqnarray}

\textbf{The Evaluation of $\frac{d^n}{dt^n} \mathfrak{X}^{\bm{C}_N,\bm{D}_N,\bm{X}_{1,P}}$}

Since we have
\begin{eqnarray}\label{e3-3:exp: higher order derivative of mathfrak X}
\lefteqn{ \mathfrak{X}^{\bm{C}_N,\bm{D}_N,\bm{X}_{1,P}} = }\nonumber \\
&&\sum\limits_{\substack{k_c, i_c,q_c,k_d,i_d,q_d,\\ k_{1},i_{1}}}\frac{(\beta^{[2]}(\lambda_{k_c},\lambda_{k_d},\lambda_{k_1}))^{([0,q_d,0])}}{ ([0!,q_d!,0!])} \times \bm{P}_{k_c,i_c}(\bm{N}_{k_c,i_c}^{q_c}-\bm{N}_{k_d,i_d}^{q_d})\bm{N}_{k_d,i_d}^{q_d}\bm{Y}\bm{P}_{k_{1},i_{1}}\nonumber \\
&+&\sum\limits_{\substack{k_c, i_c,q_c,k_d,i_d,q_d,\\ k_{1},i_{1}}}\frac{(\beta^{[2]}(\lambda_{k_c},\lambda_{k_d},\lambda_{k_1}))^{([q_c,0,0])}}{ ([q_c!,0!,0!])} \times \bm{N}_{k_c,i_c}^{q_c}(\bm{N}_{k_c,i_c}^{q_c}-\bm{N}_{k_d,i_d}^{q_d})\bm{P}_{k_d,i_d}\bm{Y}\bm{P}_{k_{1},i_{1}}\nonumber \\
&+&\sum\limits_{\substack{k_c, i_c,k_d,i_d,q_d,\\ k_{1},i_{1}}}\frac{(\beta^{[1]}(\lambda_{k_c},\lambda_{k_1}))^{([q_d,0])}}{ ([q_d!,0!])} \times \bm{P}_{k_c,i_c}\bm{N}_{k_d,i_d}^{q_d}\bm{Y}\bm{P}_{k_{1},i_{1}}\nonumber \\
&-&\sum\limits_{\substack{k_c, i_c,q_c,k_d,i_d,\\ k_{1},i_{1}}}\frac{(\beta^{[1]}(\lambda_{k_d},\lambda_{k_1}))^{([q_c,0])}}{ ([q_c!,0!])} \times \bm{N}_{k_c,i_c}^{q_c}\bm{P}_{k_d,i_d}\bm{Y}\bm{P}_{k_{1},i_{1}},
\end{eqnarray}
to derive the $n$-th derivative of the operator-valued expression
\[
\mathfrak{X}^{\bm{C}_N,\bm{D}_N,\bm{X}_{1,P}}(t),
\]
with respect to \( t \), based on Eq.~\eqref{e3-3:exp: higher order derivative of mathfrak X}, we have to analyze the impact of differentiating terms where \( \bm{X} + t\bm{Y} \) enters implicitly via the spectral data: \( \lambda_{k_d}(t) \), \( \bm{P}_{k_d,i_d}(t) \), and \( \bm{N}_{k_d,i_d}(t) \).

By applying same techniques, Faà di Bruno formula and the chain rule of derivatives, as the evaluation of $\frac{d^n}{dt^n} \mathfrak{X}^{\bm{C}_P,\bm{D}_P,\bm{X}_{1,P}}$, we have the final $n$-th derivative results:
\begin{eqnarray}
\frac{d^n}{dt^n} \mathfrak{X}^{\bm{C}_N,\bm{D}_N,\bm{X}_{1,P}}(t)\Big|_{t=0}
&=& 
\sum_{\substack{k_c, i_c,q_c,k_d,i_d,q_d,\\ k_{1},i_{1}}}
\sum_{\substack{m + j_1 + j_2 = n}} 
\sum_{\pi \in \mathcal{P}(m)} 
\frac{1}{q_d!} \cdot 
\frac{m!}{\prod_{j=1}^m b_j! j!^{b_j}} 
\cdot \frac{1}{j_1! j_2!} \nonumber \\
&&\cdot \left( \partial_{\lambda_{k_d}}^{|\pi|} (\beta^{[2]})^{([0,q_d,0])} \right)(\lambda_{k_c}, \lambda_{k_d}, \lambda_{k_1})
\cdot \prod_{j=1}^m \left( \lambda_{k_d}^{(j)} \right)^{b_j} \nonumber \\
&&\cdot \bm{P}_{k_c,i_c}
\cdot \left( - \bm{N}_{k_d,i_d}^{q_d,(j_1)} \right)
\cdot \bm{N}_{k_d,i_d}^{q_d,(j_2)}
\cdot \bm{Y}
\cdot \bm{P}_{k_1,i_1} \nonumber \\
&+& 
\sum_{\substack{k_c, i_c,q_c,k_d,i_d,q_d,\\ k_{1},i_{1}}}
\sum_{\substack{m + j_1 + j_2 = n}} 
\sum_{\pi \in \mathcal{P}(m)} 
\frac{1}{q_c!} \cdot 
\frac{m!}{\prod_{j=1}^m b_j! j!^{b_j}} 
\cdot \frac{1}{j_1! j_2!} \nonumber \\
&&\cdot \left( \partial_{\lambda_{k_d}}^{|\pi|} (\beta^{[2]})^{([q_c,0,0])} \right)(\lambda_{k_c}, \lambda_{k_d}, \lambda_{k_1})
\cdot \prod_{j=1}^m \left( \lambda_{k_d}^{(j)} \right)^{b_j} \nonumber \\
&&\cdot \bm{N}_{k_c,i_c}^{q_c}
\cdot \left( - \bm{N}_{k_d,i_d}^{q_d,(j_1)} \right)
\cdot \bm{P}_{k_d,i_d}^{(j_2)}
\cdot \bm{Y}
\cdot \bm{P}_{k_1,i_1} \nonumber \\
&+& 
\sum_{\substack{k_c, i_c, k_d, i_d, q_d,\\ k_1, i_1}}
\frac{(\beta^{[1]}(\lambda_{k_c}, \lambda_{k_1}))^{([q_d, 0])}}{q_d!}
\cdot \bm{P}_{k_c,i_c}
\cdot \bm{N}_{k_d,i_d}^{q_d,(n)}
\cdot \bm{Y}
\cdot \bm{P}_{k_1,i_1} \nonumber \\
&-& 
\sum_{\substack{k_c, i_c,q_c,k_d,i_d,\\ k_{1},i_{1}}}
\sum_{\substack{m + j_1 = n}} 
\sum_{\pi \in \mathcal{P}(m)} 
\frac{1}{q_c!} \cdot 
\frac{m!}{\prod_{j=1}^m b_j! j!^{b_j}} 
\cdot \frac{1}{j_1!} \nonumber \\
&&\cdot \left( \partial_{\lambda_{k_d}}^{|\pi|} (\beta^{[1]})^{([q_c,0])} \right)(\lambda_{k_d}, \lambda_{k_1})
\cdot \prod_{j=1}^m \left( \lambda_{k_d}^{(j)} \right)^{b_j} \nonumber \\
&&\cdot \bm{N}_{k_c,i_c}^{q_c}
\cdot \bm{P}_{k_d,i_d}^{(j_1)}
\cdot \bm{Y}
\cdot \bm{P}_{k_1,i_1}.
\end{eqnarray}

\textbf{Evaluation of $\frac{d^n}{dt^n} \mathfrak{X}^{\bm{C}_N,\bm{D}_N,\bm{X}_{1,N}}$}

Since we have
\begin{eqnarray}\label{e3-4:exp: higher order derivative of mathfrak X}
\lefteqn{ \mathfrak{X}^{\bm{C}_N,\bm{D}_N,\bm{X}_{1,N}} = }\nonumber \\
&&\sum\limits_{\substack{k_c, i_c,q_c,k_d,i_d,q_d,\\ k_{1},i_{1},q_1}}\frac{(\beta^{[2]}(\lambda_{k_c},\lambda_{k_d},\lambda_{k_1}))^{([0,q_d,q_1])}}{ ([0!,q_d!,q_1!])} \times \bm{P}_{k_c,i_c}(\bm{N}_{k_c,i_c}^{q_c}-\bm{N}_{k_d,i_d}^{q_d})\bm{N}_{k_d,i_d}^{q_d}\bm{Y}\bm{N}_{k_{1},i_{1}}^{q_1}\nonumber \\
&+&\sum\limits_{\substack{k_c, i_c,q_c,k_d,i_d,q_d,\\ k_{1},i_{1},q_1}}\frac{(\beta^{[2]}(\lambda_{k_c},\lambda_{k_d},\lambda_{k_1}))^{([q_c,0,q_1])}}{ ([q_c!,0!,q_1!])} \times \bm{N}_{k_c,i_c}^{q_c}(\bm{N}_{k_c,i_c}^{q_c}-\bm{N}_{k_d,i_d}^{q_d})\bm{P}_{k_d,i_d}\bm{Y}\bm{N}_{k_{1},i_{1}}^{q_1}\nonumber \\
&+&\sum\limits_{\substack{k_c, i_c,k_d,i_d,q_d,\\ k_{1},i_{1},q_1}}\frac{(\beta^{[1]}(\lambda_{k_c},\lambda_{k_1}))^{([q_d,q_1])}}{ ([q_d!,q_1!])} \times \bm{P}_{k_c,i_c}\bm{N}_{k_d,i_d}^{q_d}\bm{Y}\bm{N}_{k_{1},i_{1}}^{q_1}\nonumber \\
&-&\sum\limits_{\substack{k_c, i_c,q_c,k_d,i_d,\\ k_{1},i_{1},q_1}}\frac{(\beta^{[1]}(\lambda_{k_d},\lambda_{k_1}))^{([q_c,q_1])}}{ ([q_c!,q_1!])} \times \bm{N}_{k_c,i_c}^{q_c}\bm{P}_{k_d,i_d}\bm{Y}\bm{N}_{k_{1},i_{1}}^{q_1},
\end{eqnarray}
to derive the $n$-th derivative of the operator-valued expression, 
\[
\mathfrak{X}^{\bm{C}_N,\bm{D}_N,\bm{X}_{1,N}}(t),
\]
with respect to \( t \) based on Eq.~\eqref{e3-4:exp: higher order derivative of mathfrak X}, we have the following final result:
\begin{eqnarray}
\frac{d^n}{dt^n} \mathfrak{X}^{\bm{C}_N,\bm{D}_N,\bm{X}_{1,N}}(t)\Big|_{t=0}
&=& 
\sum_{\substack{k_c, i_c,q_c,k_d,i_d,q_d,\\ k_{1},i_{1},q_1}}
\sum_{\substack{m + j_1 + j_2 = n}} 
\sum_{\pi \in \mathcal{P}(m)}
\frac{1}{q_d! q_1!} \cdot 
\frac{m!}{\prod_{j=1}^m b_j! j!^{b_j}} \cdot 
\frac{1}{j_1! j_2!} \nonumber \\
&&\cdot 
\left( \partial_{\lambda_{k_d}}^{|\pi|}(\beta^{[2]})^{([0,q_d,q_1])} \right)
(\lambda_{k_c}, \lambda_{k_d}, \lambda_{k_1})
\cdot 
\prod_{j=1}^m \left( \lambda_{k_d}^{(j)} \right)^{b_j} \nonumber \\
&&\cdot \bm{P}_{k_c,i_c}
\cdot \left( -\bm{N}_{k_d,i_d}^{q_d,(j_1)} \right)
\cdot \bm{N}_{k_d,i_d}^{q_d,(j_2)}
\cdot \bm{Y}
\cdot \bm{N}_{k_1,i_1}^{q_1} \nonumber \\
&+& 
\sum_{\substack{k_c, i_c,q_c,k_d,i_d,q_d,\\ k_{1},i_{1},q_1}}
\sum_{\substack{m + j_1 + j_2 = n}} 
\sum_{\pi \in \mathcal{P}(m)}
\frac{1}{q_c! q_1!} \cdot 
\frac{m!}{\prod_{j=1}^m b_j! j!^{b_j}} \cdot 
\frac{1}{j_1! j_2!} \nonumber \\
&&\cdot 
\left( \partial_{\lambda_{k_d}}^{|\pi|} (\beta^{[2]})^{([q_c,0,q_1])} \right)
(\lambda_{k_c}, \lambda_{k_d}, \lambda_{k_1})
\cdot 
\prod_{j=1}^m \left( \lambda_{k_d}^{(j)} \right)^{b_j} \nonumber \\
&&\cdot \bm{N}_{k_c,i_c}^{q_c}
\cdot \left( -\bm{N}_{k_d,i_d}^{q_d,(j_1)} \right)
\cdot \bm{P}_{k_d,i_d}^{(j_2)}
\cdot \bm{Y}
\cdot \bm{N}_{k_1,i_1}^{q_1} \nonumber \\
&+& 
\sum\limits_{\substack{k_c, i_c,k_d,i_d,q_d,\\ k_{1},i_{1},q_1}}
\frac{(\beta^{[1]}(\lambda_{k_c},\lambda_{k_1}))^{([q_d,q_1])}}{q_d! \cdot q_1!}
\cdot \bm{P}_{k_c,i_c}
\cdot \bm{N}_{k_d,i_d}^{q_d,(n)}
\cdot \bm{Y}
\cdot \bm{N}_{k_{1},i_{1}}^{q_1} \nonumber \\
&-& 
\sum_{\substack{k_c, i_c,q_c,k_d,i_d,\\ k_{1},i_{1},q_1}}
\sum_{\substack{m + j_1 = n}} 
\sum_{\pi \in \mathcal{P}(m)}
\frac{1}{q_c! q_1!} \cdot 
\frac{m!}{\prod_{j=1}^m b_j! j!^{b_j}} \cdot 
\frac{1}{j_1!} \nonumber \\
&&\cdot 
\left( \partial_{\lambda_{k_d}}^{|\pi|} (\beta^{[1]})^{([q_c,q_1])} \right)
(\lambda_{k_d}, \lambda_{k_1})
\cdot 
\prod_{j=1}^m \left( \lambda_{k_d}^{(j)} \right)^{b_j} \nonumber \\
&&\cdot \bm{N}_{k_c,i_c}^{q_c}
\cdot \bm{P}_{k_d,i_d}^{(j_1)}
\cdot \bm{Y}
\cdot \bm{N}_{k_1,i_1}^{q_1}.
\end{eqnarray}

\end{example}

\bibliographystyle{IEEETran}
\bibliography{SpecialCase_and_MOI_Bib}

\end{document}